\input amstex
\documentstyle{amsppt}
\magnification=\magstep1
\NoBlackBoxes
\pagewidth{30truecc}
\baselineskip=16truept
\input epsf

\def\A{{\Cal A}}
\def\D{{\Cal D}}
\def\F{{\Cal F}}
\def\G{{\Cal G}}
\def\H{{\Cal H}}
\def\J{{\Cal J}}
\def\K{{\Cal K}}
\def\L{{\Cal L}}

\def\Z{{\Cal Z}}
\def\bn{\text{\bf n}}
\def\cee{{\Bbb C}}
\def\cb{\operatorname{cb}}
\def\op{\operatorname{op}}

\def\MIN{\operatorname{MIN}}
\def\NEW{\operatorname{NEW}}
\def\SEP{\operatorname{SEP}}
\def\OLD{\operatorname{OLD}}
\def\de{\definition}
\def\endde{\enddefinition}
\def\re{\remark}
\def\endre{\endremark}

\def\ep{\varepsilon}
\def\la{\langle}
\def\ra{\rangle}
\def\iitem{\itemitem}
\def\defeq{\overset\text{\rm def}\to=}
\def\chix{{\raise.5ex\hbox{$\chi$}}}
\def\tT{\skew3\tilde{\tilde T}}
\topmatter
\title $M$-Complete approximate identities in operator spaces\endtitle
\author A. Arias and H.P. Rosenthal\endauthor
\address Alvaro Arias, Department of Mathematics, University
of Texas at San Antonio, TX 78249-0664
\endaddress 
\email arias\@math.utsa.edu\endemail
\address Haskell Rosenthal, Department of Mathematics, University
of Texas at Austin, Austin, TX 78712-1082\endaddress
\email rosenthl\@math.utexas.edu\endemail
\abstract
This work introduces the concept of an $M$-complete approximate identity 
(M-cai) for a given operator subspace $X$ of an operator space $Y$. 
M-cai's generalize central approximate identities in ideals in 
$C^*$-algebras, for it is proved that if $X$ admits an M-cai in $Y$, then 
$X$ is a complete $M$-ideal in $Y$. 
It is proved, using ``special'' M-cai's, that if $\J$ is a nuclear ideal in 
a $C^*$-algebra $\A$, then $\J$ is completely complemented in $Y$ for any 
(isomorphically) locally reflexive operator space $Y$ with 
$\J\subset Y\subset\A$ and $Y/\J$ separable. 
(This generalizes the previously known special case where $Y=\A$, due to 
Effros-Haagerup.) 
In turn, this yields a new proof of the Oikhberg-Rosenthal Theorem that 
$\K$ is completely complemented in any separable locally reflexive 
operator superspace, $\K$ the $C^*$-algebra of compact operators on $\ell^2$. 
M-cai's are also used in obtaining some special affirmative answers to the 
open problem of whether $\K$ is Banach-complemented in $\A$ for any 
separable $C^*$-algebra $\A$ with $\K\subset\A\subset B(\ell^2)$. 
It is shown that if conversely $X$ is a complete $M$-ideal in $Y$, then $X$ 
admits an M-cai in $Y$ in the following situations:
(i)~$Y$ has the (Banach) bounded approximation property; 
(ii)~$Y$ is 1-locally reflexive and $X$ is $\lambda$-nuclear for some 
$\lambda \ge1$; 
(iii)~$X$ is a closed 2-sided ideal in an operator algebra $Y$ (via the 
Effros-Ruan result that then $X$ has a contractive algebraic approximate 
identity). 
However it is shown that there exists a separable Banach space $X$ which is an 
$M$-ideal in $Y=X^{**}$, yet $X$ admits no $M$-approximate identity in $Y$.
\endabstract


\endtopmatter

\document

\head Table of Contents\endhead 

{\parindent=40pt\hsize=5.5truein
\line{\indent\hphantom{\S1. }  Introduction \dotfill 2}
\leftline{\indent\S1. Permanance properties of $M$-(Complete)}
\leftline{\indent \hskip.75truein approximate identities
\dotfill 9}
\line{\indent\S2. Complementation results \dotfill 19}
\line{\indent\S3. Examples and complements \dotfill 34}
\line{\indent\hphantom{\S1. } Appendix \dotfill 49}
\line{\indent\hphantom{\S1. } References \dotfill 52}
}
\newpage
\head Introduction\endhead

Let $\K$ denote the $C^*$-algebra of compact operators on separable 
infinite dimensional Hilbert space $H$. 
Consider the following open 

\de{Problem A}
Let $X\subset Y$ be separable operator spaces, and let $T:X\to\K$ be a 
completely bounded (linear) operator. 
Does $T$ admit a bounded linear extension $\tilde T:Y\to \K$? 
That is, can we find a bounded $\tilde T$ completing the following diagram?
$$\matrix Y&&\cr
\raise4ex\hbox{$\bigcup$}&\hbox{\epsfysize=.6truein\epsfbox{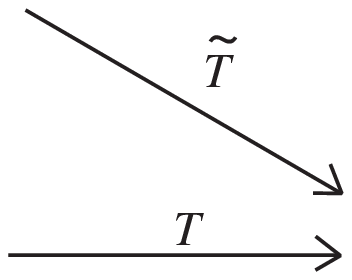}}\cr
\noalign{\vskip-10pt}
X&&\K
\endmatrix
\tag 0.1$$
(See \cite{Pi} or \cite{Ro} for the definition and basic properties of 
operator spaces that we use here.)
An interesting example of E.~Kirchberg yields that one cannot, in general, 
complete this diagram with a completely bounded $\tilde T$ \cite{Ki}; 
by a result of \cite{OR}, it follows there are even locally reflexive 
separable operator spaces where this is the case (in fact, where $Y=C_1$, 
the space of trace class operators). 
\endde

However, the following result is proved in \cite{OR}. 

\proclaim{Theorem 1} 
Assume in $(0.1)$ that $T$ is a complete surjective isomorphism and $Y$ 
is locally reflexive (with $Y$ separable). 
Then there exists a completely bounded $\tilde T$ completing the diagram 
$(0.1)$. 
\endproclaim

We give here a new proof of Theorem 1, using also another structural result 
obtained in \cite{OR}, as well as positive solutions to Problem~A in 
special cases. 
Our methods involve the new concept of an $M$-complete approximate identity 
(an M-cai) for a given operator space $X$ contained in another space $Y$; 
this is a uniformly bounded net $(T_\alpha)$ of (linear) operators from 
$Y$ to $X$ satisfying certain conditions (see Definition~1.1). 
For example, if $X$ is an ideal in a $C^*$-algebra $Y$ and $(x_\alpha)$ 
is a central approximate unit for $X$ in $Y$ consisting of 
positive contractions, then letting $T_\alpha (y) = x_\alpha y$ for all 
$y\in\A$, $(T_\alpha)$ is an M-cai for $X$ in $Y$ (see Proposition~1.4). 

Theorem 1 may be regarded as a ``quantized'' version of a result discovered 
by A.~Sobczyk in (1941): 

\proclaim{Sobczyk's Theorem {\rm \cite{S}}}
Let $X\subset Y$ be separable Banach spaces and $T:X\to c_0$ be 
a given bounded operator. 
There exists a bounded extension $\tilde T:Y\to c_0$ of $T$ with 
$\|\tilde T\| \le 2\|T\|$; moreover ``$2$'' is the best constant here, 
for general $Y$. 
\endproclaim

Many proofs have been given since \cite{S} appeared, cf. \cite{Pe}, \cite{V},
\cite{HWW}, and \cite{Ro}. 
We give yet another proof here, which perhaps explains why this isomorphic 
result (i.e., constant~2) really results from the application of two 
isometric results whose quantized versions 
form the basis for our approach to Problem~A and Theorem~1. 
One of these is a classical theorem concerning  $C(\Omega)$, the space 
of continuous functions on a compact Hausdorff space, 
published in 1933, namely

\proclaim{Borsuk's Theorem {\rm \cite{B}}}
Let $\Omega$ be a compact Hausdorff space and $K$ be 
a closed metrizable subset. 
There exists a norm-one linear operator $L:C(K)\to C(\Omega)$ so that 
$\pi L(f)=f$ for all $f\in C(K)$, where $\pi f= f|K$ for all $f\in C(\Omega)$. 
That is, we have the diagram 
$$\matrix &&C(\Omega)\cr
\noalign{\vskip4pt}
&\hbox{\epsfysize=.5truein\epsfbox{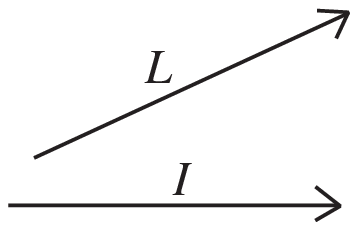}}
&\raise4ex\hbox{\Big\downarrow\rlap{$\scriptstyle \pi$}}\cr
\noalign{\vskip-10pt}
C(K)&&C(K)\ .
\endmatrix
\tag 0.2$$
\endproclaim

Now of course $C(\Omega)$ is a commutative unital $C^*$-algebra; if $K$ is 
a closed subset of $\Omega$ and $\J_K = \{f\in C(\Omega) :f|K=0\}$, then 
$\J_K$ is an ideal in $C(\Omega)$ and every (closed) ideal is of this form.  
Moreover $K$ is metrizable iff $C(K) = C(\Omega)/\J_K$ is separable. 
In view of the Gelfand-Neumark Theorem \cite{GN}, Borsuk's result may 
then  be  reformulated as follows:

\proclaim{Theorem}
Let $\A$ be a unital commutative $C^*$-algebra and  $\J$ a (closed) ideal 
in $\A$ with $\A/\J$ separable. 
Then there exists a contractive (linear) lift $L:\A/\J \to\A$ of 
$I_{\A/\J}$; that is, the following diagram holds. 
$$\matrix &&\A\cr
\noalign{\vskip4pt}
&\hbox{\epsfysize=.5truein\epsfbox{R-li-arrows.eps}}
&\raise4ex\hbox{\Big\downarrow\rlap{$\scriptstyle \pi$}}\cr
\noalign{\vskip-10pt}
\A/\J&&\A/\J\ .
\endmatrix
\tag 0.3$$
\endproclaim

\noindent 
Let us note also that if $L$ is a contractive linear map satisfying (0.3), 
then $\J$ {\it is contractively co-complemented in\/} $\A$ via the map 
$P = I-L\pi$; that is, $P$ is a projection onto $\J$ such that 
$\|I-P\| = \|L\pi\| =1$. 

We now apply Borsuk's  Theorem and the injectivity of $\ell^\infty$, to 
obtain a 

\demo{Proof of Sobczyk's Theorem}
Let $X$, $Y$ and $T$ be as in the statement, and regard 
$c_0\subset\ell^\infty$; note that $\ell^\infty$ is a commutative 
$C^*$-algebra and $c_0$ is an ideal in $\ell^\infty$. 
Now $\ell^\infty$ is isometrically injective (by an easy application of 
the Hahn Banach theorem). 
Hence we may choose a linear extension $\tT :Y\to\ell^\infty$ 
of $T$ with 
$$\|\tT\| = \|T\|\ .
\tag 0.4$$ 
Now let $\A$ be the $C^*$-subalgebra of $\ell^\infty$ generated by $c_0$, 
{\bf 1}, and $\tT(Y)$. 
Then since $c_0\subset \A\subset \ell^\infty$, $c_0$ is an ideal in $\A$ 
and of course $\A$ is commutative; hence by Borsuk's theorem, $c_0$ is 
co-contractively complemented in $\A$. 
Thus we may choose a projection $P$ mapping $\A$ onto $c_0$ with 
$\|I-P\| =1$. 
Thus 
$$\|P\| \le 2\ .
\tag 0.5$$
Then letting $\tilde T = P\tT$, $\tilde T$ is an extension  
of $T$ satisfying 
$$\|\tilde T\| \le \|P\|\, \|\tT\| \le 2\|T\|
\tag 0.6$$
by (0.4) and (0.5).\qed  
\enddemo 

We obtain here the following ``quantized'' version of Borsuk's Theorem
(see Theorem~2.1 and the Theorem of the Appendix).

\proclaim{Theorem 2}
Let $\J$ be a nuclear (2-sided closed) ideal in a $C^*$-algebra $\A$.
Let $\lambda \ge1$ and let $Y$ be a closed linear subspace of $\A$ with 
$\J\subset Y$ so that $Y/\J$ is separable and $Y$ is $\lambda$-locally 
reflexive. 
Then for every $\ep >0$, there exists a completely bounded lift 
$L:Y/\J\to Y$ of $I_{Y/\J}$ with $\|L\|_{\cb} < \lambda +\ep$. 
That is,
$$\matrix &&Y\cr
\noalign{\vskip4pt}
&\hbox{\epsfysize=.5truein\epsfbox{R-li-arrows.eps}}
&\raise4ex\hbox{\Big\downarrow\rlap{$\scriptstyle \pi$}}\cr
\noalign{\vskip-10pt}
Y/\J&&Y/\J
\endmatrix
\tag 0.6$$
holds, where $\pi$ is the quotient map. 
Moreover if $\lambda=1$, $L$ may be chosen to be a complete isometry. 
\endproclaim

E.~Effros and U.~Haagerup establish this result for the case $Y= \A$ in 
\cite{EH}  (when necessarily $\lambda=1$). 
Although our proof uses a basic idea in their discussion, the latter is 
isometric, and does not adapt to the case $\lambda>1$, which however is 
crucial in order to recapture Theorem~1, via the following result obtained 
in \cite{OR}. 

\proclaim{Theorem  3 {\rm (Theorem 1.1 of \cite{OR})}}
Let $Y$ a separable operator space, $X$ a subspace of $Y$, and $T:X\to B(H)$ 
a complete isomorphic injection of $X$ be given. 
There exists a complete isomorphic injection $T':Y\to B(H)$ extending $T$. 
That is, letting $X'= T(X)$ and $Y' = T'(Y)$, $X'\subset Y'\subset B(H)$, 
$T,T'$ are complete isomorphisms so that 
$$\CD
Y @>T'>> Y'\\
\bigcup @.\bigcup\\
X @>T>> X'\ .
\endCD
\tag 0.7$$
\endproclaim

We now obtain Theorem 1 in the same spirit as our proof of Sobczyk's theorem.

\demo{Proof of Theorem 1} 
Let $X\subset Y$ be separable operator spaces with $Y$ locally reflexive 
and let $T:X\to \K$ be a complete surjective isomorphism. 
Now $\K\subset B(H)$ ($H=\ell^2$, say), and $\K$ is an {\it ideal\/} in 
$B(H)$. 
Letting $X'=\K$, choose 
$$Y' \supset \K\ \text{ with }\ Y'\subset B(H)\ \text{ and }\ 
\tilde T:Y\to Y'$$
satisfying the conclusion of Theorem~3. 
Then of course $Y'$ is separable {\it and\/}
locally reflexive (since $\tilde T$ is a complete isomorphism): 
Theorem~2 then yields a completely bounded lift $L:Y'/\K\to Y'$ of 
$I_{Y'/\K}$. 
It follows that setting $P= I_{Y'} - LT'$, $P$ is a completely bounded 
projection of $Y'$ onto $\K$, whence $\tilde T\defeq PT'$ is a completely 
bounded operator satisfying (0.1).\qed
\enddemo

\re{Remarks} 
(a) It is proved in \cite{OR} that $T'$ may be chosen satisfying (0.7) with 
$$\|T'\|_{\cb} \le 3\|T\|_{\cb}\ \text{ and }\ \|T'\|_{\cb}\|T'{}^{-1}\|_{\cb}
\le 12 \|T\|_{\cb}\|T^{-1}\|_{\cb} +6\ .
\tag 0.8$$
The proof of Theorem 1 then yields the existence of absolute positive 
constants $A$ and $B$ (with $A\le 108$ and $B<55$) so that if $Y$ is 
$\lambda$-locally reflexive, then $\tilde T$ may be chosen satisfying (0.1) 
with 
$$\|\tilde T\|_{\cb} \le (A\gamma + B\lambda) \|T\|_{\cb}\ ,
 \ \text{ where }\ \gamma = \|T\|_{\cb} \|T^{-1}\|_{\cb}\ .
 \tag 0.9$$
What are the optimal values of these constants? 
Our estimates (as well as the constants in (0.8)) are surely far from best  
possible. 
We must have however that $A+B\ge2$, even in the case where $\lambda=1$ and 
$T$ is a complete isometry. 
(Actually, it seems likely that in this case, $A+B\ge3$; also that there 
exists such a $T$ so that $Y'$ can not be chosen 1-locally reflexive.) 

(b) N. Ozawa \cite{O} has also (independently) obtained another proof of 
Theorem~1, somewhat along the same lines as our argument.
\endre 

Return now to Problem A, which is easily seen to be a special case of the open 

\de{Problem B} 
Let $\J\subset \A$ be a (closed 2-sided) ideal in a $C^*$-algebra $\A$ 
with $\A/\J$ separable. 
Is $\J$ complemented in $\A$?
\endde

Again, this can be rephrased as a lifting problem, namely, does there exist 
a (bounded linear) lift $L:\A/\J\to \A$ of $I_{\A/\J}$? 
(This problem dates to at least 1974, when it appeared in \cite{A}.) 

To see that Problem A is a special case, consider $\K$ as an ideal in $B(H)$ 
and let $X,Y$ and $T$ be as in the statement of Problem~A. 
Now $B(H)$ is isometrically injective in the operator space category. 
Thus there exists a linear extension $T':Y\to B(H)$ of $T$ satisfying 
$\|T'\|_{\cb} = \|T\|_{\cb}$. 
Now let $\A$ denote the $C^*$-algebra generated by $\K$ and $T'(Y)$. 
Then $\A$ is separable; were $P$ a bounded linear projection from $\A$ 
onto $\K$, $\tilde T\defeq PT$ would be a bounded extension of $\tilde T$ 
satisfying (0.1). 
That is, {\it Problem~A  is equivalent to the special case of 
Problem~B when $\J= \K\subset \A\subset B(H)$.} 

Problem B has an affirmative answer when $\A/\J$ has the (Banach) bounded 
approximation property, by rather deep work of A.~Ando (\cite{An}; see also 
Theorem~2.1 of \cite{HWW} for an exposition). 
Most of the affirmative known results actually yield that $I_{\A/\J}$ 
has a completely positive contractive lift $L$; cf. \cite{CE}, \cite{EH}, 
\cite{ER}. 
(The example in \cite{Ki} does have a contractive lift but no completely 
bounded one.) 
The methods of the present paper recapture Ando's theorem in the special 
setting of Problem~A, via the following result (see Theorem~2.8 and 
Corollary~2.9). 

\proclaim{Theorem 4} 
Let $\A$ be a $C^*$-subalgebra of $B(H)$ with $\K\subset\A$ and assume for some 
$\lambda\ge1$ that $(\K,\A)$ has $\lambda$-extendable local liftings. 
Then for every $\ep>0$ and separable $Y$ with $\K\subset Y\subset \A$, 
there exists a lift $L:Y/\K\to Y$ of $I_{Y/\K}$ with $\|L\| <\lambda+\ep$. 
In particular, this holds if $\A/\K$ has the $\lambda$-bap or $\A$ is 
$\lambda$-extendably locally reflexive. 
\endproclaim

\noindent
(We say that $(\K,\A)$ has $\lambda$-extendable local liftings 
($\lambda$-ell's) if for all $\ep>0$, and finite-dimensional subspaces 
$E$ of $\A/\K$, there exists a linear operator $L:\A/\K\to \A^{**}$  
with $\|L\|<\lambda+\ep$ so that $L(E)\subset \A$ and $L|E$ is a lift 
of the identity injection of $E$ into $\A/\K$.  
See Propositions~2.6 and 2.7 for general permanence properties.) 

\proclaim{Corollary}
If $(\K,B(H))$ has $\lambda$-ell's for some $\lambda\ge1$, Problem~A has an 
affirmative answer.
\endproclaim

\noindent
(The case where $\A$ is extendably locally reflexive in Theorem~4 was 
previously obtained in \cite{OR}. 
The corollary thus extends the consequence obtained there: 
Problem~A has an affirmative answer if $B(H)$ is an extendably locally 
reflexive Banach space.)

We now discuss the methods and setting of our results. 
What is the Banach space technology which yields the theorems of Borsuk, 
Effros-Haagerup and our generalizations thereof?
Why do Banach and operator space hypotheses intervene in the algebraic 
setting of our Theorem~2 and Theorem~4, and what is the appropriate operator 
space setting of these results? 
The answer to the first question lies in the concept of an $M$-ideal, as 
pioneered in \cite{AE}; see \cite{HWW} for a comprehensive reference. 

We briefly recall the relevant notions. 

\de{Definition 0.1} 
Let $X\subset Y$ be Banach/operator spaces. 

(a) $X$ is called an $M$-summand in $Y$ if there exists a closed linear 
subspace $Z$ of $Y$ with 
$$X\oplus Z =Y
\tag 0.10i$$
so that 
$$\|x+z\| = \max \{\|x\|,\|z\|\}\ \text{ for all $x\in X$ and $z\in Z$.}
\tag 0.10ii$$
In case these are operator spaces, $X$ is called a complete $M$-summand 
if $Z$ satisfying (0.10i) also satisfies 
$$\split
&\|(x_{ij} + z_{ij})\| = \max\{ \|(x_{ij})\|,\|(z_{ij})\|\} 
\text{ for all $n$ and $n\times n$-matrices}\\
&(x_{ij})\text{ and } (z_{ij})\text{ of elements of $X$ and $Z$ respectively.}
\endsplit\tag 0.10$'$ii$$

(b) $X$ is called an $L$-summand if $Z$ can be chosen satisfying (0.10i) 
so that 
$$\|x+z\| = \|x\|+ \|z\|\ \text{ for all }\ x\in X \text{ and }z\in Z\ .
\tag 0.11$$

(c) $X$ is called an $M$-ideal (resp. complete $M$-ideal) in $Y$ if 
$X^{**}= X^{\bot\bot}$ is an $M$-summand 
(resp. complete $M$-summand) in $Y^{**}$. 
\endde

It turns out that $M$ (resp. $L$) summands for $X$ are unique; if 
$X\oplus Z$ is the corresponding $M$ (resp. $L$) decomposition, the 
projection $P$ from $Y$ onto $X$ with kernel $Z$ is called the $M$ (resp. 
$L$) projection onto $X$. 
Also, $X$ is an $M$-ideal in $Y$ if and only if $X^\bot$ is an $L$-summand 
in $Y^*$ (see \cite{HWW}; also see \cite{ER2} for the case of complete 
$M$-ideals). 

Now $M$-summands and $M$-ideals are very special in the general Banach 
space setting. 
However the following remarkable result shows their importance. 

\proclaim{Theorem} 
Let $\A$ be a $C^*$-algebra and $\J$ be a closed linear subspace. 
Then $\J$ is an $M$-ideal in $\A$ iff $\J$ is a (2-sided) ideal in $\A$ 
iff $\J$ is a complete $M$-ideal in $\A$. 
\endproclaim

\noindent 
(See \cite{HWW} for a proof and complete references; for the remarkable 
theorem  that $M$-ideals are algebraic ideals, see \cite{AE} and \cite{SW}.) 

Why is an ideal in a $C^*$-algebra an $M$-ideal? 
The commutative case is rather transparent. 
For then $\A = C_0(\Omega)$ for some locally compact Hausdorff space 
$\Omega$ and $\J = \J_K$ for some closed subset $K$ of $\Omega$ (and 
Borsuk's theorem of course could be formulated in this possibly 
non-unital setting also). 
But then identifying $\A^*$ with $M(\Omega)$, the space of regular complex 
Borel measures ``on'' $\Omega$, $\J^\bot = M(K)$ and then letting 
$Z= M(\Omega\sim K)$, $\J^\bot \oplus Z$ is an $L$-decomposition of $M(K)$.

The non-commutative case is certainly not so transparent. 
The highly motivated case of $\K$ in $B(H)$ was established by Diximer in 
1950 \cite{Di}. 
Note that it seems one must at least consider $\K^\bot \subset B(H)^*$, 
a rather huge object! 
Our approach here yields the general non-commutative case, via 
$M$-approximate identities, a notion defined only, in this setting, on 
the $C^*$-algebra $\A$ itself; one has no need to ``look'' at $\A^*$ or 
$\A^{**}$, to then ``see'' the $M$-ideal property. 

The results of our paper are all cast in the general setting of 
Banach/operator spaces and (complete) $M$-ideals. 
We shall see that our complementation results also use the property 
of certain ``special'' $M$-complete approximable identities, and not 
just the general M-cai concept. 
Our methods may also be used to recapture several of the results given 
in the initial paper \cite{Ro}. 

To more thoroughly answer the second of our ``motivating'' questions, we 
now proceed with a more detailed discussion and outline of our results. 

Various refinements of the concept of an $M$-approximate identity are 
given in Definition~1.1. 
Theorem~1.1 then establishes that if $X$ admits an M-ai (resp. M-cai) in 
$Y$, $X$ is an $M$-ideal (resp. complete $M$-ideal) in $Y$. 
Moreover if $X$ admits a strong M-cai $(T_\alpha)$, then $T_\alpha^{**}$ 
converges in the weak*-operator topology (W*-OT) on $Y^{**}$ to the 
$M$-projection mapping $Y^{**}$ onto $X^{**}$. 
We show in Proposition~1.4 that central approximate units yield strong 
contractive M-cai's for ideals $\J$ in $C^*$-algebras $\A$. 
A by-product  of Corollary~1.5: the central approximate unit $(u_\alpha)$ 
may be chosen so that setting $U_\alpha a = u_\alpha a$ for all $a\in\A$, 
then 
$$\lim_\alpha (\|U^*_\alpha y^*\| + \|(I-U_\alpha^*) y^*\|) 
= \|y^*\|\ \text{ for all }\ y^*\in \A^*\ .
\tag 0.12$$

Section one concludes with a permanence property of ``good'' M-cai's 
(see Definition~1.3), which has the consequence: 
if $X$ admits a good M-cai in $Y$, then $Z\oplus_{\op}X$ is a complete 
$M$-ideal in $Z\otimes_{\op}Y$ for all operator spaces $Z$ (Proposition~1.7). 
It remains an open problem if this permanence property holds for general 
complete $M$-ideals. 

In \S2, we introduce the notion of a special M-cai for an operator space 
$X\subset Y$. 
We show in Proposition~2.3 that if $X$ is an ideal in a $C^*$-algebra $Y$ 
and $(x_\alpha)$ is a positive contractive central approximate unit 
for $X$ in $Y$ (the $x_\alpha$'s being positive contractions in $X$), then 
defining $U_\alpha (y) = \sqrt{x_\alpha} y\sqrt{x_\alpha}$ for all 
$\alpha$, $(U_\alpha)$ is a special M-cai for $X$ in $Y$. 
Thus Theorem~2 of this Introduction is a special case (for $\lambda>1$) 
of our Theorem~2.4: 
{\it If $X$ is an approximately injective subspace of a $\lambda$-locally 
reflexive subspace $Y$ so that $X$ admits a special M-cai in $Y$, with 
$Y/X$ separable, then for all $\ep>0$, there is a lift $L:Y/X\to Y$ of 
$I_{Y/X}$ with $\|L\|_{\cb} <\lambda+\ep$.}

In the Theorem of the Appendix, we obtain a complete isometric extension 
of the Effros-Haagerup lifting result, without the special M-cai assumption; 
namely, the lift $L$ may be chosen completely contractive provided $X$ 
is an approximately injective complete $M$-ideal in $Y$, when $Y$ is 
1-locally reflexive and $Y/X$ is separable. 

Our proof of Theorem 2.4 yields that its conclusion holds if we replace its 
hypothesis that $Y$ is locally reflexive by the assumption that $X$ is 
{\it locally complemented\/} in $Y$; that is, for some $\gamma\ge1$, 
$X$ is $\gamma$-completely complemented in $Z$ for all $X\subset Z\subset Y$ 
with $Z/X$ finite dimensional. 
(See Theorem~$2.4'$ and and the following Remark.) 
Thus it follows that 
{\it if $\K\subset Y\subset B(\ell^2)$ with $Y$ separable, then $\K$ is 
completely complemented in $Y$ provided $\K$ is locally complemented in $Y$.} 

In Definition 2.2, we give the Banach space concept of 
{\it extendable local liftings\/} (ell's) for a pair of Banach spaces 
$X\subset Y$; (this is the same as the special case $(\K,\A)$ stated above). 
We observe in Proposition~2.7 that $(X,Y)$ has $\lambda$-ell's 
if $(X,Y)$ has $\lambda$-local liftings and  
e.g., $Y/X$ has the $\lambda$-bap or $Y$ is $\lambda$-extendably 
locally reflexive. 
Theorem~4 of this Introduction is then a special case of our 
Theorem~2.8, which yields that 
{\it if $X\subset Y$ are operator spaces with $X$ approximately  injective 
and $Y/X$ separable, so that $X$ admits a special M-cai in $Y$ and $(X,Y)$ 
has extendable local liftings, then $X$ is complemented in $Y$.} 

Section 3 gives further applications of the general complementation results 
in Section~2. 
The easy Proposition~3.2(a) yields that if $X_1,X_2,\ldots$ are given 
operator spaces and $X = (X_1\oplus X_2\oplus \cdots)_{c_0}$  and 
$Y = (X_1\oplus X_2 \oplus\cdots)_{\ell^\infty}$, then $X$ admits a 
(canconical) strong special M-cai in $Y$; 
moreover if the $X_j$'s are approximately injective, so is $X$. 
Corollary~3.3 then yields that if $X\subset Z\subset Y$ with $Z/X$ 
separable and $Z$ locally reflexive, $X$ is completely complemented in $Z$. 
Moreover, if $Z$ is 1-locally reflexive, $X$ is completely co-contractively 
complemented in $Z$. 
This yields the discoveries in \cite{Ro} that if the $X_j$'s are all 
1-injective Banach spaces, $X$ has the 2-Separable Extension Property; in 
particular, $c_0(\ell^\infty)$ has the  2-SEP. 
If moreover, the $X_j$'s are all 1-injective operator spaces and 
$X\subset Z\subset Y$ with $Z/X$ separable and $Z$ $\lambda$-locally 
reflexive, then $X$ is completely $(\lambda+\ep)$-cocomplemented in $Z$ 
for all $\ep>0$. 

We also recapture the main result obtained in \cite{Ro} concerning the 
Complete Separable Extension Property (CSEP), namely that for all $n\ge1$, 
$Z\defeq c_0(M_{n,\infty}\oplus M_{\infty,n})$ has the 2-CSEP. 
That is, for all separable operator spaces $X\subset Y$ and completely 
bounded maps $T:X\to Z$, there is an extension $\tilde T:Y\to Z$ with 
$\|\tilde T\|_{\cb} \le 2\|T\|_{\cb}$. 
(This is a full quantized extension of Sobczyk's Theorem 
(see Corollary~3.2). 
In turn, this is obtained via an interesting recent operator space 
extension of the Banach local reflexivity principle due to L.~Ge and 
D.~Hadwin \cite{GH} and the following application (via an elementary 
result in \cite{Ro}): 
For $Z$ as above, $Z^{**}$ is 1-locally reflexive 
(Proposition~3.7). 

Section 3 also treats the question of when the converse of Theorem~1.1 holds. 
Precisely, {\it suppose $X\subset Y$ are Banach (resp. operator) spaces 
with $X$ an $M$-ideal (resp. complete $M$-ideal) in $Y$. 
Does $X$ admit an M-ai (resp M-cai) in $Y$?}

In Theorem 3.1, we prove that this is indeed true in the case where $X$ is 
an ideal (closed 2-sided) in a (possibly) non-self-adjoint operator 
algebra $Y$. 
Effros-Ruan prove in \cite{ER1} that then $X$ is an $M$-ideal in $Y$ 
precisely when $X$ has a contractive approximate identity. 
We show directly that then $X$ admits a strong contractive M-cai 
in $Y$.

In Theorem 3.11, we obtain the (perhaps surprisingly general) result that 
(assuming $X$ is a complete $M$-ideal in $Y$), 
{\it $X$ admits an M-cai in $Y$ provided $Y$ has the Banach bounded 
approximation  property. 
We also obtain the same conclusion if $Y$ is 1-locally reflexive and 
$X$ is a finitely injective operator space.} 
Moreover we obtain that the M-cai $(T_\alpha)$ may be chosen to consist 
of finite rank operators when $Y$ has the bap or $X $is $\lambda$-nuclear 
and $Y$ is 1-locally reflexive. 
($\lambda$-finitely injective operator spaces are defined in 
Definition~3.2 (just before the statement of Theorem~3.11); these include 
$\lambda$-nuclear and $\lambda$-injective operator spaces). 

Theorem 3.11 uses an extension of the Banach local reflexivity 
principle due to S.~Bellenot \cite{Be}, 
which we formulate and prove in the operator space setting in 
Lemma~3.12, as well as its consequence, an extension of the above mentioned 
result of \cite{GH}, which we obtain in Lemma~3.13. 

Finally, we give an example of a Banach space $X$ which is an $M$-ideal 
in $X^{**}$, yet $X$ admits no $M$-approximate identity in $X^{**}$ 
(Proposition~3.16). 
The example is at the ``surface'' modulo some rather deep known results; 
namely, $X$ is a subspace of $c_0$ failing the compact bounded 
approximation property. 

We do not know of a separable pair $(X,Y)$ forming a counterexample
(as of this writing!). 
However we conjecture that if $X$ is as in Proposition~3.16, then there 
exists $X\subset Y\subset X^{**}$ with $Y$ separable, yet $X$ admits no 
M-ai in $Y$. 
This conjecture, however, appears to lie considerably below the 
surface of known results, unlike 3.16. 

Initial stages of this work were carried out by the second named author 
during a visit to the Mathematics Laboratory of the Faculty of Sciences 
of Marseilles at Saint-Jerome. 
It is his pleasure to thank 
the members of the Mathematics Equipe and especially the 
Operator Spaces Groupe de Travaille  at Saint-Jerome for the warm hospitality 
and mathematical encouragement shown him, with particular thanks to 
Christian Samuel. 
Further stages of this work were accomplished during the 1999 Summer 
Analysis Workshop at Texas A\&M University. 
It is the pleasure of both authors to thank the Workshop participants 
and organizers for their warm support during this visit.

\head \S1. $M$-(Complete) approximate identities in (complete) $M$-ideals
\endhead

We begin with the definition of the basic concept introduced in this work. 

\de{Definition 1.1} 
A. Let $X\subset Y$ be Banach/operator spaces, and let 
$(T_\alpha)_{\alpha\in\D}$ be a uniformly bounded net in $B(Y)$. 
$(T_\alpha)_{\alpha\in\D}$ is an {\it $M$-approxi-mate identity\/} (M-ai)
for $X$ in $Y$ if 
\roster
\item"(i)" $T_\alpha Y\subset X$ for all $\alpha$. 
\item"(ii)" $T_\alpha x\to x$ for all $x\in X$.
\item"(iii)" $\varlimsup_\alpha \|T_\alpha u + (I-T_\alpha)v\| \le 
\max \{\|u\|\, \|v\|\}$ for all $u,v\in Y$.
\endroster
$(T_\alpha)_{\alpha\in\D}$ is an {\it $M$-complete approximate identity\/}
(M-cai) if $(T_\alpha)$ satisfies  (i), (ii), and 
\roster
\item"(iii)$'$" for all $n$ and all $n\times n$ matrices 
$(u_{ij})$, $(v_{ij})$  in $Y$, 
$$\varlimsup_\alpha \|(T_\alpha u_{ij} + (I-T_\alpha) v_{ij})\| 
\le \max \{ \|(u_{ij})\|\, \|(v_{ij})\|\} \ .$$
\endroster

B. $(T_\alpha)$ is called a {\it strong\/} M-ai (resp. strong M-cai) if 
in addition we have 
\roster
\item"(iv)" $T_\alpha^{**} x^{**} \buildrel {w^*}\over \to 
 x^{**}$ for all $x^{**}\in X^{**}$
\endroster
and
\roster
\item"(v)" $\varlimsup_\alpha \|T_\alpha^* y^*\|\le \|y^*\|$ for all 
$y^* \in Y^*$
\endroster
resp.
\roster
\item"(v$'$)" The $T_\alpha$'s are completely bounded and 
$$\varlimsup_\alpha \|I_{\K} \otimes T_\alpha^* (\tau) \|\le \|\tau\|
\text{ for all }
\tau \in \K\otimes_{\op} Y^*\ .$$ 
\endroster

C. $(T_\alpha)$ is called a {\it contractive\/} M-ai (resp. contractive 
M-cai) if $(T_\alpha)$ is an M-ai (resp. M-cai) for $X$ in $Y$ so that 
$\|T_\alpha\| \le1$ (resp. $\|T_\alpha\|_{\cb} \le 1$) for all $\alpha$. 
\endde

\re{Remarks} 
1. Evidently if $(T_\alpha)$ is a contractive M-ai (resp. a contractive M-cai)
satisfying (iv), then $(T_\alpha)$ is a strong M-ai (resp. strong M-cai) 
for $X$ in $Y$. 

2. All these concepts are hereditary in the following sense: 
if $X\subset Z\subset Y$ and $(T_\alpha)$ is an $M$-approximate identity 
of one of the various kinds, for $X$ in $Y$, then also $(T_\alpha|Z)$ 
is an M-ai of the {\it same kind\/} for $X$ in $Z$.

3. We do not know the answers to the following questions. 
Let $X\subset Y$ be Banach/operator spaces. 
If $X$ admits an M-ai in $Y$, does $X$ admit a contractive M-ai in $Y$? 
If $X$ admits an M-cai in $Y$, does $X$ admit an M-cai 
$(T_\alpha)$ with 
\roster
\item"(i)" the $T_\alpha$'s completely bounded? 
\item"(ii)" $\sup_\alpha \|T_\alpha\|_{\cb} <\infty$? 
\item"(iii)" $\|T_\alpha\|_{\cb}\le 1$ for all $\alpha$ (i.e., so that 
$(T_\alpha)$ is a contractive M-cai)?
\endroster
\endre

Our first result provides basic motivation for introducing the concept
in Definition~1.1. 

\proclaim{Theorem 1.1} 
Let $X\subset Y$ be Banach (resp. operator) spaces and assume $X$ admits 
an M-ai (resp.  M-cai) $(T_\alpha)$ in $Y$. 
Then $X$ is an $M$-ideal (resp. complete $M$-ideal) in $Y$ and in fact 
$(T_\alpha^*)$ converges in the $W^*$-OT to the $L$-projection $P$ on $Y^*$ 
with kernel $X^\bot$.

If $(T_\alpha)$ is a strong M-ai (resp. M-cai), then $(T_\alpha^{**})$ 
converges in the $W^*$-OT on $Y^{**}$ to the $M$-projection $P^*$ onto 
$X^{**}$ (resp. $I_{\K} \otimes T_\alpha^{**}$ converges in the 
$W^*$-OT on $\K\otimes_{\op} Y^{**}$ to $I_{\K}\otimes P^*$). 
\endproclaim

\re{Remark}
We note that in the operator space case, it is {\it not\/} assumed in the 
definition of an M-cai $(T_\alpha)$ that the $T_\alpha$'s are completely 
bounded. 
If, however, we assume that in fact $\sup_\alpha \|T_\alpha\|_{\cb} <\infty$ 
and $(T_\alpha)$ is an M-cai,
we obtain that $(I_{\K}\otimes T_\alpha^*)$ converges in the $W^*$-OT 
on $\K \otimes_{\op} Y^*$ to $I_{\K}\otimes P$. 
\endre

\demo{Proof of Theorem 1.1} 
We first deduce: There is an $L$-decomposition 
$X^\bot \oplus W$ for $Y^*$, and if $P$ denotes the projection  onto $W$ 
with $\ker X^\bot$, then $T_\alpha^* \to P$ in the  $W^*$-OP. 
By first passing to a  sub-net, we may assume: 
$$w^*\!\!\!\!-\!\!\lim_\alpha T_\alpha^* y^* \text{ {\it exists}  for all }y^* \in Y^* \ .
\tag 1.1$$ 
(Later, we will show this is not needed.) 
Now we have 
$$T_\alpha^* x^\bot=0 \text{ and hence }(I-T_\alpha^*) x^\bot =x^\bot\ 
\text{ for all }\alpha, \text{and all }x^*\in X^*\ .
\tag 1.2$$ 
Trivial: because $\la T_\alpha^* x^\bot ,y\ra = \la x^\bot,
T_\alpha y\ra =0$ for all $y\in Y$. 
Thus 
$$\lim_\alpha (I-T_\alpha^*) y^* \in X^\bot\ , \text{ for all }y^* \in Y^*\ . 
\tag 1.3$$
Indeed, if $x\in X$, then 
$$\lim_\alpha \la (I-T_\alpha^*)y^* ,x\ra 
= \lim_\alpha \la y^* ,(I-T_\alpha)x\ra =0\ .$$
Since the $T_\alpha^*$'s are {\it uniformly bounded\/}, the operator 
$$Q = \lim_{\alpha} (I-T_\alpha^*)\ \text{ is bounded}
\tag 1.4$$
(where the net converges in the $W^*$-OT).
But we have: $y^*\in Y^* \implies Qy^* \in X^\bot$ by (1.3) and 
$y^*\in X^\bot \implies Qy^* = y^*$ by (1.2), hence $Q$ is indeed 
a projection onto 
$X^\bot$. 
Let $P = I-Q$ and $W=PX^* = \ker Q$.

Now we prove: 
$X^\bot \oplus W$ is an $L$-decomposition of $Y^*$. 
Let $x^\bot \in X^\bot$, $w\in W$; $\ep >0$. 
We may choose $u,v$ norm-1 elements of $Y$ so that 
$$\|x^\bot\| + \|w\| \le (1+\ep)(\la x^\bot,u\ra +\la w,v\ra)
\tag 1.5$$
(and the right hand terms are actually non-negative real numbers).
Now 
$$\la x^\bot,u\ra + \la w,v\ra 
= \la x^\bot,u\ra + \lim_\alpha \la T_\alpha^* w,v\ra\ .
\tag 1.6$$
For all $\alpha$, 
$$\align 
&\la x^\bot +w,(I-T_\alpha)u+T_\alpha v\ra\tag 1.7\\
&\qquad  = \la x^\bot ,u\ra + \la w,(I-T_\alpha)u\ra 
+ \la w,T_\alpha v\ra\ \text{ by (1.2)}\\
&\qquad  = \la x^\bot,u\ra + \la (I-T_\alpha^*)w,u\ra 
+ \la T_\alpha^* w,v\ra\ .
\endalign$$
But $\lim_\alpha (I-T_\alpha)^* w = Qw=0$, whence we have proved, 
by (1.6) and (1.7), that 
$$\align 
\la x^\bot ,u\ra + \la w,v\ra 
&= \lim_\alpha \la x^\bot +w,(I-T_\alpha)u+T_\alpha (v)\ra \tag 1.8\\
&\le \|x^\bot +w\| \lim_\alpha \|(I-T_\alpha)u+T_\alpha(y)\|\\
&\le \|x^\bot +w\|\ \text{(by Definition 1.1iii).}
\endalign$$
Hence by (1.5), $\|x^\bot \| + \|w\| \le (1+\ep)(\|x^\bot +w\|$. 
Since $\ep>0$ is arbitrary; this shows: $P$ is indeed an $L$-projection.

Thus by uniqueness of such, we conclude by the way, that we did not need to 
take a subnet; and so our original net satisfies:
$$P = \lim  T_\alpha^*\ \text{ in the $W^*$-OT}\ .
\tag 1.9$$
We have now proved: $X$ is indeed an $M$-ideal with $X^{**} \oplus W^\bot  
= Y^{**}$ the $M$-decomposition for $Y^{**}$. 

It now follows immediately from results of Effros-Ruan that $X$ is 
a complete $M$-ideal in $Y$ in case $(T_\alpha)$ is an M-cai. 
Indeed, fixing $n$, then defining $T_\alpha^n = I_{M_n} \otimes T_\alpha$ 
on $M_n (Y)$, we obtain that $(T_\alpha^n)$ is an M-ai for $M_n(X)$ 
in $M_n(Y)$, hence $M_n(X)$ is an $M$-ideal in $M_n(X)$, whence $X$ 
is indeed a complete $M$-ideal in $Y$ by \cite{ER2}.

We now proceed to the final assertion of the Theorem. 
Assume then that $(T_\alpha)$ is a strong M-ai (resp. M-cai). 

Let $w^\bot\in W^\bot$. 
We claim 
$$T_\alpha^{**} w^\bot \to0 \quad \text{ in the } w^* \ \text{topology}.
\tag 1.10$$ 
Now once (1.10) is proved, we have by (iv) that for all $x^{**}\in X^{**}$ 
and $w^\bot\in W^\bot$, 
$$w^*\!\!\!\!-\!\!\lim_\alpha T_\alpha^{**} (x^{**} + w^\bot)= x^{**}, 
\tag 1.11$$ 
whence $T_\alpha^{**}\to P^*$ in the $W^*-OT$, and of course $P^*$ is the 
$M$-projection onto $X^{**}$.

Suppose (1.10) were false. 
Then for some $y^* \in Y^*$, 
$$\varlimsup_\alpha |\la y^*,T_\alpha^{**}w^\bot\ra |\ne0\ .
\tag 1.12$$
Now let $y^* = x^\bot +w$, $x^\bot \in X$, $w\in W$. 
But then trivially (since $T_\alpha^* x^\bot=0$ for all $\alpha$), 
$$\varlimsup_\alpha |\la w,T_\alpha^{**}w^\bot\ra |\ne0\ .
\tag 1.13$$
Now (by passing to a subnet if necessary and taking obvious normalizations), 
we may assume without loss of generality that $\|w\|=1=\|w^\bot\|$ 
and there is a $\delta >0$ so that 
$$|\la w,T_\alpha^{**}w^\bot\ra | \ge \delta\ \text{ for all }\alpha\ .
\tag 1.14$$
Now choose $x\in X$, $\|x\| =1$, so that 
$$|\la w,x\ra | > 1-\frac{\delta}2\ .
\tag 1.15$$
Choose $\alpha_0$ so that also 
$$|\la T_\alpha^* w,x\ra| > 1-\frac{\delta}2\ \text{ for all } \alpha 
>\alpha_0\ .
\tag 1.16$$
Thus combining (1.14) and (1.16), we find for such $\alpha$ that 
$$|\la T_\alpha^* w,w^\bot\ra| + |\la T_\alpha^* w,x\ra| > 1+\frac{\delta}2\ .
\tag 1.17$$ 
Finally, for each such $\alpha$, choose scalars 
$\theta_\alpha$ and $\psi_\alpha$ 
of modulus one so that 
$$| \la T_\alpha^* w,w^\bot\ra | 
= \la T_\alpha^* w,\theta_\alpha w^\bot \ra
\text{ and }\la T_\alpha^* w,x\ra = \la T_\alpha^* w,\psi_\alpha x\ra\ .
\tag 1.18$$
Hence we have by (1.17) and (1.18) that 
$$\align
1+\frac{\delta}2 & < |\la T_\alpha^* w,\theta_\alpha w^\bot + \psi_\alpha x
\ra |\tag 1.19\\
&\le \|T_\alpha^* w\|\ \|\theta_\alpha w^\bot + \psi_\alpha x\|\\
& = \|T_\alpha^* w\|
\endalign$$
because $x\in X^{**}$ and $X^{**} \oplus W^\bot$ is an $M$-decomposition.

Thus finally
$$\varlimsup_\alpha \|T_\alpha^* w\| \ge 1+\frac{\delta}2\ \text{ but } 
\|w\|=1\ ,
\tag 1.20$$ 
contradicting (v).

Finally, in the complete-category, it follows that for all $\mu \in \K\otimes 
Y^*$, 
$$(I_{\K} \otimes T_\alpha^*) \mu \to (I_{\K} \otimes P)(\mu) 
\text{ weakly.}
\tag 1.21$$
But $\K\otimes Y^*$ is a norm-dense linear sequence of $\K\otimes_{\op} Y^*$, 
whence (v$'$) then indeed yields that (1.21) holds for all 
$\mu\in\K\otimes_{\op} Y^*$, yielding the final assertion.\qed
\enddemo

\proclaim{Corollary 1.2}
Suppose that $X$ admits a strong M-ai (resp. a strong M-cai) in $Y$. 
Then also $X$ admits a strong M-ai (resp. M-cai) $(U_\alpha)$ such that if $P$ 
is the $L$-projection with kernel $X^\bot$, then $U_\alpha^* \to P$ in 
the SOT on $Y^*$. 
In particular, 
$$\varlimsup_\alpha \|U_\alpha^* y^*\| + \|(I-U_\alpha^*) y^*\| 
= \|y^*\|\ \text{ for all } y^* \in Y^*\ .
\tag 1.22$$ 
If $X$ admits a strong M-cai, then $(U_\alpha)$ may be chosen so that 
$I_{\K}\otimes U_\alpha^* \to I_{\K} \otimes P^*$ in the SOT on 
$\K\otimes_{\op} Y^*$.
\endproclaim

\demo{Proof} 
Let $(T_\alpha)$ be a strong M-ai (resp. M-cai) for $X$ in $Y$. 
Theorem 1.1 shows that 
$$T_\alpha^* y^* \to Py^* \ \text{ weakly, for all } \ y^*\in Y^*\ .
\tag 1.23$$ 
Hence there exists a net $(U_\alpha)$ of ``far-out'' convex combinations
of the $T_\alpha$'s, so that 
$$U_\alpha^* y^* \to Py^* \ \text{ in norm, for all } y^*\in Y^*\ .
\tag 1.24$$
Of course this yields that $U_\alpha^* \to P$ SOT. 
In particular, for $y^*\in Y^*$, we have since $P$ is an $L$-projection,  
that for all $y^* \in Y^*$, 
$$\align 
\|y^*\| & = \|Py^*\| + \|(I-P)y^*\| \tag 1.25\\
& = \lim_\alpha \Big( \|U_\alpha^* y^* \| + \|(I- U_\alpha^* ) y^*\|\Big)\ ,
\endalign$$ 
proving (1.22).

Now it is easily seen that  $(U_\alpha)$ is a strong M-ai (resp. M-cai).
Finally, if $(T_\alpha)$ is a strong M-cai, then 
$$(I_{\K}\otimes  T_\alpha^*)  \to (I_\kappa \otimes P)\ \text{ in the 
WOT on } \K\otimes_{\op} Y^*\ ,
\tag 1.26$$
whence again the $U_\alpha$'s may be chosen as above to satisfy the final 
assertion of the Corollary.\qed
\enddemo

Let $X\subset Y$ be Banach (resp. operator) spaces. 
An inspection of the proof of Theorem~1.1 shows that its conclusion holds under
the following modified assumptions on the net of operators 
$(T_\alpha)$. 

\de{Definition 1.2} 
Let $(T_\alpha)$ be a net of contractions (resp. complete contractions) on $Y$. 
$(T_\alpha)$ is called a {\it weak contractive M-ai\/} 
({\it resp. weak contractive M-cai\/})
if (i) and (ii) of Definition~1.1 hold, but in (iii) and (iii$'$), we restrict 
the $u$'s (resp. the $u_{ij}$'s) to lie in $X$. 
\endde

\proclaim{Theorem 1.1$'$} 
The conclusion of Theorem 1.1 holds provided $X$ admits a weak contractive 
M-ai (resp. a weak contractive M-cai) in $Y$. 
\endproclaim

\demo{Proof} 
We may  argue as in the proof of Theorem 1.1, which essentially proceeds 
from first principles. 
Alternatively, we may use the following characterization of $M$-ideals 
given in Theorem~2.2 of \cite{HWW} (the ``restricted 3-ball property''): 
{\it $X$ is an $M$-ideal in $Y$ provided for any $y$, $x_1,x_2,x_3$ in 
$Ba\, Y$ and $\ep>0$, there exists an $x\in Ba\,X$ with} 
$$\|x_i + y-x\| \le 1+\ep \ \text{ for all }\ 1\le i\le 3\ .$$
Now letting $(T_\alpha)$ be a weak contractive M-ai for $X$ in $Y$, 
simply choose $\alpha$ so that 
$$\|x_i + (I-T_\alpha)(y)\| \le 1+\ep\ \text{ for all }\ 1\le i\le 3\ .$$
{\it Since\/} the $T_\alpha$'s are contractions and map $Y$ into $X$, 
$x\defeq T_\alpha y$ satisfies the above criterion, whence $X$ is an 
$M$-ideal in $Y$. 
It may then be directly verified that $T_\alpha^* \to P$ in the $W^*$-OT, 
where $P$ is the $L$-projection with kernel $X^\bot$. 
Finally, assuming that $(T_\alpha)$ is a weak contractive M-cai, we 
have that for all $n$, $(I_n\otimes T_\alpha)$ is a weak contractive 
M-ai on $M_n\otimes Y$ 
($I_n$ the identity on $M_n$), whence  $M_n(X)$ is an $M$-ideal in $M_n(Y)$ 
and so as before, via the result in \cite{ER2}, $X$ is a complete 
$M$-ideal in $Y$. 
Finally, assuming the additional hypotheses that 
$T_\alpha^{**} x^*\to x^{**}\ w^*$ for all $x^{**}\in X^{**}$, we 
obtain the final conclusion of Theorem~1.1.\qed
\enddemo

We have chosen the stronger concept given in Definition~1.1 
in the contractive case, since this is 
what occurs in the case where $X$ is an ideal in a $C^*$-algebra $Y$. 
Before dealing with this remarkable special case, we consider briefly 
when our methods yield that a Banach/operator space $X$ is an 
$M$-ideal/complete $M$-ideal in $X^{**}$. 

\proclaim{Corollary 1.3}
Let $X$ be a Banach (resp. operator) space and let $(T_\alpha)$ be a net 
of weakly compact operators on $X$. 
Suppose either of the following  two hypotheses: 
\roster
\item"(I)" $(T_\alpha^{**})$ is an M-ai (resp. M-cai) for $X$ in 
$X^{**}$. 
\item"(II)" 
The $T_\alpha$'s are contractions so that both 
\vskip1pt
\iitem{\rm (i)} $T_\alpha x\to x$ for all $x\in X$
\iitem{\rm (ii)} $\varlimsup_\alpha \|T_\alpha x+ (I-T_\alpha^*) y^{**}\| 
\le \max \{\| x\|,\|y^{**}\|\}$ for all $x\in X$ and $y^{**}\in Y^{**}$.
\iitem{\rm (ii$'$)} In the operator-space setting, the $T_\alpha$'s are 
complete contractions so that for all $n$ and $n\times n$ matrices 
$(x_{ij})$ and $(y_{ij}^{**})$ in $M_n(X)$ and $M_n(Y^{**})$ respectively, 
$$\varlimsup_\alpha \|(T_\alpha x_{ij} + (I-T_\alpha^{**}) y_{ij}^{**})\| 
\le \max \{ \|(x_{ij})\|, \|(y_{ij})\|\}\ .$$
\endroster
Then $X$ is an $M$-ideal (resp. complete $M$-ideal) in $X^{**}$. 
Moreover then $T_\alpha^*\to I_{X^*}$ in the WOT. 

\endproclaim

\re{Remarks} 
1. Of course II is weaker than I, if we assume the $T_\alpha$'s are all 
contractions (resp. complete contractions). 

2. It follows immediately that if $(T_\alpha)$ satisfies I or II, then also 
there exists a net $(\tilde T_\alpha)$ of weakly compact operators on $X$ 
satisfying I or II, so that in addition $\tilde T_\alpha^* \to I_{X^*}$ 
in the SOT. 
\endre

In particular, we recover the facts (in this setting) that if $X$ is 
separable and satisfies I or II, $X^*$ is separable; if the $T_\alpha$'s 
are all compact, then $X^*$ has the bounded (resp. metric in case~II) 
compact approximation property, and finally if the $T_\alpha$'s are 
finite rank, then $X^*$ has the bounded (resp. metric in case~II) 
approximation property. 
In particular, if $X$ is separable and the $T_\alpha$'s are finite rank, 
then $X^*$ is separable with the bounded (and hence metric) approximation 
property. 

This suggests the conjecture: 
if $X$ admits a weak M-ai in $X^{**}$ consisting of finite rank operators, 
then also $X$ admits a weak M-ai in $X^{**}$ consisting of contractive 
finite rank operators.

\demo{Proof of Corollary 1.3} 
Actually, everything but the final statement follows immediately from 
Theorem~1.1 in case I holds, and Theorem~1.1$'$ in case~II. 
(Simply note that since the $T_\alpha$'s are weakly compact, 
$T_\alpha^{**} X^{**}\subset X$ for all $\alpha$). 

Now in both cases, we obtain that setting $Y=X^{**}$ and $S_\alpha = 
T_\alpha^{**}$, then $S_\alpha^* \to P$ in the $W^*$-OT, where $P$ is 
the $L$-projection on $Y^*$ with kernel $X^\bot$. 
But in this case, we must have that $P(Y^*) = X^*$, i.e., the $L$ 
decomposition of $Y^* = X^{***} = X^\bot \oplus X^*$ 
(see \cite{HWW}). 
But then after taking the various identities into account, we simply have 
that if $x^* \in X^* \subset X^{***}$, $T_\alpha^{***} x^* = T_\alpha^*x^*$, 
whence $T_\alpha^{***} x^* \to x^*$ weak* simply means that 
$\la x^*,T_\alpha^* x^*\ra \to \la x^{**}, x^*\ra$ for all $x^{**}\in X^{**}$;
i.e., $T_\alpha^* x^*\to x^*$ weakly.\qed
\enddemo

We now pass to the strongly motivating case of ideals in $C^*$-algebras. 
The proof uses certain standard arguments in $C^*$-algebras, for which we 
nevertheless give details for the sake of completeness. 

\proclaim{Proposition 1.4} 
Let $\J$  be an ideal in a $C^*$-algebra $\A$. 
Then there is a strong contractive M-cai $(T_\alpha)$ for $\J$ in $\A$. 
\endproclaim

\demo{Proof} 
We may assume that $\A$ is unital, by simply adjoining an identity, 
denoted $1$. 
For once the result is proved here, its hereditary character yields the 
non-unital case. 
Choose a net $(x_\alpha)_{\alpha\in\A}$ of elements of $\J$ satisfying the 
following properties: 
$$\gather
0\le x_\alpha \le 1\ \text{ for all } \ \alpha \tag 1.27\\
x_\alpha x \to x\ \text{ for all }\ x\in \J\ .\tag 1.28\\
x_\alpha y- yx_\alpha \to0\ \text{ for all }\ y\in \A\ .\tag 1.29
\endgather$$
(Such a net $(x_\alpha)_{\alpha\in\gamma}$ is called a {\it central 
approximate unit\/} for $\J$. 
For the existence of such nets, see W.B.~Arveson \cite{Ar} and 
C.~Akemann-G.~Pederson \cite{AP}. See also our proof of Theorem~3.1 below.) 

Now define $T_\alpha :\A\to \A$ by $T_\alpha (y) = x_\alpha y$ for all 
$y\in \A$ and $\alpha$. 
We claim that $(T_\alpha)$ is the desired net. 
Now (i), (ii), and (v) of Definition~1.1 are immediate, 
since $\J$ is an ideal, (1.27) holds, 
and $T_\alpha$ is a complete contraction for all $\alpha$. 
Using standard facts about $C^*$-algebras, we also have (iv). 

Indeed, $\A^{**}$ is in fact a von-Neumann algebra acting on a certain 
Hilbert space $\H$, in which the $w^*$-topology on bounded sets (w.r.t.\ 
$\A^*$) coincides with the weak operator topology with respect to $B(\H)$. 
Now letting $e$ denote the unit-element of $\J^{**}$ (which exists since 
$\J$ is a von-Neumann algebra), it follows that $x_\alpha\to e$ $w^*$, 
i.e., $x_\alpha\to e$ in the WOT on $B(\H)$, whence given $x^{**}\in \J^{**}$, 
also $x_\alpha\cdot x^{**} \to e\cdot x^{**} = x^{**}$ in the WOT, 
i.e., $T_\alpha (x^{**}) = x_\alpha \cdot x^{**}\to x^{**} w^*$, proving (iv).

Now the results in \cite{Ar} yield that for all $a \in \A$, 
$$\sqrt{x_\alpha} a - a\sqrt{x_\alpha} \to 0\ \text{ and }\ 
\sqrt{1-x_\alpha} a - a\sqrt{1-x_\alpha} \to 0\ .
\tag 1.30$$ 
(This follows also immediately from the known inequality: 
$\|\sqrt{x}\, a-a\sqrt{x}\,\| \le 2\sqrt{\|a\|} \, \|xa-ax\|^{1/2}$ for all 
$x,a$ in a $C^*$ algebra with $x\ge0$; cf. \cite{D, page 73}.) 
It then follows that for any $a\in \A$, 
$$ x_\alpha a - \sqrt{x_\alpha}\, a\sqrt{x_\alpha} \to0
\tag 1.31i$$
and 
$$(1-x_\alpha) a - \sqrt{1-x_\alpha} a \sqrt{1-x_\alpha}\to0
\tag 1.31ii$$
Indeed, by (1.25), $\sqrt{x_\alpha} (\sqrt{x_\alpha} a-a\sqrt{x_\alpha})\to0$, 
and $\sqrt{(1-x_\alpha}) (\sqrt{1-x_\alpha} a- a\sqrt{1-x_\alpha})\to0$, 
yielding (1.31). 

Now for each $\alpha$, define $U_\alpha$ and $V_\alpha$ on $\A$ by 
$$U_\alpha a = \sqrt{x_\alpha} a\sqrt{x_\alpha}\ \text{ and }\ 
V_\alpha a  = \sqrt{1-x_\alpha} a \sqrt{1-x_\alpha}\text{ for all } a\in\A\ .
\tag 1.32$$
Then note that $U_\alpha$ and $V_\alpha$ are also complete contractions and 
moreover $U_\alpha \A \subset \J$ and $(I-V_\alpha) \A\subset \J$ for all 
$\alpha$. 
(See Remark~1 following the proof for the last assertion.) 

Next, define $S:\A\oplus \A\to \A$ by $S(u\oplus v) = u+v$ for all $u,v\in\A$.
then endowing $\A\oplus\A$ with the $L^\infty$-direct sum norm, 
and fixing $\alpha$, we claim that 
$$S\circ (U_\alpha \oplus V_\alpha)\text{ is a complete contraction.}
\tag 1.33$$
This follows immediately from the matrix formula: for $u,v\in\A$, 
$$U_\alpha u + V_\alpha v = (\sqrt{x},\sqrt{1-x}) 
\pmatrix u&0\\ 0&v\endpmatrix 
\pmatrix \sqrt x \\ \sqrt{1-x}\endpmatrix
\tag 1.34$$ 
and the easily seen fact that 
$\|({\sqrt{x}\atop \sqrt{1-x}})\| =  \|(\sqrt{x},\sqrt{1-x})\| =1$.

It now follows that $(T_\alpha)$ is a strong M-cai; that is that (iii$'$) 
of Definition~1.1 holds (since we have verified all the other conditions). 
Indeed, rephrasing (1.31), we have that for all $a\in \A$, 
$$T_\alpha a - U_\alpha a \to 0\ \text{ and }\  (I-T_\alpha)a-V_\alpha a\to0\ .
\tag 1.35$$
Hence for any $n$ and $(a_{ij}),(b_{ij})$ in $M_n(\A)$, we have that 
$$\align 
\varlimsup_\alpha \| (T_\alpha (a_{ij}) + (I-T_\alpha) (b_{ij}))\|
&  = \varlimsup_\alpha \|(U_\alpha (a_{ij}) + V_\alpha (b_{ij}))\|\tag 1.36\\
& \le \max \{\|(a_{ij})\|,\|(b_{ij})\|\}
\endalign$$
by (1.33).\qed
\enddemo

\re{Remarks}
1. It is trivial that $U_\alpha\A\subset\J$ for all $\alpha$, 
since if $0\le x\le1$
belongs to $\J$, then $\sqrt{x}\in\J$, whence $\sqrt{x}a\sqrt{x}\in \J$ 
for any $a\in\A$ since $\J$ is an ideal. 
A less trivial fact; also $(I-V_\alpha)\A \subset \J$ for all $\alpha$. 
This is so because for $x$ as above, if $a\in\A$, then also $a-\sqrt{1-x} a
\sqrt{1-x}\in \J$. 

Here is a simple proof of this fact: 
If suffices to show that 
$$\sqrt{1-x} a- a\sqrt{1-x} \in \J\ .
\tag $*$ $$
For then it follows that 
$$\sqrt{1-x} (\sqrt{1-x} a -a\sqrt{1-x}) = (1-x) a-\sqrt{1-x}a \sqrt{1-x}
\in \J\ .
\tag $**$ $$
But $xa\in\J$, whence $a-\sqrt{x} a\sqrt{1-x}\in \J$ as desired. 
But in fact, we have that for any continuous function $f:[0,1]\to\cee$, 
$f(x) a-af(x)\in\J$!
Indeed, the family $\F$ of all such functions $f$ is clearly a closed linear 
subspace of $C([0,1])$ which contains the constants and all powers of $t$; 
$t\to t^n$, whence $\F$ contains all polynomials, so $\F =C([0,1])$ 
by the Weierstrass approximation theorem.

%

2. In Proposition 2.2 of the next section, we obtain a stronger form of 
M-cai's in this context. 
\endre

Applying Corollary 1.2, we thus obtain 

\proclaim{Corollary 1.5} 
Let $\J$ be an ideal in $\A$ a $C^*$-algebra. 
There exists a central approximation unit $(u_\alpha)$ in $\J$ with 
$0\le u_\alpha \le 1$ for all $\alpha$ so that setting $U_\alpha  a=u_\alpha a$
for all $a$ in $\A$, and letting $P$ be the $L$-projection on $\A^*$ with 
kernel $\J^{\bot}$, then $U_\alpha^*\to P$ in the SOT; in particular 
(1.22) holds (where ``$Y$'' $=\A$). 
\endproclaim

We may crystallize some of the ideas in the above proof via the following 
notion. 

\de{Definition 1.3} 
Let $X\subset Y$ be operator spaces and $(U_\alpha)_{\alpha\in\D}$ be a net 
of operators on $Y$. 
Say that $(U_\alpha)$ is a {\it good\/} M-cai if the following conditions hold:
\roster
\item"(i)" $U_\alpha Y\subset X$ for all $\alpha$.
\item"(ii)" $U_\alpha x\to x$ for all $x\in X$.
\item"(iii)" There exists a net $(V_\alpha)_{\alpha\in\D}$ of operators 
on $Y$ so that 
\iitem{\rm (a)} $U_\alpha (y) + V_\alpha (y) \to y$ for all $y\in Y$,
\iitem{\rm (b)} $(I-V_\alpha) Y\subset X$ for all $\alpha$,
\iitem{\rm (c)} $S\circ (U_\alpha \oplus V_\alpha)$ is a complete contraction 
for all $\alpha$, where $S:Y\oplus Y\to Y$ denotes the sum operator 
$S(u\oplus v) = u+v$ and $Y\oplus Y$ is a complete $L^\infty$-direct sum. 
\endroster
\endde

\re{Remark} 
Condition (b) of 1.3 yields that good M-cai's $(U_\alpha)$ 
are hereditary; that is, 
if $Z$ is an operator space with $X\subset Z\subset Y$, then also 
$V_\alpha Z \subset Z$ for all $\alpha$ 
and hence $(U_\alpha|Z)$ is a good M-cai for $X$ in $Z$.
\endre

Now the proof of Proposition 1.4 yields that $X$ admits a good M-cai in 
$Y$ if $Y$ is a $C^*$-algebra and $X$ is an ideal in $Y$. 
It also easily yields

\proclaim{Corollary 1.6} 
Suppose $(U_\alpha)$ is a good M-cai for $X$ in $Y$. 
Then $(U_\alpha)$ is a contractive M-cai.
\endproclaim 

We conclude this section with a permanence property for good M-cai's. 
It's motivation comes from the following open problem. 
If $X\subset Y$ are operator spaces with $X$ a complete $M$-ideal in $Y$, 
is $Z\otimes_{\op} X$ a complete $M$-ideal in $Z\otimes_{\op}Y$, 
for all operator spaces $Z$?

\proclaim{Proposition 1.7}
Let $X,Y$ and $Z$ be operator spaces with $X\subset Y$ and suppose $X$ 
admits a good M-cai in $Y$. 
Then $Z\otimes_{\op} X$ admits a good M-cai in $Z \otimes_{\op}Y$. 
Hence $Z\otimes_{\op} X$ is a complete $M$-ideal in $Z\otimes_{\op} Y$.
\endproclaim

\demo{Proof} 
Let $(U_\alpha)_{\alpha\in\D}$ be a good M-cai for $X$ in $Y$. 
We claim that $(I_Z\otimes U_\alpha) \defeq (\tilde U_\alpha)$ 
is then a good M-cai for $Z\otimes_{\op} X$ in $Z\otimes Y$. 
Let $(V_\alpha)_{\alpha\in\D}$ satisfy 1.3(iii) and set $\tilde V_\alpha 
= I_Z\otimes V_\alpha$ for all $\alpha$. 
Now it is immediate that setting $\tilde X = Z\otimes_{\op}X$ and 
$\tilde Y = Z\otimes_{\op} \tilde Y$, then $(\tilde U_\alpha)$ and 
$(\tilde V_\alpha)$ satisfy (i) and (iii)(b) of 1.3. 
It also follows that (iii)(c) holds. 
Indeed, let $\tilde S:\tilde Y\oplus \tilde Y\to\tilde Y$ be the sum operator. 
Then 
$$\tilde S\circ (\tilde U_\alpha \oplus \tilde V_\alpha) = I_Z\otimes 
S\circ (U_\alpha \oplus V_\alpha)\ ,$$
hence $\tilde S\circ (\tilde U_\alpha \oplus \tilde V_\alpha)$ is a complete 
contraction since $S\circ (U_\alpha\otimes V_\alpha)$ has this property.

It remains to check the approximation condition (ii) and (iii)(a). 
But we easily obtain that for $x'\in Z\otimes X$ (the algebraic tensor 
product), $U'_\alpha x' \to x'$, whence since $Z\oplus X$ is dense in 
$Z\otimes_{\op} X$ and the $U'_\alpha$'s are complete contractions, (ii) 
holds. 
The identical density argument establishes (iii)(c) (since again 
$\tilde S\circ (\tilde U_\alpha \oplus \tilde V_\alpha)$ is a complete 
contraction for all $\alpha$).\qed
\enddemo 

\re{Remark} We may also introduce a weaker version of good M-cai's and 
obtain a similar permanence property. 
Given $X\subset Y$ spaces, let us say that $(U_\alpha)$ satisfies $(*)$ 
if $(U_\alpha)$ satisfies (i), (ii), and (iii)(a),(b) of 1.3, but instead 
of (iii)(c), we have 
\vskip1pt
\iitem{(c$'$)} $S\circ (U_\alpha| X\oplus V_\alpha) : X\oplus Y\to Y$ 
is a complete contraction for all $\alpha$.
\endre

Then we obtain again that if $(U_\alpha)$ satisfies $(*)$, then $(U_\alpha)$ 
is a weak contractive M-cai for $X$ in $Y$, hence $X$ is a complete 
$M$-ideal in $Y$. 
Moreover if $Z$ is an arbitrary operator space,  
since $X$ admits a family satisfying 
$(*)$, so does $Z\otimes_{\op} X\subset Z\otimes_{\op} Y$, whence again, 
$Z\otimes_{\op} X$ is a complete $M$-ideal in $Z\otimes_{\op} Y$. 

\head \S2. Complementation results \endhead

The main (motivating) result of this section is as follows. 

\proclaim{Theorem 2.1} 
Let $\J \subset Y\subset \A$ with $\J$ an  approximately 
injective ideal in 
a $C^*$-algebra $\A$ and $Y$ a $\lambda$-locally reflexive operator space with 
$Y/\J$ separable. 
Then for every $\ep>0$, there exists a completely bounded lift 
$L:Y/\J \to Y$ of $I_{Y/\J}$ with $\|L\|_{\cb} < \lambda +\ep$. 
\endproclaim 

When  $Y =\A$, $\lambda=1$ (necessarily); our result then generalizes (up to 
$\ep$) the theorem of E.~Effros and U.~Haagerup, which yields that then, 
assuming $\A$ is unital, there exists a completely positive lift 
$L:\A/\J\to Y$ of $I_{\A/\J}$ \cite{EH}. 
We give an isometric operator-space generalization of the Effros-Haagerup 
lifting theorem in the Appendix. 

To prove Theorem 2.1, we use the stronger properties of the M-cai's for ideals 
in $C^*$-algebras obtained in the proof of Proposition~1.4. 

\de{Definition 2.1} 
Let $X\subset Y$ be operator spaces and $(U_\alpha)_{\alpha\in\D}$ 
a net of operators on $Y$ be given. 
Say that $(U_\alpha)$ is a special M-cai if the following conditions hold.
\roster
\item"(i)" $U_\alpha Y\subset X$ for all $\alpha$. 
\item"(ii)" $U_\alpha x\to x$ for all $x\in X$.
\item"(iii)" For all $y\in Y$, $U_\alpha^{n+1} (y) - U_\alpha^n(y)\to0$ 
as $n\to \infty$, uniformly in $\alpha$.
\item"(iv)" For every positive integer $k$, there exists a net 
$(V_\alpha^{(k)})_{\alpha\in\D}$ of operators on $Y$ so that 

\iitem{\rm (a)} $U_\alpha^k (y) + V_\alpha^{(k)} (y)\to y$ for all $y\in Y$.
\iitem{\rm (b)} $S\circ (U_\alpha^k \oplus V_\alpha^{(k)})$ is a complete 
contraction for all $\alpha$, where $S:Y\oplus Y\to Y$ denotes the 
``sum'' operator, $S(u\oplus v) = u+v$ for all $u,v\in Y$, and $Y\oplus Y$
is a complete $L^\infty$-decomposition. 
\iitem{\rm (c)} $(I-V_\alpha^{(k)}) Y\subset X$ for all $\alpha$.
\endroster
\endde

Our next result yields that $(U_\alpha)$ is a special M-cai precisely when 
$(U_\alpha)$ satisfies (iii) and all powers $(U_\alpha^k)$ of $(U_\alpha)$ 
are good M-cai's as given in Definition~1.3.
\proclaim{Proposition 2.2} 
Let $X\subset Y$ be operator spaces and $(U_\alpha)$ a net of operators 
on $Y$. 
\roster
\item"A." If $(U_\alpha)$ satisfies condition (i), (ii) and (iv) of 
Definition~2.1, then for all $k$, $(U_\alpha^k)$ is a contractive 
M-cai for $X$ in $Y$.
\item"B." If $(U_\alpha)$ is a special M-cai, then also $(U_\alpha^k)$ is a 
special M-cai for all $k$. 
\item"C." If $Z$ is a closed subspace of $Y$ with $X\subset Z$, then 
$(U_\alpha|_Z)$ is a special M-cai for $X$ in $Z$.
\endroster
\endproclaim

\demo{Proof} 
Condition (iv)(b) implies that the 
$U_\alpha$'s are complete contractions, hence 
so are the $U_\alpha^k$'s. 
Clearly $(U_\alpha^k)$ satisfies (i) for any $k$. 
We easily see that $(U_\alpha^k)$ satisfies (ii) by induction on $k$. 
$k=1$ follows by definition. 
Assuming this holds for $k$, then for $x\in X$, 
$$\align 
\| U_\alpha^{k+1}x-x\| 
&= \|U_\alpha^{k+1} (x) - U_\alpha^k (x) + U_\alpha^k (x) - U_\alpha(x)\|\\
&\le \|U_\alpha^k\|\, \|U_\alpha (x) - (x)\| 
+ \|U_\alpha^k (x) - U_\alpha (x)\|\ .
\endalign$$
Hence $\lim_\alpha \|U_\alpha^{k+1} (x) -x\| =0$ as desired.

To finish proving A, we need only verify that $(U_\alpha^k)$ satisfies 
condition (iii$'$) of Definition~1.1. 
Choosing $V_\alpha^{(k)}$ as in 2.1(iv), then given $n$ and 
$(u_{ij})$, $(v_{ij})$ $n\times n$ matrices in $Y$, we have that 
$$\lim_\alpha \|\big((I-U_\alpha^k) (v_{ij}) - V_\alpha^{(k)}(v_{ij})\big)\| =0
\ \text{ by (v)(a)},$$
whence 
$$\align 
\varlimsup_\alpha \|\big(U_\alpha^k (u_{ij}) + (I-U_\alpha^k)(v_{ij})\big)\|
&= \varlimsup_\alpha \|\big(U_\alpha^k (u_{ij}) + V_\alpha^{(k)}(v_{ij})\big)\|
\tag 2.1\\
&\le \max\| (u_{ij}) \|\, \|(v_{ij})\|\ .
\endalign$$

Now to prove 2.2B, we only need to verify that $(U_\alpha^k)$ satisfies (iii). 
We first observe that if $(U_\alpha)$ is special, then for any positive 
integer $k$ and $y\in Y$, 
$$(U_\alpha^{n+k} - U_\alpha^n)(y) \to 0\ \text{ as }\ n\to\infty\ ,
\tag 2.2$$ 
uniformly in $\alpha$.

\noindent
To see this by induction: 
the case $k=1$ follows by definition. 
Assuming valid for $k$, then for $y\in Y$, 
$$U_\alpha^{n+k+1} (y) - U_\alpha^{n+k} (y) 
= [U_\alpha^{n+k+1} (y) - U_\alpha^{n+k}(y)] 
+ [U_\alpha^{n+k} (y) - U_\alpha^n (y)]\ .
\tag 2.3$$ 
But $(U_\alpha^{n+k+1} - U_\alpha^{n+k})(y)\to 0$ as $n\to\infty$, 
uniformly in $\alpha$, by definition, and $U_\alpha^{n+k} (y)-U_\alpha^n(y)
\to 0$ as $n\to\infty$, uniformly in $\alpha$, by the induction hypothesis. 

Now it follows directly that $(U_\alpha^k)$ satisfies (iii) for all $k$. 
Indeed, given $y\in Y$, we have for all $\alpha$ that 
$$(U_\alpha^k)^{n+1} (y) - (U_\alpha^k)^n(y) 
= U_\alpha^{kn+k}(y) - U_\alpha^{kn}(y)\to 0\ \text{ as }\ n\to\infty
\tag 2.4$$
uniformly in $\alpha$, by (2.3). 
Thus B is proven. 
C is immediate from (iv)(c), since this implies $V_\alpha^{(k)}Z\subset Z$ 
for all $\alpha$ and $k$.\qed
\enddemo

An inspection of the proof of Proposition~1.4 now yields our main example 
of this phenomenon. 

\proclaim{Proposition 2.3} 
Let $\J$ be an ideal in a $C^*$-algebra $\A$.  
Then there is a special M-cai $(U_\alpha)_{\alpha\in\D}$ for $\J$ in $\A$. 
\endproclaim

\demo{Proof}
As before, we may assume that $\A$ is unital. 
For if say $\J\subset \A_0\subset \A$, with $\J$ an ideal in 
$\A_0$ non-unital and $\A$ is just $\A_0$ with unit adjoined, $\J$ is an ideal 
in $\A$ and then by 2.2(c), $(U_\alpha|_{\A_0})$ serves as the special M-cai 
for $\J$ in $\A_0$.

For any $0\le y\le 1$ in $\A$, define the operator $T_y$ on $\A$ by 
$$T_y (A) = yAy\ \text{ for all }\ A\in\A\ .
\tag 2.5$$ 
Next, let $(x_\alpha)$ be a central approximate unit for $\J$ (i.e., 
$(x_\alpha)$ satisfies (1.27)--(1.29). 
Define $U_\alpha$ by (1.32); that is, 
$$U_\alpha = T_{\sqrt{x_\alpha}}\ \text{ for all }\ \alpha\ .
\tag 2.6$$
We claim that $(U_\alpha)$ is a special M-cai for $\J$ in $\A$. 

We first note that for any positive integer $k$, 
$$(x_\alpha^k)\text{ is a central approximate unit for } \J\ .
\tag 2.7$$
This follows easily by induction on $k$. 
$k=1$ is simply the definition. 
Suppose proved for $k$. 
Then of course $0\le x_\alpha^{k+1}\le 1$ for all $\alpha$. 
Given $x\in\J$. 
$$\eqalign{
&\eqalign{ x_\alpha^{k+1}x - x_\alpha x
	&= x_\alpha^{k+1} x-x_\alpha^k x + x_\alpha^k x- x_\alpha x\cr
	&= x_\alpha^k (x-x_\alpha x) + x_\alpha^k x-x_\alpha x\cr}\cr
&\to 0\ \text{ by induction hypothesis.}\cr}$$
Similarly, given $y\in Y$, 
$$\align 
&x_\alpha^{k+1}y - yx_\alpha^{k+1} 
= x_\alpha (x_\alpha^k y - yx_\alpha^k) + (x_\alpha y - yx_\alpha)x_\alpha^k\\
&\to 0\ \text{ by induction hypothesis.}
\endalign$$

Now we verify (i)--(iv) of Definition~2.1 for $(U_\alpha)$. 
(i) holds via Remark~1 following the proof of Proposition~1.4, (ii) holds 
by (1.31i) and the fact that $(x_\alpha)$ is an approximate identity in $\J$. 
Now fix $k$ a positive integer, and define 
$$V_\alpha^{(k)} = T_{(1-x_\alpha^k)^{1/2}}\ \text{ for all }\ \alpha\ .
\tag 2.8$$
Now (1.31(i),(1.31ii) (applied to ``$x_\alpha$'' $= x_\alpha^k$), 
yield that (iva) holds, while (ivb) holds via (1.32) (again applied to 
''$x_\alpha$'' $=x_\alpha^k$ for all $\alpha$). 
Also, (ivc) holds, again via Remark~1 following the proof of 1.4. 

It finally remains only to verify (iii). 
In fact, we have the stronger condition
$$\|U_\alpha^{n+1}-U_\alpha^n\| \to0\ \text{ as } \ n\to\infty\ ,
\text{ uniformly in }\alpha\ .
\tag 2.9$$ 
Indeed, for any $0\le y\le  1$ in $\A$, we have for all $A\in\A$ that 
$$T_y^n(A) = T_{y^n} (A) = y^n Ay^n\ .
\tag 2.10$$ 
Hence 
$$T_y^{n+1}(A) - T_y^n(A) = y^{n+1}A(y^{n+1}-y^n) + (y^{n+1}-y^n)Ay^n\ .
\tag 2.11$$
Now thanks to the (elementary) operational calculus, we have that 
$$\lim_{n\to\infty} \|y^{n+1}-y^n\| =0\ ,\ \text{ uniformly over $y\in\A$ 
with } 0\le y\le 1\ ,
\tag 2.12$$
Thus given $\ep>0$, we may choose $n$ so that 
$$\|y^{n+1}-y^n\| \le \frac{\ep}2\ \text{ for all }\ 0\le y\le 1\ ,\ y\in\A\ ,
\tag 2.13$$
But then 
$$\align
\|T_y^{n+1}(A)-T_y^n(A)\| &\le \|y^{n+1}\|\, \|A\|\, \|y^{n+1}-y^n\|\tag 2.14\\ 
&\qquad + \|y^{n+1}-y^n\|\, \|A\|\, \|y^n\| \\
&\le \ep\|A\|\ .
\endalign$$
Thus 
$$\|T_y^{n+1}-T_y^n\| \le\ep\ ,\ \text{ for all } y\in\A\ ,\ 0\le y\le 1\ .
\tag 2.15$$
But of course $U_\alpha^n = T_{x_\alpha^{1/2}}^n$ for all $\alpha$, 
proving (2.9).\qed
\enddemo

\re{Remark}
Buried in the proof of this result, as well as the proof of Proposition~1.4, 
is the following elementary fact: 
For any $0\le y\le1$ in $\A$, $S\circ (T_{\sqrt y} \oplus T_{\sqrt{1-y}})$ 
is a completely positive contraction. 
\endre

(The complete positivity is evident, upon explicitly writing this map as 
$$u\oplus v \to \sqrt{y}\ u\sqrt{y} + \sqrt{1-y}\ v\sqrt{1-y}\ .$$
Hence, assuming that $\A$ is unital with identity $I$, we need only compute 
$\|S\circ (T_{\sqrt y} \oplus T_{\sqrt{1-y}}) (I\oplus I)\|$; but of course 
this equals 
$$\|\sqrt y \sqrt y + \sqrt{I-y}\sqrt{I-y}\| 
= \|y+I-y\| =1\ .\ )$$

We may now give the proof of Theorem~2.1, via the following more general 
result. 

\proclaim{Theorem 2.4} 
Let $\lambda\ge1$, and let $X\subset Y$ be operator spaces with $X$ 
approximately injective, $Y$ $\lambda$-locally reflexive, and $Y/X$ separable.
Assume that $X$ admits a special M-cai in $Y$. 
Then for all $\ep>0$, there exists a completely bounded lift $L:Y/X\to Y$ 
of $I_{Y/X}$ with $\|L\|_{\cb} <\lambda+\ep$. 
\endproclaim

\re{Remarks} 
1. $L$ is called a lift of $I_{Y/X}$ if letting $\pi :Y\to Y/X$ be the 
quotient map, then the following diagram commutes:
$$\epsfysize=.75truein\epsfbox{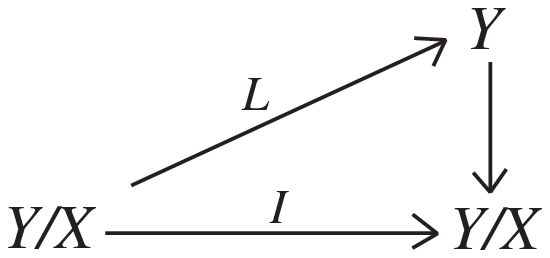}\ .$$
We note that the conclusion is equivalent to the assertion that for all 
$\ep>0$, $X$ is completely $(\lambda+\ep)$-cocomplemented in $Y$; 
that is, there exists a linear projection $P$ mapping $Y$ onto $X$ with 
$\|I-P\|_{\cb}\le\lambda+\ep$.

2. An operator space $X$ is defined to be {\it approximately injective\/} if 
for all finite-dimensional operator spaces $E\subset F$, linear maps 
$T:E\to X$, and $\ep>0$, there exists a linear map $\tilde T:F\to X$ 
extending $T$ with $\|\tilde T\|_{\cb} < \|T\|_{\cb} +\ep$. 
(This is equivalent to the definition given in \cite{EH} when 
$X$ is a $C^*$-algebra.)
We recall that  nuclear operator spaces are approximately injective, 
while approximately injective 1-locally reflexive operators spaces are 
nuclear (cf. \cite{EOR}). 

3. Theorem 2.4 immediately yields the following generalization. 
{\it Let $\lambda$, $X$ and $Y$ satisfy the assumptions of 2.4, but delete 
the hypothesis that $Y/X$ is separable. 
Then for all $\ep>0$, separable operator spaces $Z$ and completely bounded 
maps $T:Z\to Y/X$, there is a lift $\tilde T:Z\to Y$ with 
$\|\tilde T\|_{\cb} \le (\lambda+\ep) \|T\|_{\cb}$.}
Indeed, let $Y_0$ equal the closed linear span of $X$ and $T(Z)$. 
Then $Y_0$ is also $\lambda$-locally reflexive and, since special M-cai's 
are hereditary (Proposition~2.2C), $X$ admits a special M-cai in $Y_0$. 
Thus by Theorem~2.4, there exists a lift $L:Y_0/X\to Y_0$ of $I_{Y_0/X}$ 
with $\|L\|_{\cb} <\lambda+\ep$. 
But then, $\tilde T\defeq L\circ T$ is a lift of $T$ with $\|\tilde T\|\le 
(\lambda+\ep) \|T\|_{\cb}$.\qed
\endre

The following lemma yields the crucial tool for the proof. 

\proclaim{Lemma 2.5}
Let $X$ and $Y$ satisfy the hypotheses of Theorem~2.4.  
Let $Y_1\subset Y_2$ be linear subspaces of $Y$ with $X\subset Y_i$ and 
$E_i \defeq Y_i/X$ finite dimensional, $i=1,2$. 
Let $L_1:E_1 \to Y_1$ be a lift of $I_{E_1}$ and set $\gamma=\|L_1\|_{\cb}$. 
Then given $\ep>0$, there exists a lift $L_2:E_2\to Y_2$ of $I|_{E_2}$ with 
$$\|L_2|_{E_1}-L_1\| <\ep 
\tag 2.16$$
and 
$$\|L_2\|_{\cb} < \max \{\gamma,\lambda\} +\ep\ .
\tag 2.17$$
\endproclaim

\demo{Proof} 
Let $\ep>0$. 
We first note that there exists a lift $L:E_2\to Y_2$ of $I_{E_2}$  with 
$$\|L\|_{\cb} < \lambda +\ep\ .
\tag 2.18$$ 
Indeed, the hypotheses yield that $X$ is a complete $M$-ideal in $Y$, and 
hence $X^{**}$ is completely co-contractively complemented in $Y_2^{**}$. 
Then Proposition 2.6.ii yields yields that $X$ is
completely  $(\lambda+\ep)$-cocomplemented in $Y$,
(see also Sublemma~3.11 of \cite{Ro}). 
Now choosing $P$ a projection from $Y_2$ onto $X$ with $\|I-P\|_{\cb}<
\lambda$ and letting $G= (I-P)X$, it follows that $\pi_{|G}$ maps $G$ 1--1
onto $Y/X$ and $\|(\pi_{|G})^{-1}\|_{\cb} < \lambda+\ep$, whence 
$(\pi_{|G})^{-1}$ is the desired lift $L$. 

Let $(U_\alpha)_{\alpha\in\D}$ be a special M-cai for $X$ in $Y$. 
Now, choose $k$ so that 
$$\|(U_\alpha^{k+1} - U_\alpha^k) L_1\| <\ep\ \text{ for all }\ \alpha\ .
\tag 2.19$$
We may do this, since $L_1(E_1)$ is a finite-dimensional subspace of $Y$, 
using (iii) of Definition~2.1. 
Now choose $(V_\alpha^{(k)})$ as in (iv) 
of Definition~2.1.

Since $L_1$ and $L|_{E_1}$ lift $I_{E_1}$, we have that 
$$(L_1-L)(E_1) \subset X\ .
\tag 2.20$$
Hence since $(U_\alpha^k)$ is a special M-cai, choose an $\alpha$ so that 
$$\|(L_1-L|_{E_1}) - U_\alpha^k(L_1-L|_{E_1}) \|<\frac{\ep}3
\tag 2.21$$
and 
$$\|V_\alpha^{(k)}L|_{E_1} + U_\alpha^k L|_{E_1} -L|_{E_1}\| 
< \frac{\ep}3\ \text{ (by  (ivb)).}
\tag 2.22$$ 
Now note that $U_\alpha L_1(E_1)\subset X$.
Hence by the approximate injectivity of $X$ and the complete contractability 
of $U_\alpha$, we may choose $\theta:L_2\to X$ an extension of $U_\alpha L_1$ 
with 
$$\|\theta\|_{\cb} <\gamma +\ep\ .
\tag 2.23$$
Define $L_2:E_2\to Y$ by 
$$L_2 = V_\alpha^{(k)} L+ U_\alpha \theta\ .
\tag 2.24$$
Then (iv)(b) yields that 
$$\| L_2 \|_{\cb} 
\le \max\{ \|L\|_{\cb} ,\|\theta\|_{\cb}\} 
\le \max \{\lambda +\ep, \gamma+\ep\}
\tag 2.25$$
by (2.18) and (2.19).

Finally, we must estimate the norm of $L_1-L_2|_{E_1}$. 
Now we have that 
$$\left\{ \eqalign{L_2|_{E_1} & = V_\alpha^{(k)} L|_{E_1}+U_\alpha^{k+1}L_1\cr
&= V_\alpha^{(k)} L|_{E_1} + U_\alpha^k L_1 + R_1\cr
&\text{where } \|R_1 \|< \frac{\ep}3 \cr}\right.
\tag 2.26$$
by (2.19). 
But also 
$$\left\{ \eqalign{&V_\alpha^{(k)}L|_{E_1} = L|_{E_1} - U_\alpha^kL|_{E_1} 
+ R_2\cr
&\text{where }\|R_2\| <\frac{\ep}3\cr}\right.
\tag 2.27$$
by (2.22). 
Finally, 
$$\left\{\eqalign{
&L_1 = L|_{E_1} - U_\alpha^k L|_{E_1} + U_\alpha^k L_1 + R_3\cr
&\text{where }\|R_3\| < \frac{\ep}3 \cr}\right.
\tag 2.28$$
by (2.21).
Hence 
$$\eqalign{
&L_2|_{E_1} -L_1 \cr
& = L|_E - U_\alpha^k L|_{E_1} + R_2 + U_\alpha^k L_1 + R_1 
- L|_{E_1} + U_\alpha^k L|_{E_1} - U_\alpha^k L_1 - R_3\cr
& = R_1 + R_2 - R_3\ .\cr}
\tag 2.29$$
At last, 
$$\|L_2| E_1 -L_1\| \le \|R_1\| + \|R_2\| + \|R_3\| < \ep\ ,
\tag 2.30$$
completing the proof.\qed
\enddemo

We are now prepared for the 

\demo{Proof of Theorem 2.4} 
Let $0<\ep<1$ and choose finite-dimensional spaces $E_1\subset E_2\subset
\cdots$ in $Y/X$ with 
$$\overline{\bigcup_j E_j} = Y/X\ .$$
As in the first step of Lemma 2.5, choose $L_1:E_1\to Y/X$ a lift of 
$I|_{E_1}$ with 
$$\|L_1\|_{\cb} <\lambda +\frac{\ep}2\ .$$
Let $n\ge1$ and suppose $L_n :E_n\to Y/X$ has been chosen, lifting 
$I|_{E_n}$ with 
$$\|L_n\|_{\cb} < \lambda + \sum_{j=1}^n \frac{\ep}{2^j}\ .
\tag 2.31$$
Then by Lemma 2.5, we may choose $L_{n+1} :E_{n+1} \to Y/X$ lifting 
$I|_{E_{n+1}}$ with 
$$\|L_{n+1}\|_{\cb} < \|L_n\|_{\cb} + \frac{\ep}{2^{n+1}} 
\le \lambda + \sum_{j=1}^{n+1} \frac{\ep}{2^j} 
\tag 2.32\text{(i)}$$
and 
$$\|L_{n+1}|_{E_n} - L_n \| <  \frac{\ep}{2^n}\ . 
\tag 2.32\text{(ii)}$$
Now it follows that setting $Z= \bigcup_{j=1}^\infty E_j$, then 
$(L_n)$ converges pointwise to a lift $L$ of $I_Z$ satisfying 
$$\|L\|_{\cb} \le \lambda + \sum_{j=1}^\infty \frac{\ep}{2^j}
= \lambda +\ep\ .
\tag 2.33$$ 
To see this, let $z\in E_k$ for some $k$. 
Then for any $k\le m<n$, 
$$\align 
\|L_n(z) - L_m(z) \|
& = \Big\| \sum_{j=m}^{n-1} L_{j+1} (z)-L_j(z)\Big\|\\
&\le \sum_{j=m}^{n-1} \|L_{j+1}|_{E_j} - L_j\|\, \|z\|\\
&\le \frac{\ep}{2^{m-1}} \to 0\ \text{ as }\ m\to\infty\ ,
\endalign$$
hence $(L_n)$ indeed converges to a linear operator $L$ on $Z$. 
But we also see that fixing $k$, then 
since $\|L_n|_{E_k}\|_{\cb} \le\lambda+\ep$ for all $n$, by (2.31); 
also $\|L|_{E_k}\|_{\cb}\le\lambda +\ep$.
Moreover since $L_n|_{E_k}$ lifts $I_{E_k}$, so does $L$. 
Hence $L$ indeed lifts $I_Z$. 
It now remains to simply extend $L$ to all of $Y/X$ by continuity.\qed
\enddemo

\re{Remark} 
Say that a net $(U_\alpha)$ of operators on $Y$ is a 
{\it weak special M-cai\/} 
provided $(U_\alpha)$ fulfills all the conditions of Definition~2.1 
{\it except\/} that we replace (ivb) by 
\roster
\item"(ivb$'$)" $S\circ (U_\alpha^k|_X\oplus V_\alpha^{(k)})$ is a complete 
contraction from $X\oplus Y$ to $Y$ for all $\alpha$.
\endroster
(In other words, for all $k$ $(U_\alpha^k)$ satisfies condition $(*)$ given 
in the Remark at the end of Section~1, and also $(U_\alpha)$ satisfies 
condition~(iii) of Definition~2.1.) 
The {\it proof\/} 
of Theorem 2.4 yields that its conclusion holds provided we assume 
instead that $X$ admits a weak special M-cai in $Y$. 
\endremark

We next take up the problem of ensuring that $X$ is {\it complemented\/} 
in $Y$, when $X\subset Y$ are operator spaces with $Y/X$ separable 
and $X$ approximately injective. 
(It apparently remains an open question if this is always the case in 
this setting.) 
Note, however, that $X$ need not be completely complemented. 
A remarkable example of E.~Kirchberg yields a non-exact separable 
$C^*$ algebra $\A$ and an ideal $J\subset \A$ with $\J$ nuclear and $\A/\J$ 
exact \cite{K}. 
Were $\J$ completely complemented, $\A$ would be $\lambda$-exact for some 
$\lambda$; but then since $\A$ is a $C^*$-algebra, $\A$ would be 
exact (cf. \cite{Pi}). 
Another example, due to T.~Oikhberg and the second author of the 
present paper, yields an example with $X$ completely isometric to $\K$ and 
$Y/\K$ completely isometric to $c_0$ \cite{OR}.

We introduce several new concepts for our investigation.

\de{Definition 2.2} 
Let $X\subset Y$ be Banach/operator spaces and $\lambda \ge1$. 
Let $\pi$ denote the quotient map from $Y$ onto $Y/X$. 
We consider the following diagram, for $E$ a general finite-dimensional 
subspace of $Y/X$ and $i:E\to Y/X$ the identity injection. 
$$\epsfysize=.75truein\epsfbox{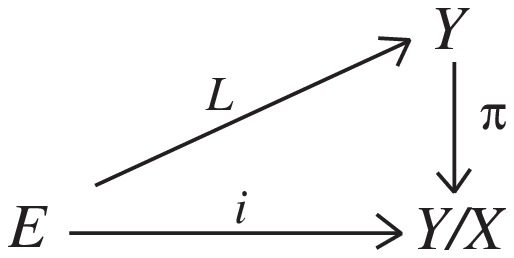}\ .
\leqno (*) $$
That is, $L$ is a lift  of $I_E$ to $Y$.

(i) $(X,Y)$ is said to have {\it $\lambda$-local liftings\/} 
($\lambda$-ll's)  if for all such $E$ and $\ep>0$, there exists a map $L$ 
satisfying $(*)$ with $\|L\| < \lambda+\ep$. 

(ii) $(X,Y)$ is said to have {\it $\lambda$-complete local liftings\/} 
($\lambda$-cll's) if for all such $E$ and $\ep>0$, there exists a map $L$ 
satisfying $(*)$ with $\|L\|_{\cb} < \lambda +\ep$. 

(iii) $(X,Y)$ is said to have {\it $\lambda$-extendable local liftings\/} 
($\lambda$-ell's)  if for all such $E$ and $\ep>0$, there exists a 
map $T:Y/X\to Y^{**}$ with $\|T|| <\lambda+\ep$ so that $L= T|E$ satisfies 
$(*)$. 

Finally, $(X,Y)$ is said to have {\it local liftings\/} (resp. 
{\it complete local liftings\/}) (resp. {\it extendable local liftings\/}), 
if there exists a $\lambda \ge1$ so that $(X,Y)$ has $\lambda$-ll's  
(resp. $\lambda$-cll's) resp. $\lambda$-ell's). 

As we show below, if $X^{**}$ is completely complemented in $Y^{**}$ and 
$Y$ is locally reflexive, then $(X,Y)$ has complete local liftings.   
The following is thus a strengthening of Theorem~2.4 (cf. the Remark 
following its proof for the definition of weak special M-cai's). 
\endde

\proclaim{Theorem 2.4$'$}
Let $X\subset Y$ be operator spaces. 
The conclusion of Theorem~2.4 holds if one replaces in its hypotheses 
the assumption that $Y$ is $\lambda$-locally reflexive, by the assumption that 
$(X,Y)$ has $\lambda$-complete local liftings, 
and that instead $X$ admits a weak special M-cai in $Y$.
\endproclaim

In fact, the proof of Theorem 2.4 gives this immediately; one only needs to 
observe that the local reflexivity assumption on $Y$ is used solely to 
produce the lift $L$ in the proof of Lemma~2.5, satisfying (2.16). 
As we show below the existence of this map follows directly from the 
assumption that $(X,Y)$ has  $\lambda$-complete local liftings.

\re{Remark} 
Let $X\subset Y$ be operator spaces. 
Then $(X,Y)$ has complete local liftings if and only if $X$ is 
{\it locally complemented\/} in $Y$; that is, there exists a $\beta\ge1$ 
so that $X$ is $\beta$-completely complemented in $Z$ for all $X\subset Z
\subset Y$ with $Z/X$ finite-dimensional (one then says $X$ is 
$\beta$-locally complemented in $Y$). 
In fact, it is easily seen that if $(X,Y)$ has $\lambda$-cll's, then $X$ 
is $(\lambda + 1+\ep)$-locally complemented in $Y$ for all $\ep>0$, 
while if $X$ is $\beta$-locally complemented in $Y$, then $(X,Y)$ has 
$(\beta+1)$-local liftings
(cf. \cite{R}, \cite{OR} for certain consequences of local complementability). 
The quantitative cll-concept is more appropriate in the context of 
the present work. 
\endre

The next two results list several easily proved permanence properties of the 
concepts introduced in Definition~2.2. 

\proclaim{Proposition 2.6} 
Let $X\subset Y$ be Banach/operator spaces and let $\lambda\ge1$.

{\rm (i)} $(X,Y)$ has $\lambda$-ll's if and only if 
$$\left\{\eqalign{
&\text{there exists a lift $L:Y^{**}/X^{**}\to Y^{**}$ of the}\cr
&\text{identity map on $Y^{**}/X^{**}$ so that $\|L\|\le\lambda$.}\cr}
\right.
\leqno(**) $$

{\rm (ii)} If $(X,Y)$ has $\lambda$-cll's, $(**)$ holds with 
$\|L\|_{\cb} \le\lambda$. 
If $(**)$ holds with $\|L\|_{\cb} \le\lambda$ and $Y$ is $\beta$-locally 
reflexive, then $(X,Y)$ has  $\lambda\beta$-cll's. 

{\rm (iii)} If $(X,Y)$ has $\lambda$-cll's, then $(X,Z)$ has 
$\lambda$-cll's for all operator spaces $Z$ with $X\subset Z\subset Y$.
\endproclaim

\proclaim{Proposition 2.7} 
Let $X,Y$ be as in 2.6. 
Assume that $(X,Y)$ has local liftings. 
Then $(X,Y)$ has extendable local liftings 
under any of the following hypotheses: 
\roster
\item"(a)" $Y/X$ has the bounded approximation property (the bap).
\item"(b)" $Y^{**}$ is an (isomorphically) injective Banach space. 
\item"(c)" $Y$ is extendably locally reflexive.
\item"(d)" $Y$ is an operator space with $Y^{**}$ an (isomorphically) 
injective operator space and $(X,Y)$ has complete local liftings. 
\endroster
In fact, let $\lambda,\beta\ge1$ and assume that $(X,Y)$ has $\beta$-ll's.
Then $(X,Y)$ has $(\beta\lambda)$-ell's provided any of the following 
hold:
\roster
\item"(a$'$)" $Y/X$ has  the $\lambda$-bap.
\item"(b$'$)" $Y^{**}$ is $\lambda$-injective.
\item"(c$'$)" $Y$ is $\lambda$-extendably locally reflexive.
\item"(d$'$)" $(X,Y)$ has $\beta$-cll's and $Y^{**}$ is a
$\lambda$-injective operator space.  
\endroster
\endproclaim

\re{Remark}
We recall (cf. \cite{OR}) 
that a Banach space $Y$ is called {\it $\lambda$-extendably locally 
reflexive\/} ($\lambda$-elr) provided for all $\ep>0$ 
and all finite-dimensional 
subspaces $G\subset Y^{**}$ and $F\subset Y^*$, there exists a linear 
operator $T:Y^{**}\to Y^{**}$ with $\|T\|<\lambda+\ep$, 
$TG\subset Y$, and $\langle Tg,f\rangle = \langle g,f\rangle$ for all 
$g\in G$ and $f\in F$.
\endre

\demo{Proof of Proposition 2.6} 
Let $\pi :Y\to Y/X$ be the quotient map. 

(i) Suppose first that there exists a lift $L$ satisfying $(**)$. 
Let $E$ be a finite-dimensional subspace of $Y/X$ and let $\ep>0$ be given. 
Regarding $Y/X\subset Y^{**}/X^{**}$, let $G= L(E)$. 
Let $Y_0=\pi^{-1}(E)$.  Then of course $X\subset Y_0$,
$Y_0/X=E$, and $G\subset Y^{**}$.  Let $F$ be $X^\bot$ relative to
$Y_0^{**}$.  Then $F$ is finite dimensional and 
for all $x^\bot \in X^\bot$
there exists an $f\in F$ with $x^\bot|_{Y_0^{**}} =
f|_{Y_0^{**}}.$

By the local reflexivity principle, choose $T:G\to Y$ with $\|T\|<1+
\frac{\ep}{\lambda}$ and 
$$\align 
&\|T\| < 1+\frac{\ep}{\lambda}\tag 2.34\\
&\langle Tg,f\rangle = \langle g,f\rangle\text{ for all $g\in G$ and $f\in F$.}
\tag 2.35
\endalign$$
Then $T\circ L|_E$ is our desired lift of $I_E$. 
Indeed, $\|T\circ L\| < (1+\frac{\ep}{\lambda}) \lambda = \lambda+\ep$, 
and for all $e\in E$, 
$$Le - TLe \in X^{**}
\tag 2.36$$ 
thanks to the definition of $F$. 
But then 
$$\align
&\pi TLe = \pi^{**} TLe = \pi^{**} Le \ \text{ by (2.35) }\\
& = e\text{ since $L$ is a lift of $I_{Y^{**}/X^{**}}$.}
\endalign$$
Now suppose $(X,Y)$ has $\lambda$-ll's and let $\D$ be the following 
directed set:
$\D = \{(E,\ep) :E$ is a finite-dimensional subspace of $Y/X$ and $\ep>0\}$, 
where $(E,\ep) \le (E',\ep')$ if $E\subset E'$ and $\ep\ge \ep'$. 
For each $\alpha = (E,\ep)\in \D$, choose $L_\alpha :E\to Y$ a lift of 
$I_E$ with $\|L_\alpha\| <\lambda +\ep$. 
By the Tychonoff theorem, we may choose a subnet 
$(L_{\alpha_\beta})_{\beta\in\D'}$ of $(L_\alpha)$ such that for all 
$e\in Y/X$, 
$$\lim_{\beta} L_{\alpha_\beta} (e) \defeq L(e)$$ 
exists weak* in $Y^{**}$. 
We easily verify that then $L:Y/X \to Y^{**}$ is a linear operator with 
$$\|L\|\le\lambda\ \text{ and }\ \pi^{**}\circ L=\chix
\tag 2.37$$
where $\chix  :Y/X\to (Y/X)^{**}$ is the canonical injection. 
It follows that letting $P:Y^{****} \to Y^{**}$ be the canonical projection, 
then $P\circ L^{**}$ is the desired lift of $(Y/X)^{**}$ into $Y^{**}$. 

(ii) 
The first assertion follows immediately by the argument given in (i), and 
indeed so does the second: we just choose $T$ as in that argument, but 
so that 
$$\|T\|_{\cb} < \beta + \frac{\lambda}{\ep}\ .
\tag 2.38$$ 
Then $\|T\circ L\|_{\cb} < (\beta +\frac{\lambda}{\ep})\lambda = 
\lambda\beta +\ep$ as desired. 

(iii)
Let $\ep>0$ and let $E$ be a finite-dimensional subspace of $Z/X$, 
regarded as a subspace of $Y/X$, and let $L:E\to Y$ be a lift of $I_E$ with 
$\|L\|_{\cb} <\lambda+\ep$. 
But then $L(E) \subset Z$!
Indeed, for $e\in E$, we must have that $\pi\circ L(e)=e$, 
which means there exists a $z\in Z$ so that $L(e)-z\in X$. 
But this says that $Le \in X+Z=Z$. 
This completes the proof of 2.6.\qed
\enddemo

\demo{Proof of Proposition 2.7} 
Of course we just need to prove the quantitative assertions. 
Let $\lambda,\beta \ge1$ and assume that $(X,Y)$ has $\beta$-ll's. 
We show that $(X,Y)$ has $(\beta\lambda)$-ell's under any of (a$'$)--(c$'$).
Let $E$ be a finite-dimensional subspace of $Y/X$ and let $\ep>0$. 
If (a$'$) holds, we may choose $T:Y/X\to Y/X$ a finite-rank operator with 
$\|T\| <\lambda +\frac{\ep}{\beta}$ and $T|E = I|E$. 
By Proposition~2.6, we may choose a lift $L:Y/X\to Y^{**}$ of 
$\chix :Y/X\to (Y/X)^{**}$ with $\|L\|\le\beta$. 
Then $L\circ T:Y/X\to Y^{**}$ satisfies:
$$\|L\circ T\| \le \|L\|\, \|T\| < \beta
\left( \lambda +\frac{\ep}{\beta}\right) = \beta\lambda +\ep\ .
\tag 2.39$$ 
But $L|E = I_E$ (regarding $Y/X \subset (Y/X)^{**}$), hence $(L\circ T) |E 
= I|_E$, and (a$'$) is thus proved. 

Now assume (b$'$) holds, and choose $L:E\to Y$ a lift of $I_E$ with 
$\|L\| <\beta + \frac{\ep}{\lambda}$. 
Since $Y^{**}$ is $\lambda$-injective, choose $\tilde L :Y/X\to Y^{**}$ an 
extension of $L$ with $\|\tilde L\|\le \lambda\|L\| <\lambda\beta+\ep$. 
This proves case~(b$'$). 
Similarly, if (d$'$) holds, we choose $L$ as above so that instead 
$\|L\|_{\cb} <\beta +\frac{\lambda}{\ep}$, and then choose $\tilde L$ with 
$\|\tilde L\|_{\cb} \le \lambda \|L\|_{\cb} < \lambda \beta +\ep$.

Finally, suppose (c$'$) holds. 
By Proposition~2.6, choose $L:Y^{**}/X^{**}\to Y^{**}$ a lift of the identity 
on $Y^{**}/X^{**}$ with $\|L\|\le\beta$. 
Let $W= L(Y^{**}/X^{**})$ and $P = L\circ \pi^{**}$. 
It follows easily that $P$ is a projection from $Y^{**}$ onto $W$ with 
kernel equal to $X^{**}$. 
Let $G= L(E)$, 
and let $F$ be a finite-dimensional subspace of $X^\bot$. 
Now by the definition of $\lambda$-extendable local reflexivity, choose 
$T:Y^{**}\to Y^{**}$ a linear operator so that 
$$\gather
\|T\| < \lambda +\frac{\ep}{\beta}\tag 2.40\\
TG\subset Y\tag 2.41 
\endgather$$
and
$$\langle Tg,f\rangle = \langle g,f\rangle\text{ for all $f\in F$ and 
$g\in G$.}
\tag 2.42 $$
Now let $T_F = \tilde T\defeq T\circ L|(Y/X)$. 
Then 
$$\|\tilde T\| \le \|T\|\, \|L\| < \left(\lambda + \frac{\ep}{\beta}\right) 
\beta = \lambda\beta +\ep\ .$$

By 2.42, identifying $(Y/X)^*$ with $X^\bot$, we have by (2.42), using that 
$L$ is a lift, that 
$$\langle e-\pi T_F e,f\rangle =0\ \text{ for all } \ f\in F\ .
\tag 2.43$$
Now let $\D$ be the set of finite-dimensional subspaces of $X^\bot$, 
directed by inclusion. 
We then obtain by (2.43) that the net $((\pi   T_F)|_E)_{F\in\D}$ converges 
to $I_E$ in the WOT. 
Since $E$ is finite dimensional, convex combinations of this net converge to 
$I|_E$ in norm. 
Thus given $\eta>0$, there exists a convex combination $S$ of the $T_F$'s 
such that $\|I_E- \pi S\|<\eta$. 
Of course $S(E)\subset Y$ and $\|S\|<\lambda\beta +\ep$ also. 
Now suppose $k=\dim E$; choose $(e_i)$ a normalized Auerbach basis for $E$; 
so also there exist normalized $f_i$'s in $X^\bot$ with 
$f_i(e_j) = \delta_{ij}$ and $e= \sum f_i (e)e_i$ for all $e\in E$. 
Now $\| \pi   Se_i-e_i\|<\eta$ for all $i$. 
Hence for each $i$, we may choose $y_i\in Y$ with $ \pi   y_i=e_i$ and 
$\|Se_i-y_i\| <\eta$. 
At last, define $\tilde S:Y/X\to Y^{**}$ by 
$$\tilde S(z) = S(z) + \sum f_i (z) (y_i-Se_i)\ .$$
Then clearly $\|\tilde S\| <\lambda\beta +\ep + k\eta$,
$\tilde S(E)\subset Y$, and if $e\in E$, 
$$ \pi   \tilde S(e) = S(e) + \sum f_i (e) e_i - \sum f_i (e) S(e_i) = e\ .$$
Thus $\tilde S|E$ is indeed a lift of $I_E$, completing  the proof 
(since $\eta>0$ is arbitrary).\qed
\enddemo

\re{Remark} 
The proof of 2.7 case (c) yields that the assumption that $(X,Y)$ 
has  extendable local liftings 
is considerably weaker than the joint assumption that 
$(X,Y)$ has local liftings and $Y$ is extendably locally reflexive.
\endre

Indeed the proof yields the following result:
{\it Suppose that $X\subset Y$ are Banach spaces with $X^{**}$ 
$\beta$-co-complemented in $Y^{**}$, and let $P:Y^{**}\to Y^{**}$ be a 
projection with kernel equal to $X^{**}$, $\|P\|\le\beta$. 
Now assume that for all finite-dimensional subspaces $E\subset Y$, 
$F\subset X^\bot$, and $\ep>0$, setting $G= PE$, there exists an operator 
$T:Y^{**}\to Y^{**}$ with $TG\subset Y$, $\|T\| <\lambda+\ep$, and 
$\langle Tg,f\rangle = \langle g,f\rangle$ for all $g\in G$, $f\in F$.
Then $(X,Y)$ has $\lambda\beta$-extendable local liftings.}

Thus if e.g., $Y$ is separable, we need only find ``extensions'' of 
local-reflexivity operators on a certain countable family of 
finite-dimensional subspaces of $Y^{**}$, rather than all its 
finite-dimensional subspaces, as in the definition of elr.

We are now prepared for the second main result of this paper.
The following qualitative special case provides its main motivation: 
{\it Suppose $\J$ is an approximately injective ideal in a $C^*$-algebra 
$\A$ with $\A/\J$ separable. 
Then $\J$ is complemented in $\A$ provided $(\J,\A)$ has 
extendable local liftings.}

\proclaim{Theorem 2.8} 
Let $X\subset Y$ be operator spaces with $X$ approximately injective. 
Assume that 
\roster
\item"(a)" $(X,Y)$ has $\lambda$-extendable local liftings.
\item"(b)" $X$ admits a weak special M-cai in $Y$. 
\endroster
Then for all operator spaces $Z$ with $X\subset Z \subset Y$ with $Z/X$ 
separable, and for all $\ep>0$, there exists a lift $L:Z/X\to Z$ of 
$I_{Z/X}$ with $\|L\| <\lambda +\ep$.
\endproclaim

\demo{Note} 
Weak special M-cai's are defined in the Remark following the proof of 
Theorem~2.4.
\enddemo

\demo{Proof} 
Let $(U_\alpha)_{\alpha\in\D}$ be a weak special M-cai for $X$ in $Y$. 
We shall define a NEW operator space structure on $Y$ with the following  
properties (where $(X,\OLD)$ denotes the given operator space structure on 
$X$, and $(X,Y)_{\NEW}$ denotes the pair $(X,Y)$ in the NEW structure). 
$$\gather
\text{The identity injection $i:(Y,\NEW)\to (Y,\OLD)$ 
is a semi-isometry.}\tag 2.44\\
(X,\NEW) = (X,\OLD)\ .\tag 2.45\\
(U_\alpha)\text{ is also a weak special M-cai for $X$ in }(Y,\NEW).\tag 2.46\\
(X,Y)_{\NEW} \text{ has $\lambda$-complete local liftings.}\tag 2.47
\endgather$$
(Recall that if $Z$ and $W$ are operator spaces, $T:Z\to W$ is a 
semi-isometry provided $T$ is a norm-preserving complete contraction.) 
Once this is accomplished, the conclusion of Theorem~2.8 follows immediately 
from Theorem~2.4$'$, Proposition~2.2 and Proposition~2.6. 
Indeed, $X$ is now an approximately injective subspace of $(Y,\NEW)$ 
satisfying (2.46) and (2.47), hence by Proposition~2.6(ii), 
$(X,Z)_{\NEW}$ has $\lambda$-cll's, and also (by the proof 
of Proposition~2.2c), $(U_\alpha|Z)$ is a weak special M-cai for $X$ in $Z$. 
Hence by Theorem~2.4$'$, for all $\ep>0$, there exists a lift 
$L:(Z/X)_{\NEW}\to (Z,\NEW)$ of $I_{Z/X}$ with 
$\|L\|_{\cb} <\lambda+\ep$. 
But then of course $\|L\| <\lambda +\ep$, and since (2.44) holds, the 
Banach norm of $L$ is the same in the NEW and OLD structures. 

We define NEW as follows (where $\pi:Y\to Y/X$ is the quotient map): 
For $\tau \in \K\otimes Y$, set 
$$\|\tau \|_{\NEW} = \max \{\|\tau \|, \|(I_{\K}\otimes \pi)(\tau)
\|_{MAX}\}\ .
\tag 2.48$$

We first note that $\|\cdot\|_{\NEW}$ is indeed an operator space 
structure on $Y$, and moreover for all $\tau\in\K\otimes Y^{**}$, 
$$\|\tau \|_{\NEW} = \max \{\|\tau\|,\|(I_{\K}\otimes\pi^{**})(\tau)
\|_{MAX}\}\ .
\tag 2.49$$
(This fact uses only the definition of NEW; none of the other assumptions 
on $(X,Y)$ are needed.) 
To see this, define $T:Y\to Y\oplus Y/X$ ($\ell^\infty$-direct sum) by 
$Ty = y\oplus \pi y$ for all $y\in Y$, and set $Y' = T(Y)$. 
It is trivial that $T:Y\to Y'$ is a surjective isometry. 
Now simply endow $Y\oplus Y/X$ with the $\ell^\infty$-direct sum operator 
space structure of $(Y,\OLD)$ and $(Y/X)_{MAX}$; and call 
this $(Y\oplus Y/X,\NEW)$. 
Now $(Y,\NEW)$ is nothing but the operator-space structure induced 
on $Y'$ by $(Y\oplus Y/X,\NEW)$.
Hence, since $(Y\oplus Y/X,\NEW)$ is an operator space, so is $(Y,\NEW)$. 
But furthermore, 
$$((Y\oplus Y/X)^{**},\NEW) = (Y^{**},\OLD) \oplus (Y^{**}/X^{**}, MAX)\ ,
\tag 2.50$$ 
and so again $(Y^{**},\NEW)$ is nothing but the operator space structure 
induced on $(Y')^{**}$ in $(Y^{**}\oplus Y^{**}/X^{**},\NEW)$, which is of 
course given by (2.49). 

Now it is trivial that (2.44) and (2.45) hold; it remains to verify (2.46). 
Since the Banach-norms in $(Y,\OLD)$ and $(Y,\NEW)$ coincide, all of the 
norm properties  of $(U_\alpha)$ remain valid in $(Y,\NEW)$, so in fact 
we only need verify that given $k$  and $V_\alpha^{(k)}$ satisfying ``weak'' 
(iv) for $(Y,\OLD)$, that also ``weak'' (iv(d)) holds in $(Y,\NEW)$. 
Precisely, we have that $S\circ (U_\alpha^k|_X\oplus V_\alpha^{(k)})$ 
is a complete contraction on $(X\oplus Y,\OLD)$, and we must verify the same 
for $(X\oplus Y,\NEW)$. 

Now by Theorem~1.1, $X$ is a complete $M$-ideal in $Y$; let then $W$ be 
the ($w^*$-closed) linear subspace of $Y^{**}$ such that $X^{**}\oplus W$ 
is a complete $M$-decomposition of $(Y^{**},\OLD)$ and let $R$ be the 
projection from $Y^{**}$ onto $W$ with kernel $X^{**}$. 
It follows that $\pi^{**}$ is a complete surjective isometry from $W$ onto 
$Y^{**}/X^{**}$. 
But then (2.49) yields that 
$$(Y^{**},\NEW) = (X^{**},\OLD) \oplus (W,MAX)
\tag 2.51$$ 
(where we take the complete $\ell^\infty$-direct sum norm in this 
decomposition). 
Now to verify (iva), it is enough to show that for all $\alpha$, 
$$\leqalignno{
&S^{**}\circ (U_\alpha^{k**}|_{X^{**}}\oplus V_\alpha^{(k)**})\text{ is a 
complete contraction from}&(2.52)\cr
&(X^{**}\oplus Y^{**},\NEW) \ \text{ to }\ (Y^{**},\NEW).\cr}$$
Let then $n$ and $(x_{ij}^{**}),(y_{ij}^{**})$ be given $n\times n$ 
matrices of elements of $X^{**}$ and $Y^{**}$ respectively; choose unique 
$\bar x_{ij}^{**}$'s and $w_{ij}$'s in $X^{**}$ and $W$ respectively so that 
$$(y_{ij}^{**}) = (\bar x_{ij}^{**}) \oplus (w_{ij})\ .
\tag 2.53$$
Then 
$$\leqalignno{
&\|(U_\alpha^{k**} x_{ij}^{**} + V_\alpha^{(k)**} y_{ij}^{**})\|_{\NEW}
&(2.54)\cr
&\qquad = \max \{\|(U_\alpha^{k**} x_{ij}^{**} + V_\alpha^{(k)**} 
\bar x_{ij}^{**})\|,
\|(RV_\alpha^{(k)**} w_{ij})\|_{MAX}\}\cr
&\qquad \le \max \{\|(x_{ij}^{**})\|,\|(\bar x_{ij}^{**})\|, 
\|RV_\alpha^{(k)**}\|\, \|(w_{ij})\|_{MAX}\cr
&\qquad = \max \{\|(x_{ij}^{**})\|_{\NEW}, \|(y_{ij}^{**})\|_{\NEW}\ .
\cr}$$
Here, in the above inequality, we have used the fact that $U_\alpha^k$ 
and $V_\alpha^{(k)}$ satisfy ``weak (iva)'', plus the crucial observation 
that since $(W,\NEW)= (W,MAX)$, 
$$\|RV_\alpha^{(k)**}|_W\|_{\cb} = \|RV_\alpha^{(k)**}|_W\| =1\ .$$
The latter holds since the Banach-norm of $V_\alpha^{(k)**}$ in the 
NEW and OLD structures is the same, namely equal to one, since weak iv(b) 
implies 
$V_\alpha^{(k)}$ is a complete contraction in OLD and hence a contraction.

Of course (2.54) now yields that (2.52) holds, completing the proof that 
$(U_\alpha)$ is a weak special M-cai in $(Y,\NEW)$. 

It remains to prove that (2.47) holds. 
Let $E$ be a finite-dimensional subspace of $Y/X$, and let $\ep>0$. 
Since $(X,Y)$ has $\lambda$-ell's, we may choose a linear operator 
$T:Y/X\to Y^{**}$ so that 
$$\gather 
\|T\| <\lambda +\ep\tag 2.55i\\
T(E) \subset Y\tag 2.55ii\endgather$$
and
$$T|_E \text{ is a lift of } I_E\ .\tag 2.55iii$$
We claim that 
$$L\defeq RT|_E 
\tag 2.56$$
is the desired lift. 
Of course $L$ is a lift of $I_{Y/X}$; the crucial point is to compute its 
cb-norm. 
But as we have pointed out above, $(Y/X,\NEW) = (Y/X,MAX)$. 
Hence
$$\|L\|_{\cb} \le \|RT\|_{\cb} = \|RT\| < \lambda +\ep 
\tag 2.57$$
as desired. 
This completes the proof of Theorem~2.8.\qed
\enddemo

\proclaim{Corollary 2.9} 
Let $\J$ be an approximately injective ideal in a $C^*$-algebra $\A$ 
and let $Y$ be a closed linear subspace of $\A$ with $\J\subset Y$ and 
$Y/\J$ separable. 
Then $\J$ is Banach-complemented in $Y$ provided $(\J,\A)$ has 
extendable local liftings. 
In particular, for a given $\lambda\ge1$, then for every $\ep>0$, there  
exists a lift $L:Y/\A \to Y$ of $I_{Y/\A}$ with $\|L\| <\lambda+\ep$ 
provided any of the following hold:
\roster
\item"(i)" $\A/\J$ has the $\lambda$-bounded approximation property.
\item"(ii)" $Y/\J$ has the $\lambda$-bounded approximation property.
\item"(iii)" $Y^{**}$ is a $\lambda$-injective Banach space.
\item"(iv)" $Y$ is $\lambda$-extendably locally reflexive.
\item"(v)" $\A$ is $\lambda$-extendably locally reflexive.
\endroster
\endproclaim

\demo{Proof} 
This is an immediate consequence of our previous work. 
First of all, $\J$ has a special M-cai in $\A$, by Proposition~2.3, and 
of course such is then a weak special M-cai. 
Secondly, since $\J$ is a (complete) $M$-ideal in $\A$, $(\J,\A)$ has 
$\lambda$-ll's by Proposition~2.6(ii), whence also $(\J,Y)$ has this property.
Thus cases (i) and (v) yield that $(\J,\A)$ has $\lambda$-ell's by 
(a$'$) and (c$'$) of Proposition~2.7, while cases (ii)--(iv) yield that 
$(\J,Y)$ has $\lambda$-ell's by (a$'$)--(c$'$) of 2.7. 
Thus Theorem~2.8 yields the conclusion of the Corollary.\qed
\enddemo

\proclaim{Corollary 2.10} 
If $(\K,B(\ell^2))$ has extendable local liftings, then $\K$ is Banach 
complemented in $Y$ for any separable operator space $Y$ with $\K\subset Y$.
\endproclaim 

\head \S3. Examples and complements\endhead

We first consider the case of (closed two sided) ideals $\J$ in non-self 
adjoint operator algebras $\A$. 
We say that a net $(u_\alpha)$ in $\J$ is a contractive approximate identity 
for $\J$ if $\|u_\alpha\| \le1$ for all $\alpha$ and $ u_\alpha x\to x$ for 
all $x\in \J$. 
A remarkable result of Effros-Ruan yields that 
{\it a closed linear subspace $\J$ of an operator algebra $\A$ is an 
$M$-ideal in $\A$ iff $\J$ is an ideal in $\A$ which admits a contractive 
approximate identity\/} \cite{ER1}. 
(The same equivalences were established earlier by R.~Smith, in the case of 
{\it uniform\/} algebras  $\A$ \cite{S}.) 
The discussion in \cite{ER1} easily yields that when this happens, $\J$ is a 
complete $M$-ideal in $\A$. 

We obtain the additional information that these conditions are 
{\it equivalent\/} to $\J$ having an M-cai in $\A$; in fact, we obtain the 
direct generalization of Proposition~1.4 to the non-self-adjoint case. 

\proclaim{Theorem 3.1} 
Let $\J$ be an ideal in an operator algebra $\A$ such that $\J$ has a 
contractive approximate identity. 
Then $\J$ admits a strong contractive M-cai in $\A$.
\endproclaim

\demo{Proof} 
We assume that $\A$ is a (closed) subalgebra of $B(H)$ for some Hilbert 
space $H$. 
We may easily reduce to the case where $\A$ is unital (in which case 
$\A$ may be assumed to be a unital subalgebra of $B(H)$). 
Indeed, if $\A$ is non-unital, simply let $\tilde{\A}$ 
be $\A$ with $I$ adjoined 
(where $\A$ is a non-unital closed subalgebra of $B(H)$). 
Then $\J$ remains an ideal in $\tilde{\A}$; if $( U_\alpha)$ is a strong 
contractive M-cai for $\J$ in $\tilde{\A}$, $( U_\alpha|_{\A})$ is such for 
$\J$ in $\A$. 
Let then $(u_\alpha)$ be a contractive (algebraic) approximate identity for 
$\J$ in $\A$. 
By passing to a subnet, and regarding $\A\subset \A^{**}\subset B(H)^{**}$, 
we may assume that
$$(u_\alpha)_{\alpha\in\D} \text{ converges weak* to an element 
$e^{**}$ of $\A^{**}$.}
\tag 3.1$$

This is nothing but the first step of the proof by Effros-Ruan that the stated 
hypotheses yield that $\J$ is an $M$-ideal in $\A$ (p.919 of \cite{ER1}). 
Now it is proved in \cite{ER1} that then 
$$\align
&e^{**}\text{ is a self-adjoint idempotent in the center of $\A^{**}$,}
\tag 3.2\\
&\text{with $e^{**} x=x$ for all $x\in \J^{**}$.}
\endalign$$

Of course this uses remarkable properties of $C^*$-algebras such as the 
fact that $\A^{**}$ is a $C^*$-subalgebra of $B(H)^{**}$, 
a von-Neumann algebra. 
We also make use of the fact that if $(v_\alpha)$ is a bounded net in 
$B(H)^{**}$ with 
$w^*\!\!\!\!-\!\!\lim_\alpha v_\alpha = v$ ($v\in B(H)^{**}$), then for any 
$w\in B(H)^{**}$ 
$$v_\alpha w  \buildrel {w^*}\over \to vw\quad
\text{ and }\quad  wv_\alpha 
\buildrel {w^*}\over \to wv.
\tag 3.3$$

Finally, the set of positive elements in $Ba(B(H))$ is $w^*$-dense in 
the set of positive elements of $Ba(B(H))^{**}$; hence we may choose a net 
$(e_\alpha)$ in $B(H)$ with 
$$0\le e_\alpha \le 1\ \text{ for all }\ \alpha\ \text{ and }\  
e_\alpha \buildrel {w^*}\over \to e^{**}.
\tag 3.4$$
It then follows that we may choose an appropriate directed set $\D$ and 
``re-labeled'' new nets $(u_\alpha)_{\alpha\in\D}$ and 
$(e_\alpha)_{\alpha\in\D}$ such that 
$$\lim_{\alpha\in\D} u_\alpha = e^{**} = \lim_{\alpha\in\D} e_\alpha\ .
\tag 3.5$$
Hence it follows that 
$$u_\alpha -e_\alpha \to 0\ \text{ weakly.}
\tag 3.6$$
But then we may find a new net of ``far out'' convex combinations of 
$(u_\alpha,e_\alpha)_{\alpha\in\D}$, say $(\tilde u_\alpha,\tilde e_\alpha)_
{\alpha\in\tilde{\D}}$ with 
$$\|\tilde u_\alpha - \tilde e_\alpha \|\to 0\ .
\tag 3.7$$
Of course the $\tilde u_\alpha$'s remain a contractive approximate identity 
in $\J$ and still $0\le \tilde e_\alpha\le1$ for all $\alpha$ with 
$\tilde e_\alpha \to e^{**}$ $w^*$. 
Thus, by re-labelling again we may assume without loss of generality that 
$$\lim_{\alpha\in\D} \|u_\alpha -e_\alpha\|=0\quad \text{ and }\quad 
w^*\!\!\!\!-\!\!\lim_{\alpha\in\D} e_\alpha = e^{**} 
=w^*\!\!\!\!-\!\!\lim_{\alpha\in\D}  u_\alpha.
\tag 3.8$$

Now moreover we have that for any $a\in \A$, 
$$\lim_\alpha u_\alpha a = e^{**}a\ \text{ and }\ 
\lim_\alpha au_\alpha = ae^{**}
\tag 3.9$$
(using 3.3). 
But since $e^{**}$ is central in $\A^{**}$, we obtain that for all $a\in\A$, 
$$\lim_\alpha u_\alpha a - au_\alpha = 0\ \text{ weakly in }\ \A\ .
\tag 3.10$$
At last, by again taking far out convex combinations in our net 
$(u_\alpha,e_\alpha)_{\alpha\in\D}$, we may assume without loss of 
generality that 
$$\lim_{\alpha\in\D} \|u_\alpha a -au_\alpha\| = 0\ \text{ for all }\ 
a\in\A
\tag 3.11$$
and still that (3.8) holds. 

Define then $U_\alpha : \A\to \A$ by 
$$U_\alpha (a) = u_\alpha a\ \text{ for all }\ a\in\A\ .
\tag 3.12$$ 

We shall now prove that $(U_\alpha)_{\alpha\in\D}$ is a strong 
contractive M-cai for $\J$ in $\A$ (by essentially the same argument as 
the proof of Proposition~1.4). 
It is trivial that $U_\alpha$ is a complete contraction for 
all $\alpha$, since $0\le \|u_\alpha\|\le 1$. 

Since the $u_\alpha$'s lie in $\J$, it is trivial that 
$U_\alpha \A\subset\J$ for all $\alpha$, and of course $U_\alpha(x)\to x$ 
for all $x\in\J$ since $(u_\alpha)$ is an approximate identity. 
We obtain that (iv) of Definition~1.1 holds just as in the proof of 
Proposition~1.4. 
Indeed, for any $x^{**} \in \J^{**}$, 
$$U_\alpha^{**} (x^{**}) = u_\alpha  x^{**} 
\buildrel {w^*}\over \to e^{**} x^{**} 
= x^{**}\quad \text{ by (3.2) and (3.3)\ .}$$

To complete the proof, it remains to verify condition (iii$'$) of 
Definition~1.1.
Now by (3.8) and (3.11), 
$$\|e_\alpha a- ae_\alpha\| \to 0\ \text{ for all } \ a\in\A\ .
\tag 3.13$$
Thus we obtain that 
$$\sqrt{e_\alpha} a - a\sqrt{e_\alpha} \to 0\ \text{ and }\ 
\sqrt{1-e_\alpha} a - a\sqrt{1-e_\alpha} \to 0\ \text{ for all }\ 
a\in \A\ .
\tag 3.14$$
(See the comment following (1.30).) 

Then just repeating the proof of (1.31i,ii), we obtain for all $a\in\A$ 
that 
$$e_\alpha a - \sqrt{e_\alpha} a\sqrt{e_\alpha} \to 0
\tag 3.15i$$
and 
$$(1-e_\alpha) a - \sqrt{1-e_\alpha} a \sqrt{1-e_\alpha} \to 0\ .
\tag 3.15ii$$
Now define operators $\tilde U_\alpha$ and $\tilde V_\alpha$ on $B(H)$ by 
$$\tilde U_\alpha y = \sqrt{e_\alpha} y \sqrt{e_\alpha}\ \text{ and }\ 
\tilde V_\alpha y = \sqrt{1-e_\alpha} y \sqrt{1-e_\alpha}\ \text{ for all } 
\ y\in B(H)\ .
\tag 3.16$$
Then by (3.8) and (3.15), for all $a\in \A$, 
$$U_\alpha a - \tilde U_\alpha a\to 0\ \text{ and }\ 
(I-U_\alpha) a-\tilde V_\alpha a \to 0\ .
\tag 3.17$$ 
But as we showed in the proof of Proposition~1.4, 
$$S\circ (\tilde U_\alpha \oplus \tilde V_\alpha) : B(H) \oplus B(H) 
\to B(H)\ \text{ is a complete  contraction,}
\tag 3.18$$ 
where $S(y\oplus z) = y+z$ for all $y,z\in B(H)$. 
Thus, given $n$ and $(a_{ij})$, $(b_{ij})$ in $M_n(\A)$, 
$$\leqalignno{
&\varlimsup_\alpha \|U_\alpha (a_{ij}) + (I-U_\alpha)(b_{ij})\|\cr
&\qquad 
= \varlimsup_\alpha \|\tilde U_\alpha (a_{ij}) + \tilde V_\alpha (b_{ij}) \|
\quad \text{by (3.17)}&(3.19)\cr
&\qquad \le \max \{\| (a_{ij})\|,\|(b_{ij}\|\}\qquad\text{by (3.18).}\cr}$$
\rightline{$\square$}
\enddemo

We next deal with $c_0$ sums of operator spaces. 
Our methods yield many previously obtained results, in this setting, 
via the following very simple result. 

\proclaim{Proposition 3.2} 
{\rm (a)} Let $X_1,X_2,\ldots$ be given operator spaces and 
let $X= (X_1\oplus X_2\oplus \cdots)_{c_0}$, 
$Y= (X_1,\oplus X_2 \oplus \cdots)_{\ell^\infty}$. 
Then $X$ admits a  strong special M-cai in $Y$. 

{\rm (b)} If the $X_j$'s are all approximately injective, so is $X$. 
\endproclaim

\demo{Proof} 
(a) Define $T_n: Y\to X$ by $T_n(y) = x_1 \oplus x_2 \oplus \cdots\oplus x_n$ 
if $y = (x_j)_{j=1}^\infty$, regarding $X_1\oplus X_2 \oplus\cdots\oplus X_n$ 
as canonically embedded in $X$. 
It is then essentially immediate that $(T_n)$ is the desired strong 
M-cai for $X$ in $Y$. 
Indeed, fixing $n$, we have since $T_n$ is a projection, that for any 
$k\ge 1$, 
$T_n^k = T_n$, $(I-T_n)^k = I-T_n$, and in fact $S\circ (T_n\oplus I-T_n): 
Y\oplus Y\to Y$ is a complete contraction, where $S:Y\oplus Y\to Y$ is 
the sum-operator. 
Moreover $T_n^{k+1} - T_n^k \equiv 0$, $T_n(Y) \subset X$, and 
$T_nx\to x$ for all $x\in X$. 
Finally, $X^{**} = (X_1^{**} \oplus X_2^{**} \oplus \cdots)_{\ell^\infty}$ 
and $X^* = (x_1^* \oplus x_2^* \oplus \cdots)_{\ell^1}$. 
Thus if $x^{**} \in X^{**}$, $x^{**} = (x_j^{**})_{j=1}^\infty$ 
$T_n^{**} (x^{**}) = x_1^{**}\oplus \cdots \oplus 
x_n^{**}\to x^{**}$ $\omega^*$, completing the proof of (a).

To prove (b), let $P_n = (T_n - T_{n-1})X$ for all $n\ge1$ (where $T_0=0$), 
i.e., $P_n$ is just the canonical projection onto the $n^{th}$ coordinate. 
Of course $P_n$ is completely contractive for all $n$. 
Now suppose $E\subset F$ are finite-dimensional operator spaces, 
$S:E\to X$ is a given linear map, and $\ep>0$ is given. 
For each $n$, choose $\tilde S_n :F\to X_n$ an extension of 
$S_n \defeq P_nS$ so that 
$$\|\tilde S_n\|_{\cb} \le (1+\ep)\|S_n||_{\cb} \le (1+\ep)\|S\|_{\cb}\ .
\tag 3.20$$
Now simply define $\tilde S : F\to Y$ by 
$\tilde S(f) = (\tilde S_n(f))_{n=1}^\infty$, $f\in F$. 
Let us first note: 
If $f\in F$, $\tilde S(f)\in X$. 
We have that {\it because\/} $E$ is finite-dimensional and $\|S_ne\|\to0$ 
for all $e\in E$, 
$$\|S_n\|_{\cb} \to 0\ .
\tag 3.21$$
Hence 
$$\left\{ \eqalign{
&\|\tilde S_n(f)\| \le \|\tilde S_n\|\, \|f\| \le (1+\ep) \|S_n\|_{\cb}\to0\cr
&\text{as }\ n\to\infty\ .\cr}
\right\}
\tag 3.22$$
Thus $S(F)\subset X$. 
Moreover
$$\|\tilde S\|_{\cb} = \sup_n \|\tilde S_n\|_{\cb} \le (1+\ep) 
\|S\|_{\cb}
\tag 3.23$$
and of course $\tilde S|E=S$, completing the proof.\qed
\enddemo

An example of A.M.~Davie (as refined by W.~Lusky) yields a sequence 
of finite-dimensional Banach spaces $(X_1,X_2,\ldots)$ and a separable 
$Z$ with $X\subset Z\subset X^{**}$ so that $X$ is uncomplemented in $Z$, 
where $X= (X_1\oplus X_2\oplus \cdots)_{c_0}$ 
(cf. Proposition~2.3 of \cite{HWW} for an exposition and the relevant 
references). 
Thus although $X$ admits a strong special M-cai in $Z$ and $Z$ is 1-locally 
reflexive (using the MIN structure), there is {\it no\/} bounded lift 
$L:Z/X\to Y$ of $I_{Z/X}$. 
When the $X_j$'s are all approximately injective and $Z$ is locally 
reflexive, however, we do obtain a (completely) bounded lift via the 
next result.  

\proclaim{Corollary 3.3} 
Let $X_1,X_2,\ldots$, $X$, and $Y$ be as in Proposition 3.2(a), and 
suppose $X\subset Z\subset Y$ with $Z/X$ separable. 

{\rm (a)} $X$ is a complete $M$-ideal in $Y$. 
Hence if the $X_n$'s are all reflexive, $X$ is a complete $M$-ideal in $X^{**}$.

{\rm (b)} If the $X_j$'s are all approximately injective and $Z$ is 
$\lambda$-locally 
reflexive, or more generally, if $(X,Z)$ admits $\lambda$-complete local 
liftings, then  
$X$ is completely $(\lambda+\ep)$-cocomplemented in $Z$ for all $\ep>0$. 
Moreover if $Z$ is 1-locally reflexive, $X$ is completely 
co-contractively complemented in $Z$.

{\rm (c)} If $(X,Y)$ or $(X,Z)$ admits $\lambda$-extendable local liftings 
and the $X_j$'s are approximately injective, $X$ is 
$(\lambda+\ep)$-cocomplemented in $Z$ for all $\ep>0$. 
In particular, this is the case if $Y$ or $Z$ is $\lambda$-extendably locally 
reflexive, or $Y^{**}$ or $Z^{**}$ is a $\lambda$-injective Banach 
space, or $Y/X$ or $Z/X$ has the $\lambda$-bap.
\endproclaim

\demo{Proof}
(a) follows immediately from Theorem~1.1. 
(b) follows immediately from Proposition~3.2(b) together with Theorem~2.4 
in case $\lambda>1$ and $Z$ is $\lambda$-locally reflexive or the 
Theorem in the Appendix in case  $\lambda=1$ (since then  
$X$ is approximately injective and 1-locally reflexive, $X$ is nuclear). 
Theorem~2.8 yields (a) under the $\lambda$-cll's hypothesis. 
Finally, (c) follows immediately from Theorem~1.1 and Theorem~2.8 and 
the (elementary) Proposition~2.7.\qed
\enddemo

In turn, Corollary 3.3 yields as special cases certain theorems discovered 
by the second author of the present paper.

\proclaim{Corollary 3.4 {\rm \cite{Ro}}} 
Let $X_1,X_2,\ldots$ be 1-injective Banach spaces, and $X= (X_1\oplus X_2
\oplus\cdots )_{c_0}$, $Y= (X_1\oplus X_2 \oplus\cdots)_{\ell^\infty}$. 

{\rm (a)} $X$ has the $2$-$\SEP$. In particular, $c_0(\ell^\infty)$ has 
the $2$-$\SEP$.

{\rm (b)} Suppose the $X_j$'s are all 1-injective operator spaces, and let 
$Z$ be an operator space with $X\subset Z\subset Y$ and $Z/X$ separable. 
Suppose that $Z$ is $\lambda$-locally reflexive, or more generally, that 
$(X,Z)$ admits $\lambda$-complete local liftings. 
Then $X$ is completely $(\lambda+\ep)$-cocomplemented in $Z$ for all $\ep>0$. 
If $Z$ is 1-locally reflexive, $X$ is completely co-contractively 
complemented in $Z$.
\endproclaim

\demo{Proof} 
Part (b) follows immediately from Corollary~3.3, since the $X_j$'s 
are thus all approximately injective. 

Part (a) follows from the last statement in (b). 
Indeed, let $\tilde X\subset\tilde Y$ be separable Banach spaces and 
$T:\tilde X\to X$ be a bounded linear operator. 
Then $Y$ is a 1-injective Banach space, hence there exists $\tilde T:\tilde Y
\to Y$ extending $T$. 
Letting $Z$ denote the closed lienar span of $X$ and $\tilde T(\tilde X)$, 
$Z$ satisfies the hypothesis of the final statement in (4), 
since of course we have $Z$ endowed with MIN, which is thus 1-locally 
reflexive. 
Whence there is a linear projection $P:Z\to X$ with $\|I-P\|\le1$, 
so $\|P\|\le 2$; $P\circ\tilde T$ is thus an extension of $T$ with 
$\|T\| \le 2\|T\|$.\qed
\enddemo

\remark{Remark}
Corollary 3.4(a) is obtained (up to $\ep>0$) as Theorem~1.1 of \cite{Ro}. 
Corollary~3.4(b) is a special case of Theorem~3.4 of \cite{Ro}; the 
quantitative result obtained there is not as good. 
However the full qualitative result in \cite{Ro} is more general than 
3.4(b), for it is assumed in \cite{Ro} that the $X_j$'s are 
$\lambda$-injective operator spaces, for some $\lambda\ge1$. 
Thus if $\lambda >1$, it need not be so that the $X_j$'s are 
approximately injective, so the methods of the present paper do not apply. 
\endremark

We next recapture the main results in \cite{Ro} concerning the CSEP
(using also some recent work of L.~Ge and P.~Hadwin \cite{GH}). 
We first recall a concept introduced in \cite{Ro}. 

\de{Definition 3.1} 
Let $C\ge1$. 
A family $\Z$ of operator spaces is said  to be of $C$-finite matrix type 
if for any finite-dimensional operator space $G$, there is an $n=\bn (G)$ 
so that 
$$\|T\|_{\cb} \le C\|T\|_n \text{ for all linear operators } 
T:G\to Z\text{ and all } z\in \Z\ .
\tag 3.24$$
Briefly, we say that $\Z$ is $C$-finite with function $\bn$; a single space 
$Z$ is called $C$-finite provided $\{Z\} $ is $C$-finite. 
\endde

\noindent 
(Recall that for operator spaces $X$ and $Y$ and $T:X\to Y$ a bounded linear 
map, $\|T\|_n = \|I_n\otimes T\|$, where $I_n$ denotes the identity map 
on $M_n$.) 

$C$-finite operator spaces are $C$-locally reflexive, thanks to the following 
interesting  operator space 
extension of the Banach local reflexivity principle, due to 
L.~Ge and  P.~Hadwin. 

\proclaim{Lemma 3.5 {\rm \cite{GH}}}
Let $Y$ be an arbitrary operator space, $\ep>0$, $n\ge1$, and $F,G$ be 
finite-dimensional subspaces of $Y$ and $Y^{**}$ respectively. 
Then there exists a linear operator $T:G\to Y$ satisfying the following:
$$\align
\text{\rm (i)}&\qquad \|T\|_n < 1+\ep \tag 3.25\\
\text{\rm (ii)}&\qquad \langle Tg,f\rangle = \langle g,f\rangle
\text{ for all } g\in G\text{ and } f\in F\ .\\
\text{\rm (iii)}&\qquad T_{|G\cap Y} = I_{|G\cap Y} .\\
\text{\rm (iv)}&\qquad T\text{ is 1--1 and } \|T^{-1}_{|T(G)}\|_n <1+\ep\ .
\endalign$$
\endproclaim

\re{Remark}
The case $n=1$ is precisely the Banach local reflexivity principle as 
formulated in \cite{JRZ}. 
(We only use (3.25)(i)--(iii) in our discussion here.)
We obtain an extension of Lemma~3.5 in Lemma~3.13 below.
\endre

We may now easily obtain the following permanence properties for 
$C$-finite families.

\proclaim{Proposition 3.6}
Let $C,\lambda\ge1$. 

{\rm (a)} Let $X_1,X_2,\ldots$ be operator spaces such that 
$\{X_1,X_2,\ldots \}$ is  $C$-finite. 
Then $(X_1\oplus X_2\oplus\cdots)_{\ell^\infty}$ is $C$-finite.

{\rm (b)} Let $X$ be an operator space, which is $C$-finite 
for all $C>\lambda$. 
Then $X$ is $\lambda$-locally reflexive. 
\endproclaim

\demo{Proof} 
(a) This is a simple consequence of Definition~3.1. 
Let $G$ be a finite-dimensional operator space, and $n=\bn (G)$ be the 
``$n$'' which works for the family $\{X_1,X_2,\ldots\}$. 
Let $T:G\to (X_1\oplus X_2\oplus\cdots)_\infty$ be a linear operator, 
$P_j$ the canonical projection onto $X_j$, and set 
$T_j = P_jT$ for all $j$. 
Then 
$$\align 
\|T\|_{\cb} & = \sup_j \|P_jT\|_{\cb}\tag 3.26\\
&\le \sup_j C\|P_j T\|_n\ \text{ by $C$-finiteness}\\
& = C\|T\|_n
\endalign$$
where the first and last equalities follow by the definition of the 
$\ell^\infty$-operator sum.

(b) Let $X$ satisfy the given hypothesis and let $G$ and $F$ be 
finite-dimensional subspaces of $X^{**}$ and $X^*$ respectively, and let 
$\ep>0$ be given. 
Let $n=\bn (G)$, $n$ the $C$-finiteness function for $X$. 
Choose $C$ with $\lambda <C<\lambda +\ep$, and then choose $\ep'$ with 
$C+\ep' C<\lambda +\ep$. 
Now by Lemma~3.5 (i.e., \cite{GH}), choose $T:G\to X$ satisfying (3.25)
(for ``$\ep$'' $= \ep'$). 
Then 
$$\align 
\|T\|_{\cb} & \le C\|T\|_n \ \text{ since $X$ is $C$-finite}\tag 3.27\\
&\le C(1+\ep')\ \text{ by (3.25)}\\
&< \lambda +\ep\ .
\endalign$$
This completes the proof.\qed
\enddemo 

We are mainly interested here in the case $C=1$ in 3.6, to obtain the 
following 

\proclaim{Proposition 3.7} 
Let $n\ge1$. Then $\ell^\infty (M_{n,\infty}\oplus M_{\infty,n})$ is  
1-locally reflexive.
\endproclaim

\demo{Proof} 
It is obtained in \cite{Ro, Proposition 2.6} that for all $j\ge 1$, 
$M_{\infty,j}$ and $M_{j,\infty}$ are 1-finite, with function $\bn(G)=j\dim G$, 
where $G$ ranges over all finite dimensional operator spaces. 
Thus by Proposition~3.6(a), $\ell^\infty (M_{j,\infty}\oplus M_{\infty,j})$ 
is also 1-finite with the {\it same\/} function $\bn$, and hence is 
1-locally reflexive by Proposition~3.6(b).\qed
\enddemo

The following is the main result concerning the CSEP, obtained in 
\cite{Ro, Corollary~2.7}. 

\proclaim{Corollary 3.8} 
$c_0(M_{n,\infty}\oplus M_{\infty,n})$ has the 2-CSEP for all $n\ge1$.
\endproclaim

\demo{Proof} 
Let $X\subset Y$ be separable operator spaces and $T:X\to Z$ be a 
completely bounded map, where $Z= c_0 (M_{n,\infty} \oplus M_{\infty,n})$. 
Then $Z^{**} = \ell^\infty (M_{n,\infty}\oplus M_{\infty,n})$ is an 
isometrically injective operator space, hence there exists a complete 
bounded extension $\tilde T:Y\to Z^{**}$ of $T$ with $\|\tilde T\|_{\cb} 
= \|T\|_{\cb}$. 
Set $\tilde Y = [\tilde T(Y),Z]$. 
$Z^{**}$ is 1-locally reflexive by  Proposition~3.7, hence so is $\tilde Y$. 
$Z$ is then an approximately injective complete ideal in $\tilde Y$ by 
Proposition~3.2 and Theorem~1.1, whence $Z$ is completely co-contractively 
complemented in $\tilde Y$ by the Theorem in the Appendix. 
Thus choosing $P:\tilde Y\to Z$ a projection with $\|I-P\|\le1$, 
$P\circ\tilde T$ is an extension of $T$ satisfying $\|P\circ \tilde T\|_{\cb}
\le 2\|T\|_{\cb}$.\qed
\enddemo

\re{Remark} 
It is proved in \cite{Ro} (see Corollary~2.5 and Proposition~2.15 of 
\cite{Ro}) that if $Z_1,Z_2,\ldots$ are operator spaces so that for some 
$\lambda$ and $C$, $Z_1,Z_2,\ldots$ have the $\lambda$-CSEP and 
$\{Z_1,Z_2,\ldots\}$ is $C$-finite, then again 
$(Z_1\oplus Z_2\oplus\cdots)_{c_0}$ has the CSEP (in fact the $C\lambda^2 
+\lambda +\ep$-CSEP, all $\ep>0$). 
This result does not follow from the methods of the present paper, in part 
because the $Z_j$'s may not be approximately injective. 
However if the $Z_j$'s have this property, the next result yields better 
quantitative information than the cited results in \cite{Ro}. 
\endremark

\proclaim{Corollary 3.9} 
Let $C\ge 1$ and let $X_1,X_2,\ldots$ be approximately injective operator 
spaces with $\{X_1,X_2,\ldots\}$ $C$-finite; let $X= (X_1\oplus X_2\oplus 
\cdots)_{c_0}$ and $Y= (X_1\oplus X_2\oplus \cdots)_\infty$. 
Then if $X\subset Z\subset Y$ with $Z$ separable, $X$ is completely 
$(C+\ep)$-cocomplemented in $Z$ for all $\ep>0$.
\endproclaim

\demo{Proof} 
$Y$ is $C$-locally reflexive by Proposition~3.6(b), hence the Corollary 
follows by Corollary~3.3(b).\qed
\enddemo

We continue now with further study of the converse of Theorem~1.1. 
We shall show that this is valid in the setting of the appendix, namely 
the case where $X$ is a nuclear complete $M$-ideal in a 1-locally 
reflexive operator space $Y$. 
In fact our result here holds in the more general situation of 
$\lambda$-finitely injective complete $M$-ideals; these include the 
$\lambda$-nuclear ones (cf.\ \cite{KR} for some results concerning the 
latter). 

\de{Definition 3.2} 
Let $\lambda \ge1$ and $X$ an operator space be given. 

(a) $X$ is called {\it $\lambda$-finitely injective\/} if for all operator 
spaces $Y$, finite-dimensional subspaces $G$, $\ep>0$, and linear maps 
$T:G\to X$, there exists a linear extension $\tilde T:Y\to X$ with 
$\|\tilde T\|_{\cb} \le (\lambda+\ep)\|T\|_{\cb}$. 
In case $\tilde T$ can always be chosen finite rank, we shall call $X$ 
{\it $\lambda$-finite rank injective\/}. 

(b) $X$ is called {\it $\lambda$-nuclear\/} if for all finite-dimensional 
subspaces $F$ of $X$ and all $\ep>0$, there exists an $n$ and linear maps 
$U:F\to M_n$ and $V:M_n\to X$ with 
$$\|U\|_{\cb} \|V\|_{\cb} <\lambda +\ep
\tag 3.28i$$
and 
$$V\circ U = i \ , \text{where $i:F\to X$ is the identity injection.}
\tag 3.28ii$$ 
That is, we have the diagram
$$\matrix 
&M_n&\\
&\llap{$\scriptstyle U$} \nearrow\qquad \searrow\rlap{$\scriptstyle V$}&\\
&F\quad\buildrel i\over\longrightarrow\quad X&\ .
\endmatrix$$
Thus, $X$ is 1-nuclear precisely when $X$ is nuclear. 
\endde

\remark{Remark}
(a) It can be proved that in the Banach space category, the 
$\lambda$-nuclear operator spaces coincide with the $\L_\infty$ spaces. 
Precisely, $(X,\MIN)$ is $\lambda$-nuclear for some $\lambda$ if and 
only if $X$ is an $\L_\infty$-space (iff $X^{**}$ is an isomorphically 
injective Banach space), while $(X,\MIN)$ is nuclear iff $X$ is an 
$L^1$-predual (i.e., $X^{**}$ is isometrically injective). 

(b) Results in \cite{Pi} yield that $B(H)$ is {\it not\/} $\lambda$-finite 
rank injective for any $\lambda$ (we are indebted to T.~Oikhberg for 
this fact). 
Note however that trivially, if $X$ is $\lambda$-injective, $X$ is 
$\lambda$-finitely injective. 
\endre

The next simple result shows the connection between the concepts given 
in Definition~3.2 (a) and (b). 

\proclaim{Lemma 3.10} 
Let $\lambda\ge 1$ and $X$ be a $\lambda$-nuclear operator space. 
Then $X$ is $\lambda$-finite rank injective. 
\endproclaim

\demo{Proof} 
Let $F = T(G)$ and choose $n$ and $U:T(G)\to M_n$, $V:M_n\to X$ 
satisfying (3.28). 
Next, using the 1-injectivity of $M_n$, choose $S:Y\to M_n$ an extension 
of $UT$ with 
$$\|S\|_{\cb} = \|UT\|_{\cb}\ .$$
Finally, let $\tilde T = VS$. 
Then by this equality, 
$$\align \|\tilde T\|_{\cb} & \le \|V\|_{\cb} \|U\|_{\cb} \|T\|_{\cb}\\
& < (\lambda +\ep) \|T\|_{\cb}\ \text{ by (3.28ii)}
\endalign$$
and $\tilde T$ is an extension of $T$ by (3.28ii).\qed
\enddemo

The following gives a converse to Theorem~1.1 in the case of 
$\lambda$-finitely injective complete $M$-ideals in 1-locally reflexive 
operator spaces, or in the case of arbitrary ones in spaces with the bap. 

\proclaim{Theorem 3.11} 
Let $X\subset Y$ be operator spaces with $X$ a complete $M$-ideal in $Y$, 
and let $\lambda\ge1$. 
Then $X$ admits an M-cai $(T_\alpha)$ in $Y$ if either of the following hold: 
\roster
\item"{\rm (a)}" $Y$ has the Banach $\lambda$-bap.
\item"{\rm (b)}" $Y$ is 1-locally reflexive and $X$ is $\lambda$-finitely 
injective.
\endroster
In case (a), the $T_\alpha$'s may be chosen to be finite rank operators 
with $\|T_\alpha\|\le\lambda$ for all $\alpha$. 
In case (b), the $T_\alpha$'s may be chosen with $\|T_\alpha\|_{\cb}\le 
\lambda$ for all $\alpha$; moreover if $X$ is $\lambda$-nuclear, again the 
$T_\alpha$'s may be chosen to be finite rank. 
\endproclaim

\re{Remark} 
Thus if $X$ is 1-nuclear and $Y$ is 1-locally reflexive, Theorem~3.11 yields 
that $X$ admits a contractive M-cai in $Y$, consisting of finite rank 
operators.
\endre

We first require an extension of the local reflexivity concept 
(due to S.~Bellenot \cite{Be} in the Banach space category). 

\proclaim{Lemma 3.12} 
Let $X\subset Y$ be  operator spaces with $Y$ $\lambda$-locally reflexive, 
and let $G$ and $F$ be finite-dimensional subspaces of $Y^{**}$ and 
$Y^*$ respectively. 
Then for all $\ep>0$, there exists a linear operator $T:G\to Y$ satisfying 
the following: 
\roster
\item"{\rm (i)}" 
$\|T\|_{\cb} \le \lambda +\ep$. 
\item"{\rm (ii)}" 
$\langle Tg,f\rangle = \langle g,f\rangle$ for all $g\in G$, $f\in F$.
\item"{\rm (iii)}" 
$T_{|G\cap Y} = I_{|G\cap Y}$.
\item"{\rm (iv)}" 
$T(G\cap X^{**}) \subset X$. 
\endroster
\endproclaim

\demo{Comments} 
1. $\lambda$-local reflexivity may be {\it defined\/} as the existence 
of $T$'s satisfying (i) and (ii) only. 
(iii) is known as ``folk-lore''. 
(iv) is the new extension. 

2. We only use Lemma 3.12 here in the case $\lambda=1$.
\enddemo

\demo{Proof} 
We first obtain (i), (ii) and (iv). 
By the basic known equivalences, we may alternatively express $\lambda$-local 
reflexivity as follows: 

For all finite-dimensional operator spaces $G$, 
$$Ba\, \cb  (G,Y^{**}) \subset \lambda\overline{Ba\,\cb(G,Y)}^{w^*}\ .
\tag 3.29$$
(Since $Ba\,\overline{\cb(G,Y)}^{w^*} = Ba(\cb(G,Y)^{**})$ by Goldstein's 
theorem, this is equivalent to the identity map $i:\cb (G,Y^{**})\to 
\cb (G,Y)^{**}$ having the property that $\|i\|_{\cb} \le\lambda$; 
note that $\|i^{-1}\|_{\cb} = 1$ always). 
Now fixing $G\subset Y^{**}$ finite-dimensional, set $H= G\cap X^{**}$ and 
$W = \{T\in \cb(G,Y^{**}):T(H)\subset X^{**}\}$. 
Now of course we may just {\it identify\/} $\cb(G,Y^{**})$ with 
$G^* \otimes_{\op} Y^{**} = G^* \otimes Y^{**}$ algebraically. 
Since $G$ is finite-dimensional, we have 
$$W= H^\bot \otimes Y^{**} + G^* \otimes X^{**}\ .
\tag 3.30$$ 
Then 
$$W_\bot = H\otimes X^\bot\ , 
\tag 3.31$$ 
hence 
$$W_{\bot\bot} = H^\bot \otimes Y + G^* \otimes X\ .
\tag 3.32$$ 
But then $W= \overline{W_{\bot\bot}}^{w^*}$ and 
$$W_{\bot\bot} = \{T\in \cb(G,Y) : T(H)\subset X\}\ .
\tag 3.33$$

Now $I_G \in Ba(W)$, so given $\ep>0$ (since $W=(W_{\bot\bot})^{**}$, again 
applying Goldstein's theorem, and the fact that $I_G \in\overline{\lambda 
(Ba\, W_{\bot\bot})}^*$ by (3.29), for each finite-dimensional subspace 
$\alpha$ of $Y^*$ with $\alpha\supset F$, we may choose $T_\alpha \in 
W_{\bot\bot}$ satisfying (i) and (ii) for ``$\ep$'' $=\ep/2$; 
of course (iv) holds by (3.33) (where ``$T$'' $= T_j$ in (i)--(iv)). 

Now let $\D$ be the family of all finite-dimensional subspaces $\alpha$ 
of $Y^*$ with $\alpha\supset F$, directed by inclusion. 
But then it follows that 
$$T_\alpha g\to g\text{ weakly for all } g\in G\cap Y\ .
\tag 3.34$$
It follows that we may then choose a net $(\tilde T_\alpha)$ of far out convex
combinations of the $T_\alpha$'s so that 
$$\tilde T_\alpha g\to g\ \text{ strongly for all }\ g\in G\cap Y\ .
\tag 3.35$$ 

Of course it follows trivially that each $\tilde T_\alpha$ still satisfies 
(i), (ii), and (iv) (for ``$T$'' $=\tilde T_\alpha$). 

Finally, a standard perturbation argument yields that for one of 
these $\tilde T_\alpha$'s, there exists $T$ a perturbation of $\tilde T_\alpha$ 
satisfying (i)--(iv).\qed
\enddemo

We next obtain an extension of the result of Ge-Hadwin (stated as 
Lemma~3.5 above). 

\proclaim{Lemma 3.13} 
Let $X\subset Y$ be operator spaces, $\ep>0$, $n\ge1$ and $F,G$ be finite 
dimensional subspaces of $Y^*$ and $Y^{**}$, respectively. 
There exists a linear operator $T:G\to Y$ satisfying $(3.25)$ so that 
$$T(G\cap X^{**})\subset X\ .
\tag 3.36$$
\endproclaim

\demo{Proof}
Let $\G$ be the group of $n\times n$ unitary matrices and $m$ its 
associated Haar measure. 
Let $T_n = M_n^*$ (with elements regarded as $n\times n$ matrices) 
and let $F,G$ be finite-dimensional subspaces of $Y^*$ and $Y^{**}$ 
respectively. 
Since $M_n(G)$ is finite-dimensional, we may assume (by enlarging $F$ 
if necessary) that $T_n(F)$ $(1+\ep)$-norms $M_n(G)$; that is, for all 
$\tilde g$ in $M_n(G)$, there exists an $\tilde f\in T_n(F)$ with 
$\|\tilde g\|=1$ and 
$$\|\tilde g\| < (1+\ep)|\langle \tilde g,\tilde f\rangle|\ .
\tag 3.37$$
(We only need this to obtain (3.25 iv), which we didn't really need in 
our subsequent discussion.) 

Now let $\tilde Y = M_n(Y)$, $\tilde X = M_n(X)$, $\tilde G= M_n(G)$, 
and $\tilde F= T_n(F)$. 
By Lemma~3.12 applied in the Banach space category, we obtain a linear 
operator $\tilde T :\tilde G\to \tilde Y$ satisfying (i)--(iv) of 3.12
(where ``$Y$''~$=\tilde Y$, ``$X$''~$= \tilde X$, etc., and 
$\tilde Y^{**}$, $\tilde X^{**}$ are identified with $M_n(Y^{**})$ and 
$M_n(X^{**})$ respectively). 
Now define a linear operator $S$ on $\tilde G$ by 
$$S(\tilde y) = \int_{\G} \tilde T(\tilde y u)u^*\, dm(u)\ ,
\tag 3.38$$
for all $\tilde y\in \tilde G$. 
Now it follows that also $S$ satisfies (i)--(iv) of 3.12 (for 
``$T$''~$=S$, ``$Y$''~$=\tilde Y$, etc.). 
For example, to see that (ii) holds, let $\tilde g\in\tilde G$, 
$\tilde f \in\tilde F$, then 
$$\align 
\langle S\tilde g,\tilde f\rangle 
& = \int_{\G} \langle \tilde T(\tilde gu),u^* \tilde f\rangle \,dm(u)
\tag 3.39\\
& = \int_{\G} \langle \tilde g u,u^* \tilde f\rangle \,dm(u)\\
& = \int_{\G} \langle\tilde g  u u^*,\tilde f\rangle \,dm(u)\\
& = \langle\tilde g,\tilde f\rangle\ .
\endalign$$
The first equality follows since the pairing between $M_n(Y^{**})$ and
$T_n(Y^*)$ is taken so that
$\langle \tilde g b,\tilde f\rangle = \langle \tilde g, b\tilde f\rangle$
for all $\tilde g\in M_n(Y^{**})$, $\tilde f\in T_n(Y^*)$, and
$b\in M_n$.  The second equality follows since for each $u\in \G$,
$\tilde g\in\tilde G$, and $\tilde f\in\tilde F$, we have that
$\tilde g u\in\tilde G$ and $u^* \tilde f \in \tilde F$.

It then follows moreover that 
$$S \text{ is 1--1 and } \|(S_{|\tilde G})^{-1}\| < 1+\ep\ .
\tag 3.40$$
Indeed, if $\tilde g\ne0$ in $\tilde G$, choose $\tilde f\in\tilde F$ 
of norm one, satisfying (3.37); then 
$$\|\tilde S\tilde g\| > |\langle S\tilde g,\tilde f\rangle| 
> (1+\ep)^{-1} \|\tilde g\|
\tag 3.41$$
by (3.37), which yields (3.40).

Now moreover, we have for all $u_0\in\G$ and $\tilde y\in\tilde G$ that 
$$\align 
S(\tilde y u_0) & = \int \tilde T(\tilde y u_0 u)u^*\,dm (u)\tag 3.42\\
& = \int \tilde T (\tilde y v) v^* u_0\,dm(v)\\
& = S(\tilde y) u_0
\endalign$$
where the second equality holds by left translation invariance of Haar 
measure. 
But then 
$$S(\tilde y A) = S(\tilde y)A\ \text{ for all }\ A\in M_n\ .
\tag 3.43$$
In turn, (3.43) yields there is a unique linear operator $T:G\to Y$ 
with $S = I_n\otimes T$ (where $I_n=I_{M_n}$). 
It now follows immediately that $T$ satisfies the conclusion of Lemma~3.13, 
since $S$ satisfies (i)--(iv) of Lemma~3.14 and also satisfies (3.40).\qed
\enddemo

The main part of the proof of Theorem 3.11 is conveniently isolated 
in the following result.

\proclaim{Lemma 3.14} 
Let $\lambda$, $X$ and $Y$ be as in the hypotheses of Theorem~3.11, let $G$ 
be a given finite dimensional subspace of $Y$ with $G\cap X\ne0$, and 
let $n\ge1$, $\ep>0$. 
Let $\alpha= (G,n,\ep)$. 
There exists an operator $T_\alpha :Y\to X$ satisfying the following: 
\roster
\item"(i)" $\|T_\alpha\| < \lambda +\ep $ if (a) of 3.11 holds, or 
\item"(i$'$)" $\|T_\alpha\|_{\cb} < \lambda +\ep$ if (b) of (3.11) holds.  
\item"(ii)" $(T_\alpha)_{ |G\cap X} = I_{|G\cap X}$.  
\item"(iii)" $\|T_\alpha (u_{ij}) + (I-T_\alpha)(v_{ij})\| \le (1+\ep)
\max \{\|(u_{ij})\|,\|(v_{ij})\|\}$ for all $(u_{ij})$ and $(v_{ij})\in 
M_n(G)$.  
\item"(iv)" $T_\alpha$ is finite rank in case 3.11(a) holds, or 3.11(b) holds 
wtih $X$ $\lambda$-nuclear. 
\endroster
\endproclaim

\demo{Proof} 
Let $X^{**}\oplus W$ be the complete $M$-decomposition of $Y^{**}$ and 
$P$ the projection onto $X^{**}$ with kernel $W$. 
Let $G_1 = PG$ and $G_2 = (I-P)G$, and set $\tilde G= G_1\oplus G_2$. 
Then evidently 
$$G\subset \tilde G\ .
\tag 3.44$$
Now first assume 3.11(b). 
By Lemma~3.13 (applied to ``$G$''~$=\tilde G$) we may choose 
$T:\tilde G\to Y$ a linear operator so that 
$$\gather
\|T\|_{\cb} < 1+\frac{\ep}{\lambda}\ ,\tag 3.45\\
T_{|\tilde G\cap Y} = I_{|\tilde G\cap Y}\tag 3.46\\
\noalign{\text{and}}
T(G_1 ) \subset X\ .\tag 3.47
\endgather$$
Since $P$ is a complete contraction, we have that 
$$\|TP_{|G} \|_{\cb} < 1+\frac{\ep}{\lambda}\ .
\tag 3.48$$
Hence by the Definition of $\lambda$-finite injectivity, since $\lambda  
(1+\frac{\ep}{\lambda}) = \lambda+\ep$, we may choose a linear extension 
$T_\alpha :Y\to X$ of $P_{|G}$ satisfying (i$'$) of 3.14, which is moreover  
finite rank if $X$ is $\lambda$-nuclear, by Lemma~3.10. 

Suppose now that (3.11(a)) holds. 
Let $\ep' >0$, to be decided, and choose $A:Y\to Y$ a finite rank operator 
such that
$$A_{|G} = I_{|_G}
\tag 3.49$$
and
$$\|A\| < \lambda +\ep'\ .
\tag 3.50$$
Now set  $G' = A(G)$, $G'_1 = PG'$, $G_2 = (I-P)G'$, and 
$\tilde{G'} = G'_1 \oplus G'_2$. 
Of course $G'\supset G$ by (3.44). 
Then by our extension of the Ge-Hadwin result, namely Lemma~3.13, we may 
choose $T:\tilde{G'} \to Y$ a linear operator so that 
$$\|T\|_n < 1+\ep'
\tag 3.51$$
and again, (3.46) and (3.47) hold (in fact, they hold replacing 
$\tilde G$ by $\tilde{G'}$ and $G_1$ by $G'_1$). 
Now define $T_\alpha$ by 
$$T_\alpha = TPA\ .
\tag 3.52$$
Thus 
$$\|T_\alpha\| < (1+\ep')(\lambda +\ep') <\lambda +\ep 
\tag 3.53$$
if $\ep'$ is chosen with $\ep' \lambda  + \ep' + (\ep')^2 <\ep$.
It then follows that if $x\in G\cap X$, $Tx=x$, and so in case (b), 
$T_\alpha x = TPx=x$ while in case (a), $T_\alpha (x) = TPAx=x$; 
thus (ii) of 3.14 holds. 
Now note that if $v\in G$, then by (3.44) and (3.46) 
$$Tv=v\ .
\tag 3.54$$
Hence for any such $v$,  
$$(I-TP)(v) = (T-TP)(v)\ .
\tag 3.55$$
Finally, let $(u_{ij})$ and $(v_{ij})$ be elements of $M_n(G)$. 
Then 
$$\leqalignno{
&\| T_\alpha u_{ij}+(I-T_\alpha)v_{ij})\| &(3.56)\cr 
&\qquad = \|(TPu_{ij}) + ((T-TP)v_{ij})\| \cr
&\qquad\qquad \text{ (by 3.55) and the definition on $T_\alpha$)}\cr
&\qquad \le (1+\ep) \|Pu_{ij} + (I-P)v_{ij}\|\cr
&\qquad\qquad \text{ (by (3.45) in case (b), (3.51) in case (a))}\cr
&\qquad = (1+\ep) \max \{\|(Pu_{ij})\|,\|(I-P)v_{ij}\|\}\cr
&\qquad \le (1+\ep) \max \{\|(u_{ij})\|,\|(v_{ij})\|\}\ .}
$$
(The last equality holds because $X^{**}\oplus W$ is a complete 
$M$-decomposition of  $Y^{**}$, while the last inequality holds because 
$P$ and $I-P$ are complete contractions.) 
  
Thus 3.14iii holds, completing the proof of the lemma.\qed
\enddemo

We are finally prepared for the 

\demo{Proof of Theorem 3.11}

Let $\D$ be the directed set consisting of all $\alpha = (G,n,\ep)$ where 
$G$ is a finite-dimensional subspace of $Y$ with $G\cap X\ne\{0\}$, 
$n\ge1$, and $\ep>0$. 
Of course $(G,n,\ep)\le (G',n',\ep')$ if $G\subset G'$, $n\le n'$ and 
$\ep \ge \ep'$. 
For each such $\alpha$, choose $T_\alpha$ satisfying the conclusion 
of Lemma~3.14. 

Finally, let $\tilde T_\alpha= T_\alpha (1+\frac{\ep}{\lambda})^{-1}$. 
Then trivially 
$$\|\tilde T_\alpha\|_{\cb} \le (\lambda+\ep) \Big(1+\frac{\ep}{\lambda}
\Big)^{-1} = \lambda\text{ (by (3.14(i)).}
\tag 3.56$$
Since also trivially 
$$\lim_\alpha \|T_\alpha - \tilde T_\alpha\| =0\ ,
\tag 3.57$$
it suffices to prove that $(T_\alpha)_{\alpha\in\D}$ is an M-cai, for 
then (3.57) yields that $(\tilde T_\alpha)_{\alpha\in\D}$ is an M-cai. 
Now fix $n$, $x\in X$ with $x\ne0$, and $(u_{ij}),(v_{ij}) \in M_n(Y)$. 
Then letting $G = \text{sp}\{ x,u_{ij},v_{ij}: 1\le i,\, j\le n\}$, and 
$\ep>0$, $\alpha =  (G,n,\ep) \in \D$; hence for any $\beta\in\D$ with 
$\beta\ge\alpha$, $T_\beta x=x$ and (iii) of Lemma~3.14 holds. 
This completes the proof of the Theorem.\qed
\enddemo

We next briefly discuss the case when $\K(X)$ is an $M$-ideal in $B(X)$. 
Here, $X$ is a fixed Banach space, $\K(X)$ denotes the space of compact 
operators on $X$, and $B(X)$ denotes the space of compact operators on $X$. 
For a comprehensive survey of known facts, see section~VI.4 of \cite{HWW}. 
These results yield immediately (via our Theorem~1.1) that 
{\it $\K(X)$ is an $M$-ideal in $B(X)$ if and only if $\K(X)$ admits a 
(contractive) {\rm M-ai} in $B(X)$.}
(See Theorem~VI.4.17 of \cite{HWW}.) 
Let us note --- it is known that this is the case for $X=c_0$ or $\ell^p$, 
$1<p<\infty$ (generalizing the $p=2$ case).

It also follows (from known results) that this property is {\it not\/} 
hereditary; in fact, the final example of this paper yields a subspace $X $
of $c_0$ so that $\K(X)$ is not an $M$-ideal in $B(X)$ (via many 
known results, e.g., the result of C.-M.~Cho and W.B.~Johnson, 
that for $X\subset c_0$, $\K(X)$ is an $M$-ideal 
in $B(X)$ iff $X$ has the compact metric approximation property \cite{CJ}. 
(See also \cite{KW} for a remarkable generalization). 
In fact, it is known that if $\K(X)$ is an $M$-ideal in $B(X)$, then $X$ 
has a shrinking contractive approximation to the identity $(K_\alpha)$, 
consisting of compact operators, so that in fact letting $T_\alpha\in B(B(X))$
be defined as $T_\alpha (T) = K_\alpha T$ for all $\alpha$, then 
$(T_\alpha)$ is an M-cai for $\K(X)$ in $B(X)$. 
It is also known that $X$ is then an $M$-ideal in $X^{**}$. 
We now pose the following questions (which certain seem accessible via 
the technology given in \cite{HWW}). 
Does $X$ then admit an M-ai in $X^{**}$? 
In fact, can the shrinking compact approximation to the identity $(K_\alpha)$ 
be chosen as above, so that additionally $(K_\alpha^{**})$ is an M-ai for $X$ 
in $X^{**}$ (in the case $X$ is non-reflexive)? 
Note this question is simply equivalent to: 
can $(K_\alpha)$ be so chosen so that in addition, 
$$\varlimsup_\alpha \|K_\alpha^{**} x^{**} + (I- K_\alpha^{**}) y^{**}\| 
\le \max\{ \|x^{**}\| \, \|(y^{**})\|\}\text{ for all } x^{**},y^{**}\in 
X^{**}\ ?$$
(By Theorem 5.3.b of \cite{HWW}, this is indeed so if $(\K(X\oplus X)$ is an 
$M$-ideal in $\L(X\oplus X)$ ($L^\infty$-direct sum), for then $(K_\alpha)$ 
may be chosen so that $\varlimsup_\alpha \|S\circ (K_\alpha \oplus I-K_\alpha
)\| =1$, $S:X^{**} \oplus X^{**}\to X^{**}$ the ``sum'' operator.) 
Finally, we ask: can the $K_\alpha$'s be chosen so that letting  
$(T_\alpha)$ be as above, then $(T_\alpha)$ is a strong M-ai for $\K(X)$ 
in $B(X)$? 
Equivalently, so that $T_\alpha^* u^*\to u^*$ in norm for all $u^*\in 
B(X)^*$? 

We conclude with an example of a Banach space $X$ which is an $M$-ideal 
in $X^{**}$, but so that there does not exist either an M-ai or a 
weak contractive M-ai for $X$ in $X^{**}$. 
The example follows quickly from known (but rather non-trivial!) results; 
the same example (for a different purpose) is given in \cite{JO}. 

We first need a standard fact (also given in \cite{JO}). 

\proclaim{Lemma 3.15} 
Let $X$ be a closed linear subspace of $c_0$. 

{\rm (a)} Let $T:X\to X$ be a given (bounded linear) operator. 
Then either $T$ is compact or there is a subspace $Y$ of $X$ with $Y$ 
isomorphic to $c_0$ so that $T_{|Y}$ is an isomorphism. 

{\rm (b)} Let $T:X^{**}\to X$ be a given operator. 
Then $T{|_X}$ is compact. 
\endproclaim

\demo{Proof} 
(a) Any semi-normalized weakly null sequence in $c_0$ contains a 
subsequence equivalent to the usual $c_0$ basis. 
This implies $X^*$ has the Schur property (i.e., weak and norm 
sequential convergence coincide on $X^*$), whence $T$ weakly compact 
implies $T$ compact. 
But if $T$ is not weakly compact, there exists a bounded sequence $(u_j)$ 
in $(c_0)$ so that $(Tu_j)$ has no weakly convergent subsequence. 
We may then pass to a subsequence $(u_{n_j})$ of $(u_j)$ so that {\it both\/} 
$(u_{n_j})$ and $(Tu_{n_j})$ are equivalent to the summing basis of $c_0$, 
which implies that $Y\defeq [u_{n_j}]$ is isomorphic so $c_0$ and $T|Y$ 
is an isomorphism. 

(b) Suppose that $T|X$ were not compact. 
Then by part~(a), there exists a subspace $Y$ of $X$ so that $Y$ is 
isomorphic to $c_0$ and $T|Y$ is an isomorphism. 
Hence $Y^{**}$ is isomorphic to $\ell^\infty$ and 
$P\defeq (T_{|Y})^{-1}
T_{|Y^{**}}$ is a projection from $Y^{**}$ into $Y$, which contradicts the 
fact that $c_0$ is uncomplemented in $\ell^\infty$.\qed
\enddemo

Finally, we recall a rather deep result of T.~Szankowski 
(cf. \cite{LT, Theorem~2.a.7}) 
\vskip1pt
\noindent {(3.58)}\qquad
{\it there exists a subspace $X$ of $c_0$ failing the compact}

\noindent {\hphantom{(3.58)\qquad}} bounded approximation property.

\proclaim{Proposition 3.16} 
Let $X$ satisfy $(3.58)$. 
Then $X$ is an $M$-ideal in $X^{**}$, but $X$ has no M-ai or weak 
contractive M-ai in $X^{**}$. 
\endproclaim

\demo{Proof} 
Since $c_0$ is an $M$-ideal in $\ell^\infty = c_0^{**}$, $X$ is an $M$-ideal 
in $X^{**}$ (cf. \cite{HWW, page 111}). 
If $(T_\alpha)_{\alpha\in\D}$ were either an M-ai or a weak contractive 
M-ai for $X$ in $X^{**}$, we would have that 
$$\gather
(T_\alpha) \text{ is uniformly bounded}\tag 3.59\\
T_\alpha :X^{**} \subset X\text{ for all } \alpha \tag 3.60\\
\noalign{\text{and}}
T_\alpha x\to x\text{ for all } x\in X\ .\tag 3.61
\endgather$$
But then by Lemma 3.15b, $(T_\alpha)_{|X}$ is compact for all $\alpha$, 
hence by (3.59) and (3.61), $X$ has the compact bounded approximation 
property, a contradiction. 
\enddemo

We conjecture that one may also find  a separable situation in which there 
are no M-ai's for a certain $M$-ideal. 
Precisely,

\proclaim{Conjecture} 
Let $X$ satisfy $(3.58)$. 
Then there exists a separable $Y$ with $X\subset Y\subset X^{**}$ so that 
there is no uniformly bounded sequence $(T_n)$ of operators satisfying 
\roster
\item"(i)" $T_nY\subset X$ for all $n$; 
\item"(ii)" $T_nx \to x$ for all $x\in X$;  
\item"(iii)" $\varlimsup_n \|T_nx + (I-T_n)y\|\le \max\{\|x\|,\|y\|\}$ 
for all $x\in X$ and $y\in Y$.
\endroster
\endproclaim

Of course if $X$ satisfies (3.58) and $X\subset Y\subset X^{**}$, $X$ is 
an $M$-ideal in $Y$; if $X$ satisfied the Conjecture, $X$ would admit no 
M-ai or weak contractive M-ai in $Y$. 
The Conjecture, however, appears to lie much further below the surface 
(modulo known results) than Proposition~3.15.

\head Appendix. An isometric lifting theorem\endhead

We obtain here an  operator-space 
generalization of the Effros-Haagerup lifting result. 

\proclaim{Theorem} 
Let $X\subset Y$ be operator spaces with $X$ nuclear, $Y$ 1-locally reflexive,
and $Y/X$ separable. 
Assume that $X$ is a complete $M$-ideal in $Y$. 
Then there exists a completely contractive lift $L:Y/X\to Y$ of $I_{Y/X}$. 
\endproclaim

In the classical case (i.e., MIN operator structures), a nuclear operator 
space $X$ is an $L^1$-predual, and since all Banach spaces are 1-locally 
reflexive, the Theorem reduces to Ando's result that one always has 
contractive liftings of the identity on $Y/X$ when $Y/X$ is separable and $X$ 
is an $L^1$-predual which is an ideal in $Y$, \cite{An}. 
In fact, our proof of the Theorem 
is the quantized version of Ando's argument, as expressed 
in \cite{HWW, page 58}. 

We first assemble the facts needed to obtain the isometric assertions 
of the Theorem. 
Throughout, we assume that $X\subset Y$ are Banach spaces with $X$ an 
$M$-ideal in $Y$; $\pi :Y\to Y/X$ denotes the quotient map. 

\proclaim{Lemma A1} 
Given $e_0\in Y/X$, there exists a $y_0\in Y$ with $\|e_0\| = \|y_0\|$ 
and $\pi y_0 = e_0$.
\endproclaim

\demo{Proof} 
Let $y\in Y$ be such that $\pi y = e_0$. 
Then $d(y,X) = \|e_0\|$. 
By Proposition~II.1.1 of \cite{HWW} (the proximality of $M$-ideals) there 
exists an $x_0\in X$ with $\|y-x_0\| = d(y,X)$.
Then $y_0\defeq y-x_0$ is the desired element of $Y$. 
\enddemo

\proclaim{Lemma A2} 
Let $V$ be a closed linear subspace of $X$ and $L\in Y$. 
Suppose $d(L,V) \defeq \inf \{\|L-v\| : v\in V\} =1$. 
Then for all $\ep>0$, there exists $V_\ep\in V$ and 
$L_\ep \in Ba(Y)$ so that 
$$\|(L-V_\ep) - L_\ep \| \le \ep\ \text{ and }\ 
(L-V_\ep) -L_\ep \in X\ .
\tag A1$$
\endproclaim

\demo{Proof} 
This follows from Lemma 2.5 of \cite{HWW}. 
We sketch a proof for completeness. 
Since $X$ is an $M$-ideal, besides its proximality, $X$ also has 
the ``strict 2-ball'' property: 
{\it given $B_1,B_2$ closed balls in $Y$ with $\text{\rm Int}(B_1\cap B_2) 
\ne\emptyset$ and $B_i\cap X\ne\emptyset$ for $i=1,2$, then } 
$B_1\cap B_2 \cap X\ne \emptyset$. 

First choose $V_\ep \in V$ with 
$$\|L'\| < 1+\ep \ ,\ \text{ where }\ L' = L-V_\ep\ .
\tag A2$$
Now let $B_1 = B(L', 1)$. Of course $d(L',V) = d(L,V)=1$, so by the 
proximality of $X$, $B_1\cap X\ne\emptyset$. 
Now set $B_2 = B(0,\ep)$ ($=\ep BaY$).
Since $\|L'\|\ge1$, (A2) yields that 
Int$(B_1\cap B_2)\ne\emptyset$. 
Hence choosing $x\in B_1\cap B_2\cap X$, and letting $L_\ep = L'-x$, we have 
that 
$$\|L_\ep\|\le 1\ \text{ (since $x\in B(L',1)$)}
\tag A3$$
and of course 
$$(L-V_\ep )-L_\ep = L'-L_\ep = x\in X\ .\eqno\qed$$
\enddemo

A simple compactness argument yields the next result. 

\proclaim{Lemma A3} 
If $X$ is a nuclear operator space, $X^{**}$ is an isometrically injective 
operator space.
\endproclaim 

\re{Remark}
In fact, it is proved in \cite{EOR} that $X$ is nuclear iff $X$ is locally 
reflexive and $X^{**}$ is 1-injective. 
\endre

\proclaim{Lemma A4} 
Let $E$ be a finite-dimensional operator space, $Y$ an operator space, and 
let $\cb(E,Y)$ be the operator space of completely bounded maps from $E$ to 
$Y$. Then if $Y$ is 1-locally reflexive, $\cb(E,Y)^{**}$ is (canonically 
isometric to) $\cb(E,Y^{**})$. 
\endproclaim

\demo{Proof sketch} 
Of course $\cb(E,Y)$ is nothing but the space of linear maps $T$ from $E$ 
to $Y$, endowed with $\|T\|_{\cb}$. 
Thus $\cb(E,Y) = E^* \otimes_{\op} Y$. 
But $Y$ is 1-locally reflexive iff $F\otimes_{\op} Y^{**} = (F\otimes_{\op}
Y)^{**}$ isometrically for all finite-dimensional operator spaces 
$F$ (cf. \cite{EH}).\qed
\enddemo 

The next result is again a quantization of an observation in \cite{HWW}
(see page~62).

\proclaim{Lemma A5} 
Assume that $X$ and $Y$ are operator spaces with $Y$ 1-locally reflexive 
and $X$ a complete $M$-ideal in $Y$, and let $E$ be a finite-dimensional
operator space. 
Then $\cb(E,X)$ is a complete $M$-ideal in $\cb(E,Y)$.
\endproclaim

\demo{Proof} 
Let $W$ be the (weak*-closed) subspace of $Y^{**}$ so that 
$Y^{**} = X^{**}\oplus W$ is a complete $L^\infty$-decomposition for $Y^{**}$, 
and set $F= E^*$. 
Then 
$$(F\otimes_{\op} X^{**})\oplus (F\otimes_{\op} W)\text{ is a complete 
$L^\infty$-decomposition for } F\otimes_{\op} Y^{**}.$$
The result now follows upon identifying $F\otimes_{\op} X^{**}$ with 
$\cb(E,X)^{**}$ and $F\otimes_{\op} Y^{**}$ with $\cb(E,Y)^{**}$, via 
Lemma~A4.\qed
\enddemo

At last, we are prepared for the fundamental lemma yielding the proof of the 
Theorem. 

\proclaim{Lemma A6 (The Crucial Lemma)} 
Let $X\subset Y$ be operator spaces with $X$ nuclear and $Y$ 1-locally 
reflexive. 
Let $E_1\subset E$ be finite-dimensional subspaces of $Y/X$, and let 
$L_1:E_1\to Y$ be a completely contractive lift of $I_{E_1}$. 
Then given $\ep>0$, there exists a completely contractive lift 
$L_2 :E \to Y$ of $I_E$ with $\|L_2|_{E_1} -L_1\|\le\ep$. 
\endproclaim

A6 is a simple consequence of the results already assembled and the basic 

\proclaim{Sublemma}
Assuming the hypotheses of A6, there exists a lift $L_{\ep} :E\to Y$ of $I_E$ 
with $L_{\ep}|_{E_1} = L_1$ and $\|L_{\ep}\|_{\cb} <1+\ep$.
\endproclaim

\demo{Proof}
Let $P$ be the $M$-projection from $Y^{**}$ onto $X^{**}$ and $L:Y^{**}/X^{**}
\to Y^{**}$ the completely contractive lift of the identity on $Y^{**}/X^{**}$ 
induced by $P$. 
Since $L|_{E_1}$ and $L_1$ are then both lifts of $I_{E_1}$ into $Y^{**}$, 
$L(e) - L_1(e)\in X^{**}$ for all $e\in E_1$, hence 
$$(I-\Gamma) (L|_{E_1} - L_1) = 0\ .
\tag A4$$
Since $X^{**}$ is 1-injective by Lemma~A3, we may choose $\theta :E_2\to
X^{**}$ a complete contraction with 
$$\theta |_{E_1} = PL_1\ .
\tag A5$$
Next, define $\tilde L:E_2\to Y^{**}$ by 
$$\tilde L= \theta + (I-P)L|_{E_2}\ .
\tag A6$$ 
Then if $e\in E_2$, 
$$\pi^{**} \tilde L(e) = \pi^{**} (I-P) L(e) = \pi^{**} L(e)=e\ .
\tag A7$$
Thus, $\tilde L$ is a lift of $I_{E_2}$ into $Y^{**}$. 
If $e\in E_1$, then 
$$\align 
\tilde L(e) & = PL_1 (e) + (I-P) L(e)\tag A8\\
& = PL_1 (e) + (I-P) L_1(e)\ \text{ (by (A4)}\\
& = L_1 (e)\ .
\endalign$$
Since $P$ is an $M$ projection, $X^{**} \oplus (I-P)Y^{**}$ is a complete 
$M$-decomposi-tion of $Y^{**}$, where by (A6)--(A8), $\tilde L$ is a completely
contractive lift of $I_{E_2}$ extending $L_1$. 

Of course $\tilde L$ lifts into $Y^{**}$, not $Y$. 
To get a lift into $Y$, we apply our extended local reflexivity principal 
for operator spaces, Lemma~3.12. 
First let $\chix :E\to Y$ be a linear lift of $I_E$. 
Now let $G= \tilde L(E) + \chix (E)$. 
Choose $T:G\to Y$ satisfing (i), (iii) and (iv) of Lemma~3.12 (for 
$\lambda =1$). 
Define $L_\ep :E\to Y$ by 
$$L_\ep = T\tilde L\ .
\tag A9$$
Then if $e\in E_1$, $\tilde L(e) = L_1 (e) = L_\ep(e)$, where the last 
equality holds by  (3.12iii) since $L_1(e) \in Y$. 

$\|L_\ep\|_{\cb} < 1+\ep$ by (A9) and Lemma 3.12(i), since $\tilde L$ 
is a complete contraction. 
Finally, if $e\in E$, 
$$\align
T\tilde L(e) & = T(\tilde L(e) - \chix (e)) + T\chix (e)\tag A10\\
& = T(\tilde L(e) - \chix (e)) + \chix (e)
\endalign$$
since $T_{|G\cap Y} = I|_{G\cap Y}$.
But then 
$$\pi T\tilde L(e) = \pi \chi (e) =e\ ,
\tag A11$$
since $\tilde L(e) - \chix (e) \in X^{**}\implies T(\tilde L(e) 
- \chix (e)) \in X$ by 3.12(iv). 
Thus $L_\ep$ satisfies the conclusion of the Sublemma.\qed
\enddemo

\demo{Proof of Lemma A6} 
Let $\tilde X = \cb (E,X)$ and $\tilde Y = \cb (E,Y)$, and let 
$V = \{ T\in \tilde X :\ker T\supset E_1\}$. 
Let $L$ be a lift of $I_E$ so that $L|_{E_1} = L_1$, and let $\ep>0$. 
Let $L_\ep$ be a lift of $I_E$ satisfying the conclusion of the sublemma. 
Then $L-L_\ep \in V$. 
This proves that 
$$d(L,V) =1\ .
\tag A12$$
We now apply Lemma~A2 to $\tilde X\subset \tilde Y$;
$\tilde X$ is an $M$-ideal in $\tilde Y$ by Lemma~A5. 
Thus, choose $V_\ep\in V$ and $L_\ep \in Ba(\tilde Y)$ satisfying (A1), 
and set $L_2 = L_\ep$. 
Now it is trivial that $L-V_\ep$ is a lift of $I_E$ and moreover 
$(L-V_\ep)|_{E_1}=L_1$. 
Since $(L-V_\ep) - L_2 \in\tilde X$, 
$[(L-V_\ep)-L_2] (E) \subset X$, whence $\pi (L-V_\ep)=\pi L_2 = I_E$, i.e., 
$L_2$ is indeed a lift of $I_E$, which is completely contractive. 
Finally, 
$$\|L_1 -L_2|_{E_1}\| 
= [(L-V_\ep)-L_\ep)|_{E_1}\| 
\le \|(L-V_\ep)-L_\ep\| \le \ep\ \text{ by (A1);} 
\tag A13$$
thus $L_2$ satisfies the conclusion of the lemma.\qed 
\enddemo

\demo{Proof of the Theorem} 
Let $e_0\in E_0$ with $\|e_0\|=1$. 
Choose $y_0\in Y$ with $\pi y_0  = e_0$ and $\|y_0\| = \|e_0\|$ (by 
Lemma~A1) and set $E_0 = [e_0]$. 
Choose finite dimensional spaces $E_0\subset E_1\subset E_2 \subset \cdots$ 
of $Y/X$ with $\bigcup_{j=1}^\infty E_j$ dense in $Y/X$. 
Define $L_0 :E_0\to Y$ by $L_0 (\lambda e_0) = \lambda y_0$ for all 
scalars $\lambda$. 
It's trivial that $L_0$ is a completely contractive lift of $I_{E_0}$. 
Now let $\ep>0$ and suppose $i\ge0$ and $L_i :E_i\to Y$ a completely 
contractive lift of $I_{E_i}$ has been chosen. 
By the Crucial Lemma, we may choose $L_{i+1} :E_{i+1} \to Y$ a completely 
contractive lift of $I_{E_{i+1}}$ with 
$$\|L_{i+1}^\ep |_{E_i} - L_i^\ep\| \le \frac{\ep}{2^{i+1}} \ .
\tag A14$$ 
Then it follows (as in the proof of Theorem~2.4) that setting 
$Z= \bigcup_{i=0}^\infty E_i$, then $\lim_{i\to\infty} L_i (z) 
\defeq L(z)$ exists for all $z\in Z$, and $L$ extends to a completely 
contractive lift of $I_{Y/X}$. 
This completes the proof.\qed
\enddemo

\remark{Remark}
We also have in the above argument that 
$$\|L (e_0) -y_0\| \le \sum_{i=1}^\infty \|L_{i+1}^\ep - L_i^\ep (y_0)\| 
\le \sum_{i=0}^\infty \frac{\ep}{2^{i+1}} = \ep\ .
\tag A15$$
That is, we have also proved that 
{\it given $\ep>0$, $e_0\in Y/X$ and $y_0\in Y$ with $\|e_0\| = 1= \|y_0\|$ 
and $\pi y_0 = e_0$, then the lift $L$ may also be chosen so that 
$\|Le_0 - y_0\|\le\ep$.} 
At this level of generality, however, it is impossible to insure that $L$ 
may be chosen with $Le_0=y_0$. 
Indeed, one may give an extreme point obstruction by constructing $X$ and $Y$ 
satisfying the hypotheses of the Theorem with $e_0$, $y_0$ as above, $y_0$ 
an extreme point of $Ba\, Y$, but so that 
$e_0$ is {\it not\/} an extreme point of $Ba \, Y/X$. 
Since $L$ is an isometry, we cannot then have that $Le_0$ is an extreme 
point of $Ba\, Y$.
\endremark


\Refs\widestnumber\key{HWW} 

\ref\key AP
\by C.A. Akeman and G.K. Pedersen
\paper Ideal perturbations of elements in $C^*$-algebras
\jour Math. Scand. \vol 41 \yr 1977 \pages 117--139
\endref

\ref\key AE
\by E. Alfsen and E. Effros
\paper Structure in real Banach spaces 
\jour  Ann. of Math. \vol 96 \yr 1972 \pages 98--173
\endref

\ref\key A
\by T.B. Andersen
\paper Linear extensions, projections, and split faces 
\jour J. Funct. Anal. \vol 17 \yr 1974 \pages 161--173 
\endref

\ref\key An
\by T. Ando
\paper A theorem on non-empty intersection of convex sets and its 
applications 
\jour J. Approx. Th. \vol 13 \yr 1975 \page 158--166
\endref

\ref\key Ar
\by W. Arveson
\paper Notes on extensions of $C^*$-algebras
\jour Duke Math. J. \vol 44 \yr 1977 \pages 329--355
\endref 

\ref\key B
\by  K. Borsuk
\paper \"Uber Isomorphie der Funktional r\"aume
\jour Bull. Int. Acad. Pol. Sci. \yr 1933 \pages 1--10
\endref 

\ref\key Be
\by S.F. Bellenot
\paper Local reflexivity of normed spaces
\jour J. Funct. Anal \vol 59 \yr 1984 \page 1--11
\endref

\ref\key CM
\by C.-M. Cho and W.B. Johnson
\paper $M$-ideals and ideals in ${\Cal L}(X)$
\jour J. Operator Theory \vol 16 \yr 1986 \pages 245--260
\endref 

\ref\key CE
\by M.-D. Choi and E. Effros
\paper The completely positive lifting problem for $C^*$-algebras
\jour Ann. of Math. \vol 104 \yr 1976 \pages 585--609
\endref

\ref\key D
\by K.R. Davidson
\book $C^*$-Algebras by Example 
\publ Amer. Math. Soc. 
\publaddr Providence \yr 1996 
\endref

\ref\key Du
\by J. Dugundji
\paper An extension of Tietze's theorem
\jour Pacific J. Math. \vol 1 \yr 1951\pages 353--367
\endref

\ref\key EH
\by E. Effros and U. Haagerup 
\paper Lifting problems and local reflexivity for $C^*$-algebras 
\jour Duke Math. J. \vol 52 \yr 1985 \pages 103--128
\endref

\ref\key EOR
\by E. Effros, N. Ozawa and Z.-J. Ruan
\paper On injectivity and nuclearity for operator spaces 
\toappear
\endref

\ref\key ER1
\by E. Effros and Z.-J. Ruan
\paper On non-self-adjoint operator algebras
\jour Proc. Amer. Math. Soc. 
\vol 110 \yr 1990 \pages 915--922
\endref

\ref\key ER2
\by E. Effros and Z.-J. Ruan
\paper Mapping spaces and liftings for operator spaces 
\jour Proc. London Math. Soc. \vol 69 \yr 1994 \pages 171--197
\endref

\ref\key GH 
\by L. Ge and D. Hadwin
\paper Ultraproducts for $C^*$-algebras
\toappear
\endref

\ref\key GN
\by I.M. Gelfand and M.A. Neumark 
\paper On the imbedding of normed rings into the ring of operators in 
Hilbert space
\jour Mat. Sbornik N.S. \vol 12\issue 54 \yr 1943 \pages 197--213
\endref 

\ref\key HWW
\by P. Harmand, D. Werner and W. Werner
\paper $M$-ideals in Banach spaces and Banach algebras 
\inbook SLNM 1547 
\publ Springer-Verlag
\yr 1993
\endref

\ref\key JO
\by W.B. Johnson and T. Oikhberg
\paper Separable lifting property and extensions of local reflexivity 
\toappear
\endref

\ref\key JRZ
\by W.B. Johnson, H.P. Rosenthal and M. Zippin
\paper On bases, finite dimensional decompositions and weaker structures 
in Banach spaces 
\jour Israel J. Math. \vol 9 \yr 1977 \pages 488--506
\endref 

\ref\key K 
\by N.J. Kalton
\paper $M$-ideals of compact operators
\jour Illinois J. Math. \vol 37 \yr 1993 \pages 147--169
\endref

\ref\key KW
\by N.J. Kalton and D. Werner 
\paper Property $(M)$, $M$-ideals, and almost isometric structure of Banach 
spaces 
\jour J. f\"ur die reine und angewandte Mathematik \vol 461 \yr 1995 
\pages 137--178
\endref 

\ref\key Ki
\by E. Kirchberg  
\paper  On non-semisplit extensions, tensor products and exactness of 
group $C^*$-algebras 
\jour Invent. Math. \vol 112 \yr 1993 \pages 449--489
\endref

\ref\key KR
\by S-H. Kye and Z.-J. Ruan
\paper On local lifting property for operator spaces 
\finalinfo preprint
\endref 

\ref\key LR
\by J. Lindenstrauss and H.P. Rosenthal
\paper Automorphisms in $c_0$, $\ell_1$ and $m$ 
\jour Israel J. Math. \vol 7 \yr 1969 \pages 222--239 
\endref

\ref\key LT
\by J. Lindenstrauss and L. Tzafriri
\book Classical Banach Spaces I
\publ Springer-Verlag 
\yr 1977
\endref

\ref\key M 
\by E. Michael 
\paper Some extension theorems for continuous functions
\jour Pacific J. Math. \vol 3 \yr 1953 \pages 789--806
\endref 

\ref\key OR
\by T. Oikhberg and H.P. Rosenthal
\paper On certain extension properties for the space of compact operators 
\jour J. Funct. Anal.
\finalinfo submitted
\endref

\ref\key O 
\by N. Ozawa 
\paper A short proof of the Oikhberg-Rosenthal Theorem
\finalinfo preprint
\endref

\ref\key Pe
\by A. Pe{\l}czy\'nski 
\paper Projections in certain Banach spaces 
\jour Studia Math. \vol 29 \yr 1960 \pages 209--227 
\endref

\ref\key Pi 
\by G. Pisier
\paper An introduction to the theory of operator spaces 
\finalinfo preprint
\endref

\ref\key R
\by A. Guyan Robertson 
\paper Injective matricial Hilbert spaces 
\jour Math. Proc. Camb. Phil. Soc. \vol 110 \yr 1991 \pages 183--190
\endref

\ref\key Ro
\by H.P. Rosenthal 
\paper The complete separable extension property 
\jour J. Operator Theory 
\toappear
\endref 

\ref\key S
\by R. Smith 
\paper An addendum to ``$M$-ideal structures in Banach algebras'' 
\jour J. Funct. Anal. \vol 32 \yr 1979 \pages 269--271
\endref

\ref\key SW
\by R. Smith and J. Ward 
\paper $M$-ideal structure in Banach algebras 
\jour J. Funct. Anal \vol 27 \yr 1978 \pages 337--349
\endref

\ref\key So
\by A. Sobczyk
\paper Projection of the space $(m)$ on its subspace $(c_0)$ 
\jour Bull. Amer. Math. Soc. \vol 47 \yr 1941 \pages 938--947
\finalinfo MR3-205
\endref

\ref\key V
\by W.A. Veech
\paper Short proof of Sobcyzk's theorem
\jour Proc. Amer. Math. Soc. \vol 28 \yr 1971 \pages 627--628
\endref 

\ref\key Z
\by M. Zippin 
\paper The separable extension problem 
\jour Israel J. Math. \vol 26 \yr 1977 \pages 372--387 
\endref

\endRefs
\enddocument